\pdfoutput=1
\documentclass[11pt,letterpaper]{amsart}
\usepackage[colorlinks=true,
pdfstartview=FitV, linkcolor=cyan, citecolor=magenta,
urlcolor=blue]{hyperref}
\usepackage{amsthm,amsmath,amsfonts,latexsym,amssymb}
\usepackage[normalem]{ulem}
\usepackage{hyperref}
\usepackage{enumitem}
\usepackage{mathrsfs}
\usepackage[latin1]{inputenc}
\usepackage[T1]{fontenc}
\usepackage{ae,aecompl}
\usepackage{braket}
\usepackage{comment}
\usepackage{color}
\usepackage{graphicx}
\usepackage{subfig}
\theoremstyle{plain}
\newtheorem{theorem}{Theorem}[section]
\newtheorem{lemma}[theorem]{Lemma}
\newtheorem{proposition}[theorem]{Proposition}
\newtheorem{corollary}[theorem]{Corollary}

\newtheorem{observation}[theorem]{Observation}
\newtheorem*{thm*}{\protect\theoremname}
\newtheorem*{cor*}{\protect\corollaryname}
\theoremstyle{definition}
\newtheorem{definition}[theorem]{Definition}
\newtheorem{remark}[theorem]{Remark}
\newtheorem*{remark*}{Remark}

\long\def\@savemarbox#1#2{\global\setbox#1\vtop{\hsize\marginparwidth 
  \@parboxrestore\tiny\raggedright #2}}
\renewcommand{\d}{{\rm d}}

\newcommand{\Bc}{\mathcal B}
\newcommand{\Ec}{\mathcal E}
\newcommand{\Fc}{\mathcal F}

\newcommand{\Kc}{\mathcal K}
\newcommand{\Sc}{\mathcal S}

\newcommand{\Cb}{\mathbb C}
\newcommand{\Ab}{\mathbb A}

\newcommand{\Zb}{\mathbb Z}

\newcommand{\id}{\mathrm{id}}
\newcommand{\Span}{\mathrm{Span}}
\newcommand{\PU}{\mathsf{PU}}
\newcommand{\PO}{\mathsf{PO}}
\newcommand{\PSL}{\mathsf{PSL}}
\newcommand{\PGL}{\mathsf{PGL}}
\newcommand{\GL}{\mathsf{GL}}

\newcommand{\abs}[1]{\left|#1\right|}
\newcommand{\ip}[1]{\left\langle #1\right\rangle}
\newcommand{\norm}[1]{\left\|#1\right\|}

\DeclareMathOperator{\Bb}{\mathbb{B}}
\DeclareMathOperator{\Hb}{\mathbb{H}}

\DeclareMathOperator{\SL}{\mathsf{SL}}
\DeclareMathOperator{\Kb}{\mathbb{K}}
\DeclareMathOperator{\Gr}{\mathrm{Gr}}
\DeclareMathOperator{\Ac}{\mathcal{A}}
\DeclareMathOperator{\Pb}{\mathbb{P}}
\DeclareMathOperator{\Rb}{\mathbb{R}}

\DeclareMathOperator{\Nb}{\mathbb{N}}
\DeclareMathOperator{\Aut}{Aut}

\providecommand{\corollaryname}{Corollary}
\providecommand{\theoremname}{Theorem}

\begin{document}

\title{Entropy rigidity for cusped Hitchin representations}
\author[Canary]{Richard Canary}
\address{University of Michigan}
\author[Zhang]{Tengren Zhang}
\address{National University of Singapore}
\author[Zimmer]{Andrew Zimmer}
\address{University of Wisconsin-Madison}
\thanks{Canary was partially supported by grant  DMS-1906441 from the National Science Foundation
and grant 674990 from the Simons Foundation. Zhang was partially supported by the NUS-MOE grants R-146-000-270-133 and A-8000458-00-00. Zimmer was partially supported by grants DMS-2105580 and DMS-2104381 from the
National Science Foundation.}

\begin{abstract}
We establish an entropy rigidity theorem for Hitchin representations of geometrically finite Fuchsian groups which generalizes a theorem of Potrie
and Sambarino for Hitchin representations of closed surface groups. In the process, we introduce the class of $(1,1,2)$-hypertransverse groups and show for such a group that the Hausdorff dimension of its conical limit
set agrees with its (first) simple root entropy, providing a common generalization of results of Bishop and Jones, for Kleinian groups, and Pozzetti, Sambarino and Wienhard, for Anosov groups. We also introduce the theory of transverse representations of projectively visible groups as a tool for studying discrete subgroups of linear groups which are not necessarily Anosov or relatively Anosov.
\end{abstract}
\maketitle

\setcounter{tocdepth}{1}
\tableofcontents

\section{Introduction}

Hitchin \cite{hitchin} discovered a component of the space of  (conjugacy classes in $\PGL(d,\Rb)$ of) representations of a closed surface group $\pi_1(S)$ into $\mathsf{PSL}(d,\mathbb R)$ which is topologically a cell. Labourie \cite{labourie-invent} showed that the representations in this component, now known as Hitchin representations, are discrete
and faithful and even quasi-isometric embeddings, so share many properties with classical Fuchsian surface groups. Fock and Goncharov \cite{fock-goncharov}
showed that these representations are exactly the representations which admit equivariant continuous positive maps of the Gromov boundary of $\pi_1(S)$,
which one may identify with the limit set of any Fuchsian uniformization, into the space $\mathcal F$ of complete $d$-dimensional flags. One may
then naturally define a representation $\rho$  of a Fuchsian group $\Gamma\subset\mathsf{PSL}(2,\mathbb R)$ into $\mathsf{PSL}(d,\mathbb R)$ to be Hitchin
if there is a continuous positive $\rho$-equivariant map of the limit set of $\Gamma$ into $\mathcal F$.
Hitchin representations  of geometrically finite Fuchsian groups share many of the same properties as classical Hitchin representations
 (see Canary--Zhang--Zimmer \cite{CZZ}).

The main result of this paper is an entropy rigidity theorem for Hitchin representations of geometrically finite Fuchsian groups which generalizes
a result of Potrie--Sambarino  \cite{potrie-sambarino} from the classical setting. 
Let $\left(\mathfrak{a}^*\right)^+$ denote the set of linear functionals which
are strictly positive on the interior of the standard positive Weyl chamber for $\mathsf{PGL}(d,\mathbb R)$.
Each $\phi\in\left(\mathfrak{a}^*\right)^+$ can be written as a non-negative linear combination of the standard simple roots $\{\alpha_i\}_{i=1}^{d-1}$
which define the standard positive Weyl chamber.
If $\phi\in\left(\mathfrak{a}^*\right)^+$, we obtain a length function
on $\rho(\Gamma)$ given by $\ell^\phi(\rho(\gamma))=\phi(\nu(\rho(\gamma)))$ where $\nu(\rho(\gamma))$ is the Jordan projection of
$\rho(\gamma)$ (i.e. the logarithms of the moduli of generalized eigenvalues of $\rho(\gamma)$ in descending order). The $\phi$-entropy $h^\phi(\rho)$ of a
Hitchin representation is then the exponential growth rate of the number of conjugacy classes of hyperbolic elements whose images have $\phi$-length
at most $T$.

Given a subgroup $\Gamma\subset\PSL(d,\Rb)$, the \emph{Zariski closure} of $\Gamma$ in $\PSL(d,\Rb)$ is the intersection of its Zariski closure in $\PGL(d,\Rb)$ with $\PSL(d,\Rb)$. We recall that Sambarino \cite{sambarino-positive} showed that the Zariski closure in $\PSL(d,\Rb)$ of the image of a Hitchin representation is either
all of $\mathsf{PSL}(d,\mathbb R)$, an irreducible image
of $\mathsf{PSL}(2,\mathbb R)$ within $\mathsf{PSL}(d,\mathbb R)$, or conjugate to either $\mathsf{PSO}(d,d-1)\subset\mathsf{PSL}(2d-1,\mathbb R)$, 
$\mathsf{PSp}(2d,\mathbb R)\subset\mathsf{PSL}(2d,\mathbb R)$, or the copy of $\mathsf{G}_2$ in $\mathsf{PSL}(7,\mathbb R)$.

\begin{theorem}[see Theorem~\ref{thm:main in body}]\label{thm:main} 
If $\Gamma \subset \PSL(2,\Rb)$ is geometrically finite, $\rho : \Gamma \rightarrow \PSL(d, \Rb)$ is Hitchin and
$\phi = \sum c_j \alpha_j \in \left(\mathfrak{a}^*\right)^+$, then 
\begin{align*}
h^{\phi}(\rho) \le \frac{1}{c_1+\dots + c_{d-1}}.
\end{align*}
Moreover, equality occurs if and only if $\Gamma$ is a lattice and either 
\begin{enumerate}
\item
$\phi=c_k\alpha_k$ for some $k$.
\item
$\rho(\Gamma)$ lies in an irreducible image of $\mathsf{PSL}(2,\mathbb R)$.
\item
$d=2n-1$, the Zariski closure of $\rho(\Gamma)$ is conjugate to $\mathsf{PSO}(n,n-1)$ and  $\phi=c_k\alpha_k+c_{d-k}\alpha_{d-k}$ for some $k$.
\item
$d=2n$, the Zariski closure of $\rho(\Gamma)$ is conjugate to $\mathsf{PSp}(2n,\mathbb R)$ and $\phi=c_k\alpha_k+c_{d-k}\alpha_{d-k}$ for some $k$.
\item
$d=7$, the Zariski closure of $\rho(\Gamma)$ is conjugate to $\mathsf{G}_2$ and $\phi=c_1\alpha_1+c_3\alpha_3+c_4\alpha_4+c_6\alpha_6$ or $\phi=c_2\alpha_2+c_5\alpha_5$.
\end{enumerate}
\end{theorem} 

In the process of establishing our main result, we introduce the class of transverse subgroups of $\mathsf{PGL}(d,\mathbb R)$ which
includes all Anosov subgroups, all images of Hitchin representations,
all images of cusped Anosov representations of geometrically finite Fuchsian groups in the sense of \cite{CZZ}, all images of relatively dominated representations in the sense of Zhu \cite{feng},
all images of relatively Anosov representations in the sense of Kapovich--Leeb \cite{kapovich-leeb}, all subgroups of these groups, 
all discrete subgroups of $\mathsf{PO}(d-1,1)$, and all discrete
groups of projective automorphisms that preserve strictly convex domains with $C^1$ boundary in real projective space. We note that the definition of transverse groups and many of general results we establish about them do not assume finite generation.

We obtain upper bounds on the Hausdorff dimensions of
conical limit sets of $P_k$-transverse groups, generalizing results of Glorieux--Montclair--Tholozan \cite{GMT} and Pozzetti--Sambarino--Wienhard \cite{PSW1} from
the Anosov setting.

We further introduce the class of $(1,1,q)$-hypertransverse subgroups which include images of the $(1,1,q)$-hyperconvex Anosov representations introduced by
Pozzetti--Sambarino--Wienhard \cite{PSW1}, their subgroups, Hitchin representations of Fuchsian groups, and all discrete subgroups of $\mathsf{PO}(d-1,1)$. We show that for such subgroups (with $\sigma_2(\gamma)=\sigma_q(\gamma)$ for all
$\gamma\in\Gamma$), the Hausdorff dimension
of the conical limit set agrees with the first simple root critical exponent. This result is a common generalization
of results of Pozzetti--Sambarino--Wienhard \cite{PSW1}, for Anosov groups, and Bishop--Jones \cite{bishop-jones}, for Kleinian groups, and the proof makes use
of techniques drawn from each source. We observe that the $\phi$-entropy and the $\phi$-critical exponent of a geometrically finite Hitchin representation
agree. We conclude that the simple root entropies of Hitchin representations of geometrically finite Fuchsian groups are at most 1, and are exactly 1 only for Hitchin
representations of lattices. We combine this with convexity properties of the entropy functional to establish our main theorem.

\medskip

We now give a few definitions which allow us to give a more detailed discussion of our work. We recall that the standard {\em Cartan subspace} for
$\PGL(d,\Kb)$, where $\Kb$ is either the real numbers or the complex numbers, is given by the set of real-valued diagonal matrices with trace zero:
$$\mathfrak a=\{\mathrm{diag}(A_1,\ldots, A_d)\in\mathfrak{sl}(d,\Rb)\ |\ A_1+\cdots +A_d=0\}.$$
The space $\mathfrak{a}^*$ of linear functionals on $\mathfrak{a}$ is generated by the simple roots $\{\alpha_i\}_{i=1}^{d-1}$
where $\alpha_i(A)=A_{i}-A_{i+1}$ and the standard {\em positive Weyl chamber} is the subset where all the
simple roots are non-negative:
$$\mathfrak a^+=\{A\in \mathfrak a\ |\ A_{1}\ge A_{2}\ge\cdots\ge A_{d} \}.$$
We will be especially interested in the set  $\left(\mathfrak{a}^*\right)^+$ of linear functionals which are strictly positive on the interior of the positive Weyl chamber,
i.e. 
$$\left(\mathfrak{a}^*\right)^+=\left\{\phi\in\mathfrak{a}^*-\{0\}\ :\ \phi=\sum_{j=1}^{d-1} c_j\alpha_j\ \text{ such that }c_j\ge 0\ \forall j\right\}.$$

Given $g\in\PGL(d,\Kb)$, let $\bar{g}\in\GL(d,\Kb)$ be a representative of $g$ whose determinant has modulus $1$. Then let  
\[\lambda_1(g)\ge\dots\ge\lambda_d(g)>0\] 
denote the modulus of the generalized eigenvalues of $\bar{g}$, and let 
\[\sigma_1(g)\ge\dots\ge\sigma_d(g)>0\]
denote the singular values of $\bar{g}$. The {\em Jordan projection} and \emph{Cartan projection} 
\[\nu,\kappa:\mathsf{PGL}(d,\mathbb K)\to \mathfrak a^+\] 
are respectively given by
\[\nu(g)=\mathrm{diag}(\log\lambda_1(g),\ldots,\log\lambda_d(g))\,\text{ and }\,\kappa(g)=\mathrm{diag}(\log\sigma_1(g),\ldots,\log\sigma_d(g)).\]
We recall that if $g\in\mathsf{PGL}(d,\mathbb K)$, then $g=\ell a m$, where $\ell, m \in\mathsf{PO}(d,\mathbb K)$ and
$a\in \mathrm{exp}(\mathfrak{a}^+)$. If $\alpha_k(g)>0$, then $U_k(g)=\ell(\langle e_1,\ldots,e_k\rangle)$ is well-defined,
and is the image of the $k$-plane which is ``stretched the most'' by $A$.

If $\phi\in\left(\mathfrak{a}^*\right)^+$ and $\Gamma\subset\mathsf{PGL}(d,\mathbb K)$ is discrete, we define its \emph{$\phi$-Poincar\'e series}
$$Q_\Gamma^\phi(s)=\sum_{\gamma\in\Gamma}e^{-s\phi(\kappa(\gamma))}$$
and its {\em $\phi$-critical exponent}
$$\delta^\phi(\Gamma)=\inf\{s\ :\ Q_\Gamma^\phi(s)<+\infty\}.$$
If $\rho:\Gamma\to\PGL(d,\Kb)$ is a representation with discrete image and finite kernel, we define its {\em $\phi$-critical exponent} by $\delta^\phi(\rho)=\delta^\phi(\rho(\Gamma))$. For all $\phi\in\left(\mathfrak{a}^*\right)^+$, the $\phi$-critical exponent for any Hitchin representation is finite, see Corollary \ref{PSgen}.

If $\Gamma\subset \mathsf{PGL}(d,\mathbb K)$ is a subgroup, we say that it is {\em $P_k$-divergent} if
$\alpha_k(\kappa(\gamma_n))\to\infty$ for any sequence $\{\gamma_n\}$ in $\Gamma$ of pairwise distinct elements. Notice that since
$\alpha_k(\kappa(\gamma))=\alpha_{d-k}(\kappa(\gamma^{-1}))$, $\Gamma$ is $P_k$-divergent if and only if
it is $P_{d-k}$-divergent. (Guichard--Wienhard \cite{guichard-wienhard} refer to $P_k$-divergent groups as $\alpha_k$-divergent, while Kapovich--Leeb--Porti \cite{KLP} call them
$\tau_{\rm mod}$-regular.)

Suppose that $\theta=\{k_1<k_2<\cdots<k_r\}$ is a symmetric subset of $\Delta=\{1,\ldots,d-1\}$ (i.e. $k\in\theta$ if and only if $d-k\in\theta$). 
The set of $\theta$-flags is the set of partial flags
with subspaces in all dimensions contained in $\theta$, i.e.
$$\mathcal F_\theta(\Kb^d)=\{F^{k_1}\subset \cdots\subset F^{k_r}\subset\mathbb K^d\ :\ \mathrm{dim}(F^{k_i})=k_i\}.$$
In particular, $\mathcal F=\mathcal F_\Delta(\Kb^d)$. When the context is clear, we write $\mathcal F_\theta = \mathcal F_\theta(\Kb^d)$. 

We say that a subset $X$ of $\mathcal F_\theta$ is {\em transverse} if
whenever $k\in\theta$ and $F,G\in X$ are distinct, then $F^k$ and $G^{d-k}$ are transverse.
If $\Gamma$ is $P_k$-divergent for all $k\in\theta$, then
$$U_\theta(\gamma)=(U_k(\gamma))_{k\in\theta}$$ 
is well-defined for all but finitely many $\gamma$ in $\Gamma$ and we can define the $\theta$-limit set $\Lambda_\theta(\Gamma)$
to be the set of limit points of the set $\{U_\theta(\gamma):\gamma\in\Gamma\}$. 
We then say that $\Gamma$ is {\em $P_\theta$-transverse} if it is $P_k$-divergent for all $k\in\theta$ and
$\Lambda_\theta(\Gamma)$ is a transverse subset of $\mathcal F_\theta$. (In Kapovich--Leeb--Porti \cite{KLP}, $P_k$-transverse subgroups are called 
$\tau_{\rm mod}$-regular and $\tau_{\rm mod}$-antipodal.)

Motivated by results of Danciger--Gu\'eritaud--Kassel \cite{DGK} and Zimmer \cite{Zimmer}, to study of transverse groups we introduce and develop a theory of projectively visible groups. 
We say that a discrete subgroup $\Gamma_0$ of $\mathsf{PGL}(d,\mathbb R)$ is {\em projectively visible} if it 
preserves a properly convex domain $\Omega$ in $\mathbb P(\mathbb R^d)$, every point in its full orbital limit set 
$$\Lambda_\Omega(\Gamma_0)=\left\{z\in\partial\Omega\ |\ z=\lim\gamma_n(x)\ \text{ for some }x\in\Omega\text{ and some }\{\gamma_n\}\subset\Gamma_0\right\}$$
has a unique supporting hyperplane to $\Omega$, and any two points in $\Lambda_\Omega(\Gamma_0)$ are joined
by a projective line segment in $\Omega$. 
As a key tool in our work, we show that every $P_\theta$-transverse subgroup  is the image of 
a projectively visible group $\Gamma_0$.

\begin{theorem}[see Theorem \ref{transverse image of visible1}]
If $\Gamma\subset\PGL(d,\Kb)$ is $P_\theta$-transverse, then there exists a properly convex domain $\Omega \subset \Pb(\Rb^{d_0})$ (for some $d_0\in\mathbb N$),
a projectively visible subgroup $\Gamma_0 \subset \mathrm{Aut}(\Omega)$, a faithful  representation
$\rho: \Gamma_0 \rightarrow \PGL(d, \Kb)$ and a $\rho$-equivariant continuous map $\xi:\Lambda_{\Omega}(\Gamma_0)\to\mathcal F_\theta$ so that 
$\rho(\Gamma_0)=\Gamma$ and $\xi(\Lambda_{\Omega}(\Gamma_0)) = \Lambda_{\theta}(\Gamma)$.
\end{theorem}

The domain $\Omega$ comes equipped with a natural, projectively invariant  Finsler metric $d_\Omega$, called  the Hilbert metric. In general, $(\Omega, \d_\Omega)$ will not be Gromov hyperbolic, but it has enough hyperbolicity to play the role that the Cayley graph does when studying Anosov groups. Similarly, the restriction of the Hilbert geodesic flow to the convex hull of $\Lambda_\Omega(\Gamma_0)$ plays the role of the
Gromov geodesic flow of a hyperbolic group.

We will observe, in Proposition \ref{prop: limit set convergence action}, that if  $\Gamma$ is $P_\theta$-transverse, then it acts on its $\theta$-limit set as a convergence group, i.e.
if $\{\gamma_n\}$ is a sequence of distinct elements in $\Gamma$, then there are points $x,y\in\Lambda_\theta(\Gamma)$ and a subsequence,
still called $\{\gamma_n\}$, so that $\gamma_n(z)\to x$ for all $z\in\Lambda_\theta(\Gamma)\setminus\{y\}$.
We recall that a point $x\in\Lambda_\theta(\Gamma)$  is a  {\em conical limit point} if there exists $a,b\in\Lambda_\theta(\Gamma)$ and a sequence
$\{\gamma_n\}$ in $\Gamma$ so that $\gamma_n(x)\to a$ and $\gamma_n(y)\to b$ for all $y\in\Lambda_\theta(\Gamma)\setminus\{x\}$. The
set of conical limit points for the action of $\Gamma$ on $\Lambda_\theta(\Gamma)$ is called the {\em $\theta$-conical limit set} and is denoted $\Lambda_{\theta,c}(\Gamma)$.

In many situations, the complement of the conical limit set is countable and consists of fixed points of ``weakly unipotent'' elements. Most classically, Beardon--Maskit \cite{beardon-maskit}
proved that this characterized geometrically finite Kleinian groups, see also Bishop \cite{bishop-conical}. This property also holds for $P_\theta$-Anosov images of geometrically finite Fuchsian groups,
and images of relatively dominated and relatively Anosov representations. Our first main result is an upper bound for the Hausdorff dimension of the conical limit
set of a transverse group in terms of the simple root critical exponent. If $k\in\theta$ and $\Gamma$ is $P_\theta$-transverse,
we define $\Lambda_{k,c}(\Gamma)=\pi_k(\Lambda_{\theta,c}(\Gamma))$, where $\pi_k:\mathcal F_\theta\to\mathrm{Gr}_k(\mathbb R^d)$ is the projection map.
Our result generalizes work of Glorieux--Montclair--Tholozan \cite{GMT} and Pozzetti--Sambarino--Wienhard \cite{PSW1} in the Anosov setting.

\begin{theorem}[see Corollary \ref{prop:upper_bound}]
If $\Gamma \subset \PGL(d,\Kb)$ is $P_{k,d-k }$-transverse, then
\[{\rm dim}_H\left(\Lambda_{k,d-k,c}(\Gamma)\right) \le \delta^{\alpha_k}(\Gamma).\]
In particular, ${\rm dim}_H\left(\Lambda_{k,c}(\Gamma)\right) \le \delta^{\alpha_k}(\Gamma)$.
\end{theorem} 

As a consequence we obtain the following generalization of results of Burger \cite{burger}, Glorieux--Montclair--Tholozan \cite{GMT} and
Kim--Minsky--Oh \cite{KMO} for pairs of convex cocompact representations.

\begin{theorem}[see Theorem~\ref{thm:KMO general in body}]
\label{KMO general}
Let $\rho_1:\Gamma\to \mathsf{SO}(d_1-1,1)$ and $\rho_2:\Gamma\to \mathsf{SO}(d_2-1,1)$ be geometrically finite representations so that
$\rho_1(\alpha)$ is parabolic if and only if $\rho_2(\alpha)$ is parabolic. If we regard $\rho=\rho_1\oplus\rho_2$ as a representation into
$\mathsf{PSL}(d_1+d_2,\mathbb R)$, then
$${\rm dim}_H\left(\Lambda_{2}(\rho(\Gamma))\right) =\max\left\{{\rm dim}_H\left(\Lambda_{1}(\rho_1(\Gamma))\right), {\rm dim}_H\left(\Lambda_{1}(\rho_2(\Gamma))\right)\right\}.$$
\end{theorem}

Notice that our assumptions imply that there is a H\"older homeomorphism $\xi:\Lambda_1(\rho_1(\Gamma))\to\Lambda_1(\rho_2(\Gamma))$ which is
typically not a diffeomorphism and that $\Lambda_2(\rho(\Gamma))$ can be smoothly identified with the graph of $\xi$. So, this is another instance,
surprisingly common in Higher Teichm\"uller theory, where a H\"older homeomorphism  fails to change Hausdorff dimension.

\medskip

Following Pozzetti--Sambarino--Wienhard \cite{PSW1} we say a group $\Gamma$ is \emph{$(1,1,q)$-hypertransverse} if it is $P_\theta$-transverse for some $\theta$ 
containing $1$ and $q$, and 
$$F^1+G^1+H^{d-q}$$
is a direct sum for all pairwise distinct $F,G,H\in\Lambda_\theta(\Gamma)$. 

Examples of $(1,1,2)$-hypertransverse groups include all (images of) exterior powers of Hitchin representations and examples of $(1,1,d-1)$-hypertransverse groups include all discrete subgroups of 
$\mathsf{PO}(d-1,1)$. Further, by definition, images of $(1,1,q)$-hyperconvex representations, in the sense of Pozzetti--Sambarino--Wienhard \cite{PSW1}, and their subgroups are also $(1,1,q)$-hypertransverse groups, so the $(1,1,q)$-hyperconvex representations into $\mathsf{SU}(n,1)$, $\mathsf{Sp}(n,1)$ and $\mathsf{SO}(p,q)$ constructed in \cite{PSW1} can be used to construct $(1,1,q)$-hypertransverse groups. 

The following result generalizes results of both Pozzetti--Sambarino--Wienhard \cite{PSW1} 
and Bishop--Jones \cite{bishop-jones} 
and uses ideas from both of their proofs.

\begin{theorem}[see Theorem~\ref{thm:lower_bd}]\label{thm:BishopJones}
Suppose that $\Gamma\subset\PGL(d,\mathbb K)$ is $(1,1,q)$-hypertransverse and
\begin{align*}
\sigma_2(\gamma) = \sigma_q(\gamma) 
\end{align*}
for all $\gamma \in \Gamma$. Then 
\begin{align*}
{\rm dim}_H(\Lambda_{1,c}(\Gamma)) = \delta^{\alpha_1}(\Gamma).
\end{align*} 
\end{theorem}

In order to apply this theorem, we  first observe that the image of a Hitchin representation has limit set of Hausdorff dimension  at most  1.

\begin{proposition}[see Proposition~\ref{pairwise transverse}]
\label{hitchinlimitset}
If $\Gamma$ is a Fuchsian group and $\rho:\Gamma\to \mathsf{PSL}(d,\mathbb R)$ is a Hitchin representation, then
$$\mathrm{dim}_{H}(\Lambda_\Delta(\rho(\Gamma)))\le 1.$$
\end{proposition}

After we observe that exterior powers of Hitchin representations are $(1,1,2)$-hypertransverse, 
Theorem \ref{thm:BishopJones} and Propositions \ref{hitchinlimitset}
have the following consequence for the simple root entropies of Hitchin representations.

\begin{corollary}[see Corollary~\ref{PSgen in body} and Theorem \ref{thm:main in body}]
\label{PSgen}
 If $\Gamma$ is a Fuchsian group and $\rho:\Gamma\to\PSL(d,\mathbb R)$ is a Hitchin representation, then
$$\delta^{\alpha_k}(\rho)\le 1$$
for all $k\in\Delta$, and equality holds if $\Gamma$ is a lattice. Furthermore, if $\phi = \sum_{k\in\Delta} c_k \alpha_k \in \left(\mathfrak{a}^*\right)^+$ is non-zero, then 
\begin{align*}
\delta^{\phi}(\rho) \le \frac{1}{c_1+\dots + c_{d-1}}.
\end{align*}
\end{corollary} 

In the case when $\Gamma$ is a uniform lattice, Corollary \ref{PSgen} was previously established by Potrie--Sambarino~\cite{potrie-sambarino}
by quite different methods, and in the convex co-compact case can also be deduced from the work of Pozzetti--Sambarino--Wienhard \cite{PSW1}.
Corollary \ref{PSgen} plays a central role in the construction of the (first) simple root pressure metric for Hitchin components of Fuchsian lattices,
see Bray--Canary--Kao--Martone \cite{BCKM2}.

If $\phi\in\left(\mathfrak{a}^*\right)^+$,  then the $\phi$-critical exponent and
the $\phi$-entropy  of a Hitchin representation of a geometrically finite Fuchsian group agree. 
We define
the {\em $\phi$-length} of an element of $g \in\mathsf{PGL}(d,\mathbb R)$ as
$$\ell^\phi(g)=\phi(\nu(g)).$$
If $\Gamma\subset\mathsf{PSL}(2,\mathbb R)$ is Fuchsian and $\rho:\Gamma\to \mathsf{PGL}(d,\mathbb R)$ is a representation, we define the \emph{$\phi$-entropy} as 
$$h^\phi(\rho)=\limsup_{T\to\infty} \frac{\log \# R_T^\phi(\rho)}{T}\qquad\text{where}\qquad R_T^\phi(\rho)=\left\{ [\gamma]\in[\Gamma_{hyp}]\ :\  \ell^\phi(\rho(\gamma))\le T\right\}$$
where $[\Gamma_{hyp}]$ is the set of conjugacy classes of hyperbolic elements of $\Gamma$. 
When $\rho$ is a Hitchin representation of a geometrically finite Fuchsian group, the $\limsup$ in the definition of $h^\phi(\rho)$ holds as a limit, and is always positive and finite
(see \cite{BCKM}).

\begin{proposition}[see Proposition~\ref{entropy and critical exponent}]
\label{critandent}
If $\Gamma\subset\mathsf{PSL}(2,\mathbb R)$ is geometrically finite, $\phi\in(\mathfrak{a}^*)^+$ and $\rho:\Gamma\to \mathsf{PSL}(d,\mathbb R)$ is a Hitchin representation, then
$$\delta^\phi(\rho)=h^\phi(\rho).$$
\end{proposition}

If $\Gamma$ is geometrically finite, but not a lattice, then $\Gamma$ is contained in a lattice $\Gamma^D$ such that any Hitchin representation
$\rho:\Gamma\to\mathsf{PSL}(d,\mathbb R)$ extends to a Hitchin representation $\rho^D:\Gamma^D\to\mathsf{PSL}(d,\mathbb R)$ (see Proposition \ref{prop: double}).
We may then apply classical arguments,
which go back to Furusawa \cite{furusawa} to establish that the $\phi$-critical exponent of $\rho(\Gamma)$ is strictly less than the critical exponent of $\rho^D(\Gamma)$.
(One may also view the proof as a concrete version of an argument in Dal'bo--Otal--Peign\'e \cite[Thm.\ A]{DOP} who make use of Patterson-Sullivan measure instead of
working directly with the Poincar\'e series.)

\begin{proposition}[see Proposition~\ref{prop: entropy drop}]
Suppose that $\rho:\Gamma\to\mathsf{PSL}(d,\mathbb R)$ is a Hitchin representation of a Fuchsian group $\Gamma$, $G$ is an infinite index, finitely generated subgroup
of $\Gamma$, and $\phi\in\left(\mathfrak{a}^*\right)^+$, then
$$\delta^\phi(\rho|_G)<\delta^\phi(\rho).$$
\end{proposition}

As an immediate consequence we obtain:

\begin{corollary}
\label{PSgen2}
If $\Gamma \subset \PSL(2,\Rb)$ is geometrically finite, but not a lattice, and $\rho:\Gamma\to\PSL(d,\mathbb R)$ is a Hitchin representation, then
$$\delta^{\alpha_k}(\rho) < 1$$
for all $k\in\Delta$.
\end{corollary}

Once we have established Corollaries \ref{PSgen} and \ref{PSgen2} we may follow a similar outline of proof  as in \cite{potrie-sambarino} to establish our
main theorem. The key step is a convexity result for the behavior of entropy over the space of functionals, which generalizes  \cite[Cor.\ 4.9]{sambarino-hyperconvex},
see also \cite{sambarino-dichotomy}.

\begin{theorem}[see Theorem~\ref{quint indicator}]
\label{quintintro}
Suppose that $\Gamma$ is a geometrically finite Fuchsian group and
$\rho:\Gamma\to\mathsf{PSL}(d,\mathbb R)$ is a Hitchin representation. Then
$$\mathcal{Q}_\theta(\rho)=\{\phi\in\left(\mathfrak{a}^*\right)^+ \ |\ h^\phi(\rho)=1\}$$
is a closed subset of a  concave, analytic submanifold of  $\mathfrak{a}^*$.
Moreover, if $\phi_1,\phi_2\in\mathcal{Q}(\rho)$, then the line segment in $\mathfrak{a}^*$ between $\phi_1$ and $\phi_2$ lies in $\mathcal{Q}(\rho)$
if and only if 
$$\ell^{\phi_1}(\rho(\gamma))=\ell^{\phi_2}(\rho(\gamma))$$
for all $\gamma\in\Gamma$.
\end{theorem}

As in~\cite{potrie-sambarino}, Theorem~\ref{thm:main} also implies a rigidity result for the symmetric space critical exponent. Given a discrete subgroup $\Gamma \subset \PSL(d,\Rb)$,  
the \emph{symmetric space critical exponent}, denoted $\delta_X(\Gamma) \in [0,\infty)$, is the critical exponent of the series 
\begin{align*}
Q_{\Gamma,x_0}^X(s):=\sum_{\gamma \in \Gamma} e^{-s \d_X(\gamma(x_0), x_0)} 
\end{align*}
(which is independent of the choice of $x_0 \in X$) where $\d_X$ is the symmetric space distance on $X = \PSL(d,\Rb) / \mathsf{PSO}(d)$ scaled so that the embedding 
$\Hb^2 \hookrightarrow X$ induced by some (hence any) irreducible representation $\PSL(2,\Rb) \rightarrow \PSL(d,\Rb)$ is isometric. 

\begin{corollary}[see Corollary~\ref{cor:symmetric rigidity in body}]\label{cor:symmetric rigidity} If $\Gamma \subset \PSL(2,\Rb)$ is geometrically finite and $\rho : \Gamma \rightarrow \PSL(d, \Rb)$ is Hitchin, then 
\begin{align*}
\delta_X(\rho) \le 1.
\end{align*}
Moreover, $\delta_X(\rho) =1$ if and only if $\Gamma$ is a lattice and $\rho(\Gamma)$ lies in the image of an irreducible representation $\PSL(2,\Rb)\to\PSL(d,\Rb)$.
\end{corollary}


\section{Background}


\subsection{Linear Algebra}\label{sec: Lie}
We recall some basic notation and terminology from the Lie theory of $\PGL(d,\Kb)$, where $\Kb$ is either the real numbers or the complex numbers. We begin
by discussing a subspace of the Cartan subspace naturally associated to a symmetric subset $\theta$ of  $\Delta=\{1,\dots,d-1\}$.
Specifically, let
$$\mathfrak a_\theta=\{ A\in \mathfrak a\ |\ \alpha_k(A)=0\ \text{ if }\ k\notin\theta\}$$ 
and let $p_\theta:\mathfrak a\to \mathfrak a_\theta$ be the projection map such that $\omega_k\circ p_\theta=\omega_k$ for all $k\in\theta$
where $\omega_k\in\mathfrak{a}^*$ is the $k^{\rm th}$ {\em fundamental weight} given by
$$\omega_k(A)=A_1+\cdots+A_k.$$
Then $p_\theta^*:\mathfrak{a}_\theta^*\to\mathfrak a^*$ identifies $\mathfrak{a}_\theta^*$ as the subspace of $\mathfrak{a}^*$ spanned by $\{\omega_k\}_{k\in\theta}$. 
In particular, every $A\in \mathfrak a_\theta$ is determined by the tuple $(\omega_k(A))_{k\in\theta}$. We will be interested in the set
$(\mathfrak{a}_\theta^*)^+$ of vectors in $\mathfrak a^*_\theta$ that are strictly positive on the interior of the $\theta$-positive Weyl chamber $p_\theta(\mathfrak{a}^+)$. Explicitly,
$$(\mathfrak{a}_\theta^*)^+=\left\{ \phi=\sum_{k\in\theta} c_k\omega_k\ :\ \phi\ne0 \text{ and } c_k\ge 0 \text{ for all } k\in\theta\right\}.$$

We say that $g\in\mathsf{PGL}(d,\mathbb K)$ is {\em $\theta$-proximal} if $\alpha_k(\nu(\gamma))>0$ for all $k\in\theta$.
We then define the {\em $\theta$-Benoist limit cone} of a discrete subgroup $\Gamma$ of $\mathsf{PGL}(d,\mathbb K)$ to be
the closure of the set of rays determined by $\theta$-Jordan projections of $\theta$-proximal elements, i.e. 
$$
\mathcal B_\theta(\Gamma)=\overline{\left\{\mathbb R^+p_\theta(\nu(\gamma))\ :\ \gamma\in\Gamma\ \text{ is } \theta\text{-proximal}\right\} }\subset\mathfrak{a}_\theta.
$$
We will then be interested in the open  set of linear functionals which are strictly positive on the $\theta$-Benoist limit cone (except at $0$),
$$
\mathcal B_\theta^+(\Gamma)=\{\phi\in\mathfrak{a}_\theta^*\ :\ \phi(A)>0 \text{ for all }  A\in \mathcal B_\theta(\Gamma)-\{0\}\}.
$$
When $\theta=\Delta$, we will use the standard notation $\mathcal B(\Gamma)=\mathcal B_\Delta(\Gamma)$ and
$\mathcal B^+(\Gamma)=\mathcal B^+_\Delta(\Gamma)$. Also, if $\rho:\Gamma\to\PGL(d,\Rb)$ is a representation, we denote $\Bc_\theta(\rho)=\Bc_\theta(\rho(\Gamma))$, $\Bc_\theta^+(\rho)=\Bc_\theta^+(\rho(\Gamma))$, $\Bc(\rho)=\Bc(\rho(\Gamma))$ and $\Bc^+(\rho)=\Bc^+(\rho(\Gamma))$.

Recall that the angle $\theta \in [0,\pi/2]$ between two lines $L_1,L_2\in\mathbb P(\mathbb K^d)$ is defined by 
$$
\cos(\theta) = \frac{\abs{\ip{v_1,v_2}}}{\norm{v_1}\norm{v_2}}
$$
where $v_1 \in L_1, v_2 \in L_2$ are some (any) non-zero vectors. Further, this angle defines a distance, denoted $\d_{\Pb(\Kb^d)}$, on $\mathbb P(\mathbb K^d)$ which is induced by a Riemannian metric. 

There is a natural smooth embedding  of $\mathrm{Gr}_k(\mathbb K^d)$ into $\mathbb P(\bigwedge^k\mathbb K^d)$ which
takes a $k$-subspace with basis $\{b_1,\ldots,b_k\}$ to the line spanned by $b_1\wedge\cdots\wedge b_k$. We then endow $\mathrm{Gr}_k(\mathbb K^d)$ with the distance $\d_{\mathrm{Gr}_k(\mathbb K^d)}$ obtained by pulling back the angle metric on $\mathbb P(\bigwedge^k\mathbb K^d)$. We then give $\prod_{k\in\theta}\mathrm{Gr}_k(\mathbb K^d)$ the product
metric and give $\mathcal F_\theta $ the metric, denoted $\d_{\Fc_\theta}$, it inherits as a subset of $\prod_{k\in\theta}\mathrm{Gr}_k(\mathbb K^d)$.

We will use the following lemma a number of times (for a proof, see for instance~\cite[Prop.\ A.1]{CZZ3}). 

\begin{lemma}\label{lem: KAK}
Let $\theta\subset\Delta$ be symmetric, $F^+,F^-\in\Fc_\theta$, and $\{g_n\}$ a sequence in $\PGL(d,\Kb)$. The following are equivalent:
\begin{enumerate}
\item $\alpha_k(\kappa(g_n))\to\infty$ for all $k\in\theta$, $U_\theta(g_n)\to F^+$, and $U_\theta(g_n^{-1})\to F^-$,
\item $g_n(F)\to F^+$ uniformly on compact subsets of $\{  F\in\Fc_\theta :  F \text{ transverse to } F^-\}$.
\end{enumerate}
\end{lemma}

We will also use the following variant  (for a proof, see for instance~\cite[Lem.\ A.3]{CZZ3}) . 

\begin{lemma}\label{lem: KAK 2}
Let $\theta\subset\Delta$ be symmetric, $F^-\in\Fc_\theta$, and $\{g_n\}$ a sequence in $\PGL(d,\Kb)$. If there is an open set $\mathcal{O} \subset \Fc_\theta$ where $g_n(F) \rightarrow F^+$ for all $F \in \mathcal{O}$, then $\alpha_k(\kappa(g_n))\to\infty$ for all $k\in\theta$ and $U_\theta(g_n)\to F^+$.
\end{lemma}

\subsection{Cusped Anosov representations of geometrically finite Fuchsian groups}

Cusped Anosov representations of geometrically finite Fuchsian groups were introduced in \cite{CZZ} as natural generalizations of Anosov representations which take
parabolic elements to elements whose (generalized) eigenvalues all have modulus $1$.  These representations are also relatively Anosov in the sense of Kapovich--Leeb \cite{kapovich-leeb}
and relatively dominated in the sense of Zhu \cite{feng}.
If $\Gamma$ is a geometrically finite Fuchsian group with limit set $\Lambda(\Gamma)\subset\partial\mathbb H^2$, then
a representation $\rho:\Gamma\to\mathsf{PGL}(d,\mathbb K)$ is said to be {\em $P_\theta$-Anosov} if
there exists a continuous $\rho$-equivariant map $\xi_\rho:\Lambda(\Gamma)\to \mathcal F_\theta$ such that
 \begin{enumerate}
 \item $\xi_\rho$ is \emph{transverse}, i.e. if $x,y $ are distinct points in $\Lambda(\Gamma)$, then $\xi_\rho^k(x)\oplus\xi_\rho^{d-k}(y)=\mathbb K^d$ for all $k\in\theta$.
 \item
 $\xi_\rho$ is {\em strongly dynamics preserving}, i.e. if $b_0\in\mathbb H^2$ and
$\{\gamma_n\}$ is a sequence in $\Gamma$  such
 that $\gamma_n(b_0)\to x\in\Lambda(\Gamma)$ and $\gamma_n^{-1}(b_0)\to y\in\Lambda(\Gamma)$, then for all $k\in\theta$ and
 $V\in\mathrm{Gr}_k(\mathbb K^d)$ that is transverse to $\xi_\rho^{d-k}(y)$, we have $\rho(\gamma_n)(V)\to\xi_\rho^k(x)$.
 \end{enumerate}
(To be precise, in~\cite{CZZ} Anosov representations were defined in terms of the exponential contraction of a linear flow on a vector bundle associated to a representation and this was shown to be equivalent to the definition above.)

 If $\Gamma$ contains a parabolic element, we refer to such representations as cusped Anosov when we want to distinguish them from
 traditional Anosov representations (which cannot contain unipotent elements in their image).
We recall the properties of cusped Anosov representations we will need in our work. For any $x,y\in\Hb^2$, let $\d(x,y)$ denote the hyperbolic distance between $x$ and $y$, and for any $\gamma\in\PSL(2,\Rb)$, let $\ell(\gamma)$ denote the minimal translation distance of $\gamma$ acting on $\Hb^2$.

\begin{theorem}[Canary--Zhang--Zimmer \cite{CZZ}]
\label{PkAnosov properties}
If $\Gamma$ is a geometrically finite Fuchsian group, $\rho:\Gamma\to\mathsf{PGL}(d,\mathbb K)$ is $P_\theta$-Anosov and $b_0 \in \Hb^2$, then
\begin{enumerate}
\item
There exist $A,a>0$ so that if $\gamma\in\Gamma$, then
$$Ae^{a\d(b_0,\gamma(b_0))}\ge e^{\alpha_k(\kappa(\rho(\gamma)))}\ge \frac{1}{A}e^{\frac{\d(b_0,\gamma(b_0))}{a}}$$
for all $k\in\theta$.
\item
There exist $B,b>0$ so that if $\gamma\in\Gamma$, then
$$Be^{b \ell(\gamma)}\ge e^{\alpha_k(\nu(\rho(\gamma)))}\ge \frac{1}{B}e^{\frac{\ell(\gamma)}{b}}$$
for all $k\in\theta$.
In particular, if $\gamma\in\Gamma$ is hyperbolic, then $\rho(\gamma)$  is $\theta$-proximal.
\item $\rho$ has the \emph{$P_\theta$-Cartan property}, i.e.  whenever $\{\gamma_n\}$ is a sequence of distinct elements of $\Gamma$
such that $\gamma_n(b_0)$ converges to $z\in\Lambda(\Gamma)$, 
then $\xi_\rho(z)=\lim U_\theta(\rho(\gamma_n))$.
\end{enumerate}
\end{theorem}
Notice that part (3) is an immediate consequence of the definition and Lemma~\ref{lem: KAK}.

Bray--Canary--Kao--Martone \cite{BCKM} established counting and equidistribution results for cusped Anosov representations. We will
need the following counting result, which generalizes a result of Sambarino \cite[Thm.\ B]{sambarino-quantitative}, see also \cite{sambarino-dichotomy}, which applies in the convex cocompact case.
Notice that Theorem \ref{PkAnosov properties} implies that $\left(\mathfrak{a}_\theta^*\right)^+\subset \mathcal B_\theta^+(\rho)$ whenever $\rho$ is
$P_\theta$-Anosov. 

\begin{theorem}[{Bray--Canary--Kao--Martone \cite[Cor.\ 11.1]{BCKM}}]
\label{BCKM counting}
Suppose that  $\Gamma$ is a torsion-free, geometrically finite Fuchsian group and  $\rho:\Gamma\to\mathsf{PGL}(d,\mathbb K)$ is  $P_\theta$-Anosov.
If $\phi\in\mathcal B^+_\theta(\rho)$, then 
$${\displaystyle \lim_{t\to\infty}R_t^{\phi}(\rho)\frac{t\delta^\phi(\rho)}{e^{t\delta^\phi(\rho)}}}=1$$
where $R_t^\phi(\rho)=\#\left\{[\gamma]\in[\Gamma_{hyp}]\ :\ \phi(\nu(\rho(\gamma)))\le t\right\}$ and $[\Gamma_{hyp}]$ is the set of conjugacy classes of hyperbolic elements in $\Gamma$.
\end{theorem} 

\begin{remark*} To be precise, Corollary 11.1 in~\cite{BCKM} was stated for representations into $\SL(d,\Rb)$, but the same argument taken verbatim works for representations into $\PGL(d,\Kb)$ since the construction of the roof functions only involve the Cartan projection. 
\end{remark*}

\subsection{Cusped Hitchin representations of Fuchsian groups} In order to define Hitchin representations, we must first recall  the definition of a positive map.
Given a transverse pair of flags $F_1,F_2\in\mathcal F=\Fc_\Delta(\Rb^d)$, an ordered basis $\mathcal B=(b_1,\dots,b_d)$ for $\Rb^d$ is \emph{compatible} with 
$(F_1,F_2)$ if $b_i\in F_1^i\cap F_2^{d-i+1}$ for all $i\in\{1,\dots,d\}$. Given a basis $\Bc$, let $U_{>0}(\Bc)\subset\SL(d,\Rb)$ denote the set of unipotent elements that, when written in the basis $\Bc$, are upper triangular and all minors (which are not
forced to be 0 by the fact that the matrix is upper triangular) are strictly positive. Following Fock--Goncharov \cite{fock-goncharov},
we say that an ordered  $k$-tuple 
\[(F_1,F_2,\ldots, F_k)\] 
of flags in $\mathcal F$ is \emph{positive} if there exists an ordered basis $\mathcal B$ compatible with $(F_1,F_k)$, and elements 
$u_2,\ldots,u_{k-1}\in U_{>0}(\mathcal B)$ so that $F_i=u_{k-1}\cdots u_iF_k$ for all $i=2,\ldots,k-1$.

If $X$ is a subset of $\partial \Hb^2$, then a map $\xi:X\to\mathcal F$ is {\em positive} if $(\xi(x_1),\ldots,\xi(x_n))$ is a positive whenever 
$(x_1,\ldots, x_n)$ is a cyclically ordered subset of distinct points in $X$. If $\Gamma$ is a Fuchsian group, 
we say that a representation  $\rho:\Gamma\to\mathsf{PSL}(d,\mathbb R)$ is a \emph{Hitchin representation}  if there exists a continuous, positive,
$\rho$-equivariant map $\xi:\Lambda(\Gamma)\to \mathcal F$. When $\Gamma\subset \PSL(2,\Rb)$ is a cocompact torsion-free lattice, they agree with the representations introduced by 
Hitchin \cite{hitchin} and  studied by Labourie \cite{labourie-invent}. If  $\Gamma$ is torsion-free and convex cocompact, but not a lattice, they were studied by
Labourie--McShane \cite{labourie-mcshane}. 

We recall several well-known properties of positive tuples of flags that were observed by Fock--Goncharov \cite{fock-goncharov} 
(see also \cite[Appendix A]{SWZ} and \cite[Section 3.1--3.3]{KTZ}.)

\begin{lemma} \label{lem: basic positivity}
If $(F_1,\dots,F_k)$ is a positive tuple of flags in $\Fc$, then
\begin{enumerate}
\item $(F_2,\dots,F_k,F_1)$ is positive.
\item $(F_k,\dots,F_1)$ is positive.
\item $(F_{i_1},\dots,F_{i_\ell})$ is positive for all $1\le i_1<\dots<i_\ell\le k$.
\item $F_1^{n_1}\oplus\dots\oplus F_k^{n_k}=\Rb^d$ for any integers $n_1,\dots,n_k$ that sum to $d$.
\item $(F_1,\dots,F_n)$ is positive for all flags $F_{k+1},\dots,F_n\in\Fc$ such that $(F_1,F_i,F_k,F_{k+1},\dots,F_n)$ is positive for some $i\in\{2,\dots,k-1\}$.
\end{enumerate}
\end{lemma}

If $\Gamma$ is geometrically finite, then the following generalization of the main result in \cite{labourie-invent} is established in \cite{CZZ}.

\begin{theorem}{\rm (Canary--Zhang--Zimmer \cite[Thm.\ 7.1]{CZZ})}
If  $\Gamma\subset\mathsf{PSL}(2,\mathbb R)$ is geometrically finite and $\rho:\Gamma\to\mathsf{PSL}(d,\mathbb R)$ is a Hitchin representation, then $\rho$ is $P_\Delta$-Anosov.
\end{theorem}

\begin{remark*} To be precise, Theorem 7.1 in~\cite{CZZ} was stated for representations into $\SL(d,\Rb)$, but the same argument taken verbatim works for representations into $\PSL(d,\Rb)$, since the entire argument takes place in $\Fc(\Rb^d)$. 
\end{remark*}

\subsection{Properly convex domains}\label{sec: properly convex domains}
We briefly recall some standard facts about properly convex domains in projective space. 
A domain $\Omega\subset\mathbb P(\mathbb R^d)$ is {\em properly convex}  if it is a bounded convex subset of some affine chart $A$ for $\mathbb P(\mathbb R^d)$. 
If $x,y\in\overline{\Omega}$, let $[x,y]_\Omega$ denote the closed projective line segment in $\overline{\Omega}$ with endpoints $x$ and $y$. We also define $(x,y)_\Omega = [x,y]_\Omega -\{x,y\}$, $[x,y)_\Omega= [x,y]_\Omega - \{y\}$, and $(x,y]_\Omega = [y,x)_\Omega$. 

The \emph{radial projection maps based at $b_0 \in \Omega$} is the map 
\[\iota_{b_0}:\overline{\Omega}-\{b_0\}\to\partial\Omega\] 
defined by letting $\iota_{b_0}(z) \in \partial \Omega$ be the unique boundary point such that $z \in (b_0, \iota_{b_0}(z)]_\Omega$.

Every boundary point $x \in \partial \Omega$ of a properly convex domain is contained in a \emph{supporting hyperplane} $H$, that is: $H = \Pb(V)$ for some codimension one linear subspace $V \subset \Rb^d$, $x \in H$ and $H \cap \Omega = \emptyset$. When $x$ is contained in a unique supporting hyperplane, we say that $x$ is a \emph{$C^1$ point} of $\partial\Omega$ and denote this unique supporting hyperplane by $T_x \partial \Omega$.

A properly convex domain $\Omega$ has a natural projectively invariant
Finsler metric $\d_\Omega$, called the {\em Hilbert metric}, which is defined in terms of the cross ratio.
If $a,b\in\Omega$, then there is a projective line $\ell$ containing $a$ and $b$ which intersects $\partial\Omega$ at points $a'$ and $b'$ (ordered so
that $\{a',a,b,b'\}$ appear monotonically along $\ell$). Then
\[\d_\Omega(a,b)=\log\frac{|a'-b||b'-a|}{|a'-a||b'-b|},\]
where $|\cdot|$ denotes some (any) norm on some (any) affine chart containing $a',a,b,b'$. If $a,b\in\overline\Omega$, then the projective line segment $(a,b)_\Omega$ joining them is a geodesic
in the Hilbert metric, although geodesics need not be unique. We also let 
$$
B_\Omega(p,r) \subset \Omega
$$ 
denote the open ball of radius $r$ centered at $p \in \Omega$ with respect to the Hilbert metric.

We will use the following basic estimate several times: if $p_1,p_2, q_1, q_2 \in \Omega$, then 
\begin{equation} 
\label{eqn:Hilbert metric Hausdorff distance}
\d_\Omega^{\rm Haus}( [p_1,p_2]_\Omega, [q_1,q_2]_\Omega ) \le \max\{ \d_\Omega(p_1,q_1), \d_\Omega(p_2,q_2)\},
\end{equation}
(see for instance~\cite[Prop.\ 5.3]{islam-zimmer-lms}). In Equation~\eqref{eqn:Hilbert metric Hausdorff distance}, $\d_\Omega^{\rm Haus}$ denotes the Hausdorff distance induced by the Hilbert distance.

It is often useful to consider the dual of a properly convex domain.
Let $V=\mathbb R^d$. For any $k\in\Delta$, there is a natural  identification $\Gr_k(V^*) \cong \mathrm{Gr}_{d-k}(V)$ given by
\[\Span_{\Rb}(\alpha_1,\dots,\alpha_k)\mapsto\bigcap_{i=1}^k\ker(\alpha_i).\]
Also, we may identify $\Gr_k(V)$ with the set of $(k-1)$-dimensional projective hyperplanes in $\Pb(V)$.

The \emph{dual} of a properly convex domain $\Omega\subset\Pb(V)$ is the set
\[
\Omega^*=\left\{[f]\in\Pb(V^*):f(X)\neq 0\text{ for all }[X]\in\overline{\Omega}\right\}.
\]
We record the following standard fact for later use.

\begin{lemma}\label{lem: dual domain}
If $\Omega\subset\Pb(V)$ is a (non-empty) properly convex domain, then $\Omega^*$
is a (non-empty) properly convex domain in $\Pb(V^*)$. Furthermore, $(\Omega^*)^*=\Omega$ and ${\rm Aut}(\Omega)={\rm Aut}(\Omega^*)$ under the canonical identification $\PGL(V)\simeq\PGL(V^*)$. 
\end{lemma}

We refer the reader to Marquis \cite{marquis-around} for further discussion of the Hilbert metric on $\Omega$ and its automorphism group.

The following proposition describes the limiting behavior of divergent sequences in ${\rm Aut}(\Omega)$. 

\begin{proposition}\label{Prop: Islam-Zimmer}
Suppose that $\Omega\subset\Pb(\Rb^d)$ is a properly convex domain and $b_0\in\Omega$. Let $\{\gamma_n\}$ is a sequence in ${\rm Aut}(\Omega)$ such that $\gamma_n(b_0)\to x\in\partial\Omega$ and $\gamma_n^{-1}(b_0)\to y\in\partial\Omega$.
\begin{enumerate}
\item If $\gamma_n\to S\in \Pb(\mathrm{End}(\Rb^d))$, then $S(\Omega)\subset\partial\Omega$, $y\in \Pb(\ker S)$, and $\Pb(\ker S)\cap\Omega$ is empty. 
\item If $\alpha_1(\kappa(\gamma_n))\to\infty$ and $\gamma_n\to S\in\Pb(\mathrm{End}(\Rb^d))$, then $S(\Omega)=x$ and the $\Pb(\ker S)$ is a supporting hyperplane to $\Omega$ at $y$. In particular, $\gamma_n(b)\to x$ for all $b\in\mathbb P(\mathbb R^d)-\mathbb P(\ker S)$, and this convergence is locally uniform.
\item If $\alpha_1(\kappa(\gamma_n))\to\infty$ and $y$ is a $C^1$-point of $\partial\Omega$, then $\gamma_n\to S\in \Pb(\mathrm{End}(\Rb^d))$ with the defining property $S(\Omega)=x$ and $\Pb(\ker(S))=T_y\partial\Omega$.
\end{enumerate}
\end{proposition}

\begin{proof}
See Islam--Zimmer \cite[Prop.\ 5.6]{islam-zimmer-lms} for a proof of (1). The assumption that $\alpha_1(\kappa(\gamma_n))\to\infty$ implies that $S$ is the projectivization of a rank $1$ linear map, so (2) follows from (1). To prove (3), first observe that by taking a subsequence of $\{\gamma_n\}$, we may assume that $\gamma_n\to T\in\Pb(\mathrm{End}(\Rb^d))$. It then suffices to show that $T=S$. By (2), $T(\Omega)=x$ and $\Pb(\ker T)$ is a supporting hyperplane to $\Omega$ at $y$. Since $y$ is a $C^1$-point of $\partial\Omega$, $\Pb(\ker T)=T_y\partial\Omega$. Thus, $S=T$.
\end{proof}

\subsection{Special representations}\label{sec: special repns} We recall the skew-symmetric and symmetric tensor representations and their basic properties.

\subsubsection{Skew-symmetric tensors}\label{sec: skew-symmetric tensors}

Given a $\Kb$-vector space $V$, let $W=\bigwedge^k V$ be the vector space of skew-symmetric tensors of order $k$. Let $d = \dim_{\Kb}(V)$ and $D = \dim_{\Kb}(W)$, and let 
\[
E^k=E^k_V:\PGL(V)\to\PGL(W)
\] 
denote the representation defined by 
\[E^k(g)[v_1\wedge\dots\wedge v_k] = [(g v_1)\wedge\dots \wedge (g v_k)].\] 
It is straightforward to verify that $E^k$ is faithful and irreducible. We may also define a continuous, transverse, $E^k$-equivariant map
\[
\xi_{E^k}: \Fc_{k,d-k}(V) \rightarrow \Fc_{1,D-1}(W)
\] 
by 
\begin{align*}
\xi_{E^k}& \left( \Span(v_1,\dots,v_k), \Span(v_1,\dots,v_{d-k})\right) \\
& = \left(  [v_1\wedge\dots\wedge v_k],  \ker\left( w \in W \mapsto w \wedge v_1 \wedge \dots \wedge v_{d-k} \in \bigwedge\nolimits^d V \right)\right).
\end{align*}

In the special case when $V = \Kb^d$, the standard basis $(e_1,\dots,e_d)$ of $\Kb^d$ induces a standard basis $(e_{i_1}\wedge\dots\wedge e_{i_k})_{1\le i_1<\dots< i_k\le d}$ of $W$, and thus gives an identification $W\simeq\Kb^D$. Under this identification, we have 
\begin{align}
\label{eqn:singular_values_skew_repn} 
\alpha_1(\kappa(E^k(g))) =\alpha_k(\kappa(g)) \quad \text{and} \quad \sigma_1(E^k(g)) = (\sigma_1 \cdots \sigma_k)(g)
\end{align}
for all $g \in \PGL(d,\Kb)$.

\subsubsection{Hermitian symmetric tensors}\label{sec: symmetric tensors} Given a $\Kb$-vector space $V$, fix an (Hermitian) inner product $\ip{\cdot, \cdot}$ on $V$ and let $X \mapsto X^*$ denote the associated transpose on ${\rm End}(V)$, the space of $\Kb$-linear maps $V \rightarrow V$. Also, given $v \in V$ let $v^* \in V^*$ be the functional $v^* = \ip{\cdot, v}$. 

Let ${\rm Her}(V)$ denote the \textbf{real} vector space
$$
{\rm Her}(V) = \{ X \in {\rm End}(V) : X^* = X \},
$$ 
let $d = \dim_{\Kb}(V)$, and let $D= \dim_{\Rb}({\rm Her}(V))$. 

Next let
\[S_V:\PGL(V)\to\PGL({\rm Her}(V))\] 
denote the representation defined by  $S_V(g)(X)= g\circ X\circ g^*$. It is straightforward to verify that $S_V$ is faithful and irreducible. We may also define a continuous, transverse, $S_V$-equivariant map
\[\xi_{S_V}: \Fc_{1,d-1}(V) \rightarrow \Fc_{1,D-1}({\rm Her}(V))\] 
by 
$$
\xi_{S_V}([v],H) = \left([v\cdot v^*],  \Span \{  w\cdot v^*+v \cdot w^* : v \in V, w \in H \} \right).
$$

An element $X \in {\rm Her}(V)$ is \emph{positive definite} if $\ip{X(v),v}> 0$ for all non-zero $v \in V$. One can then verify that 
\[
\Omega_0 = \left\{ [X]\in\Pb({\rm Her}(V)) : X\text{ is positive definite}\right\}
\]
is a $S_V(\PGL(V))$-invariant properly convex domain.

In the special case when $V = \Rb^d$, the standard basis $(e_1,\dots,e_d)$ of $\Rb^d$ induces a standard basis $(e_i \cdot e_j^*+e_j \cdot e_i^*)_{1\le i\le j\le d}$ of ${\rm Her}(V)$ and in the special case when $V = \Cb^d$, the standard basis $(e_1,\dots,e_d)$ of $\Cb^d$ induces a standard basis 
$$
(e_i \cdot e_i^*)_{1 \le i \le d} \cup (e_i\cdot e_j^*+e_j\cdot e_i^*)_{1\le i <  j\le d}\cup  (ie_i\cdot e_j^*-ie_j\cdot e_i^*)_{1\le i <  j\le d}
$$ 
of ${\rm Her}(V)$ and thus gives an identification ${\rm Her}(V)\simeq\Rb^D$. Under these identifications, we have 
\begin{align}
\label{eqn:singular_values_symmetric_repn} 
\alpha_1(\kappa(S_V(g))) =\alpha_1(\kappa(g))
\end{align}
for all $g \in \PGL(d,\Kb)$. 


\section{Divergent and transverse subgroups}


In this section, we study divergent and transverse subgroups and their limit sets.
We exhibit examples and show that projectively visible subgroups are $P_{1,d-1}$-transverse.

\subsection{Properties of the limit set}
Let $\theta\subset\Delta$ be symmetric. 

\begin{proposition}\label{prop: limit set invariance}
If $\Gamma\subset\PGL(d,\Kb)$ is $P_\theta$-divergent, then $\Lambda_\theta(\Gamma)$ is $\Gamma$-invariant. 
\end{proposition}

\begin{proof} Fix $F^+ \in \Lambda_\theta(\Gamma)$ and $\gamma \in \Gamma$. By the definition of $\Lambda_\theta(\Gamma)$ and Lemma~\ref{lem: KAK} there exist $F^- \in \Lambda_\theta(\Gamma)$ and a sequence $\{\gamma_n\}$ in $\Gamma$ such that $\gamma_n(F) \rightarrow F^+$ uniformly on compact subsets of $\{ F \in \Fc_\theta : F \text{ transverse to } F^-\}$. Then
$$\gamma\gamma_n(F) \rightarrow \gamma(F^+)$$
uniformly on compact subsets of $\{ F \in \Fc_\theta : F \text{ transverse to } F^-\}$. Then Lemma~\ref{lem: KAK} implies that $U_\theta(\gamma \gamma_n) \rightarrow \gamma(F^+)$ and so $\gamma(F^+) \in \Lambda_\theta(\Gamma)$.
\end{proof}

The \emph{$P_\theta$-domain of discontinuity} for $\Gamma$, denoted $\Omega_\theta(\Gamma)$, is the set of flags in $\Fc_\theta$ that are transverse to every flag in 
$\Lambda_\theta(\Gamma)$. Since $\Lambda_\theta(\Gamma)$ is compact, observe that $\Omega_\theta(\Gamma)$ is a (possibly empty) open set. It also follows from
Proposition \ref{prop: limit set invariance} that $\Omega_\theta(\Gamma)$ is $\Gamma$-invariant. 
The following is a special case of a result from Guichard--Wienhard \cite[Thm.\ 7.4]{guichard-wienhard}. 

\begin{proposition}\label{prop: domain of discontinuity}
If $\Gamma\subset\PGL(d,\Kb)$ is $P_\theta$-divergent, then the action of $\Gamma$ on $\Omega_\theta(\Gamma)$ is properly discontinuous.
\end{proposition}

Since a hyperbolic group acts on
its Gromov boundary as a convergence group, Anosov groups act on their limit sets as convergence groups. We extend this property to transverse groups.
Recall that if $M$ is a compact metric space, a group $\Gamma$ acts as a \emph{convergence group} on $M$ if for any infinite sequence $\{\gamma_n\}$ of distinct elements in $\Gamma$, 
there exist some $x,y\in M$ and a subsequence $\{\gamma_{n_k}\}$ of $\{\gamma_n\}$ such that $\gamma_{n_k}(z)\to x$ for all $z\in M-\{y\}$. 

\begin{proposition} \label{prop: limit set convergence action}
If $\Gamma \subset  \PGL(d,\Kb)$ is $P_\theta$-transverse, then $\Gamma$ acts on $\Lambda_\theta(\Gamma)$ as a convergence group. 
\end{proposition}

\begin{proof}
Suppose $\{\gamma_n\}$ is an infinite sequence of distinct elements in $\Gamma$. By taking a subsequence, we may assume that there exists $F^+, F^-\in\Fc_\theta$ such that $U_{\theta}(\gamma_n)\to F^+$ and $U_{\theta}(\gamma_n^{-1})\to F^-$. Since $\Gamma$ is $P_\theta$-divergent, $\alpha_k(\kappa(\gamma_n))\to\infty$ for all $k\in\theta$. Since $\Gamma$ is $P_\theta$-transverse, by Lemma~\ref{lem: KAK}, $\gamma_n(F) \rightarrow F^+$ for all $F \in \Lambda_{\theta}(\Gamma) - \{F^-\}$. 
\end{proof} 

Next, we consider $P_1$-divergent subgroups $\Gamma$ that leave invariant a properly convex domain in $\Omega\subset\Pb(\Rb^d)$. 
Recall that the the orbital limit set of $\Gamma$ is
$$
\Lambda_\Omega(\Gamma) =\left\{z\in\partial\Omega\ |\ z=\lim\gamma_n(x)\ \text{ for some }x\in\Omega\text{ and some }\{\gamma_n\}\subset\Gamma\right\}.
$$
Also, recall that if $k \in \theta$, then $\pi_k : \Fc_\theta \rightarrow \Gr_k(\Kb^d)$ denotes the projection map. The next result relates the two limit sets of a $P_{1,d-1}$-divergent subgroup which preserves a properly convex domain $\Omega$. 

\begin{proposition}\label{prop: divergent groups preserving properly convex domains}
Suppose that $\Omega \subset \Pb(\Rb^d)$ is a properly convex domain and $\Gamma \subset \Aut(\Omega)$. If $\Gamma$ is $P_{1,d-1}$-divergent, then 
$$
\pi_1(\Lambda_{1,d-1}(\Gamma)) = \Lambda_\Omega(\Gamma),
$$
and, identifying $\Gr_{d-1}(\Rb^d) = \Pb(\Rb^{d*})$, 
$$
\pi_{d-1}(\Lambda_{1,d-1}(\Gamma)) = \Lambda_{\Omega^*}(\Gamma).
$$
In particular, if $F\in\Lambda_{1,d-1}(\Gamma)$, then $F^{d-1}$ is a supporting hyperplane to $\Omega$ at $F^1$.
\end{proposition} 

\begin{proof} The statements that $\pi_1(\Lambda_{1,d-1}(\Gamma)) = \Lambda_\Omega(\Gamma)$ and $\pi_{d-1}(\Lambda_{1,d-1}(\Gamma)) = \Lambda_{\Omega^*}(\Gamma)$ are dual, so we only prove the former.

Fix $F \in\Lambda_{1,d-1}(\Gamma)$ and let $\{\gamma_n\}$ be a sequence of distinct elements in $\Gamma$ with $U_{1,d-1}(\gamma_n)\to F$. Pass to a subsequence so
that $\gamma_n(b_0)\to x$ and $\gamma_n^{-1}(b_0)\to y$ (for some $b_0\in\Omega)$ and $\{\gamma_n\}$ converges to $S\in\mathbb P(\mathrm{End}(\mathbb R^d))$.
Proposition~\ref{Prop: Islam-Zimmer} part (2) implies that $S$ has rank 1 and $S(\Omega)=x$. Since $U_1(\gamma_n)\to F^1$, by Lemma \ref{lem: KAK} we must have $x=F^1$. Therefore,  $\pi_1(\Lambda_{1,d-1}(\Gamma)) \subset \Lambda_\Omega(\Gamma)$.

It remains to show that $\Lambda_\Omega(\Gamma) \subset \pi_1(\Lambda_{1,d-1}(\Gamma))$. 
Fix $x \in \Lambda_\Omega(\Gamma)$ and let $\{\gamma_n\}$ be a sequence in $\Gamma$ such that $\gamma_n(b_0) \rightarrow x$ for some $b_0 \in \Omega$. 
Passing to a subsequence we may  assume that $U_{1,d-1}(\gamma_n) \rightarrow F \in \Lambda_{1,d-1}(\Gamma)$, $\gamma_n^{-1}(b_0)\to y\in\Lambda_\Omega(\Gamma)$, and 
$\gamma_n\to S \in \Pb(\mathrm{End}(\Rb^d))$. Then by Proposition \ref{Prop: Islam-Zimmer} part (2), $S$ has rank 1 and $S(\Omega) = x$.
Since $U_1(\gamma_n)\to F^1$, we again see that $x=F^1$. 

Since $x$ is arbitrary, $\Lambda_\Omega(\Gamma) \subset \pi_1(\Lambda_{1,d-1}(\Gamma))$.
\end{proof} 

\subsection{Examples of transverse subgroups}
A large source of examples of $P_\theta$-transverse subgroups of $\PGL(d,\Kb)$ are the images of $P_\theta$-Anosov representations,  
$P_\theta$-relatively dominated  representations in the sense of Zhu \cite{feng}, $P_\theta$-asymptotically embedded representations in the
sense of Kapovich--Leeb \cite{kapovich-leeb} and any subgroup of one of these groups. Notice that these subgroups are not required to be finitely generated.

Another source of examples comes from the Klein-Beltrami model of real hyperbolic space. In particular, $\PO(d-1,1)$ preserves the properly convex domain 
$$
\Bb = \left\{ [x_1:\dots :x_{d-1}: 1] \in \Pb(\Rb^d) : \sum_{j=1}^{d-1} x_j^2 < 1\right\}
$$
whose boundary is smooth and contains no line segments. Then, any discrete subgroup of $\PO(d-1,1)$ is a $P_{1,d-1}$-transverse subgroup of $\mathsf{PGL}(d,\mathbb R)$, since it is $P_{1,d-1}$-divergent and its $P_{1,d-1}$-limit set is contained in $\{ (x, T_x \partial \Bb) : x \in \partial \Bb\}$ which is a transverse subset of $\mathcal F_{1,d-1}$.

Recall, that a discrete group $\Gamma \subset \PGL(d,\Rb)$ is \emph{projectively visible} if there exists a properly convex domain $\Omega \subset \Pb(\Rb^d)$ preserved by $\Gamma$ where the full orbital limit set $\Lambda_\Omega(\Gamma)$ satisfies: 
\begin{enumerate}
\item $(x,y)_\Omega \subset \Omega$ for all $x,y \in \Lambda_\Omega(\Gamma)$,
\item every $x \in \Lambda_\Omega(\Gamma)$ is a $C^1$ point of $\partial \Omega$. 
\end{enumerate}
In this case we also say that \emph{$\Gamma$ is a projectively visible subgroup of $\Aut(\Omega)$}. 

Notice that any discrete subgroup of $\PO(d-1,1)$ is a projectively visible subgroup of $\Aut(\Bb)$. Our next result gives some basic properties of projectively visible subgroups.

\begin{proposition}\label{prop:visible_implies_conv} 
Suppose that $\Omega \subset \Pb(\Rb^d)$ is a properly convex domain and
$\Gamma \subset \mathrm{Aut}(\Omega)$ is a projectively visible subgroup.  Fix $b_0\in\Omega$. 
\begin{enumerate}
\item $\Gamma\subset \PGL(d,\Rb)$ is a $P_{1,d-1}$-transverse subgroup and
\[\Lambda_{1,d-1}(\Gamma)=\{(x,T_x\partial\Omega):x\in\Lambda_\Omega(\Gamma)\}.\]
\item If $\{\gamma_n\}$ is a sequence in $\Gamma$ with $\gamma_n(b_0) \rightarrow x \in \Lambda_\Omega(\Gamma)$ and $\gamma_n^{-1}(b_0) \rightarrow y \in \Lambda_\Omega(\Gamma)$, then  
$$\gamma_n(F) \to (x,T_x\partial\Omega)$$
for all $F \in \Fc_{1,d-1}$ transverse to $(y,T_y \partial \Omega)$. Moreover, the convergence is locally uniform.
\item $\Gamma$ acts as a convergence group on $\Lambda_\Omega(\Gamma)$. 
\end{enumerate}
\end{proposition} 

\begin{proof} (1): First, we prove that $\Gamma$ is $P_{1,d-1}$-divergent. Let $\{\gamma_n\}$ be a sequence in $\Gamma$ of pairwise distinct elements. By taking a subsequence, we may assume that $\gamma_n(b_0) \rightarrow x \in \Lambda_\Omega(\Gamma)$, $\gamma_n^{-1}(b_0) \rightarrow y \in \Lambda_\Omega(\Gamma)$, and $\gamma_n\to S\in\Pb(\mathrm{End}(\Rb^d))$. By Proposition \ref{Prop: Islam-Zimmer} part (1), $S(\Omega)\subset\partial\Omega$, $y \in \mathbb P(\ker S)$, and 
$\mathbb P(\ker S) \cap \Omega$ is empty. Thus, if we pick any $w\in\Omega$, then
$$
S(w) = \lim_{n \rightarrow \infty} \gamma_n(w),
$$
so $S(\Omega) \subset \Lambda_\Omega(\Gamma)$. It follows that
$$
[S(w),x]_\Omega = [S(w), S(b_0)]_\Omega = S( [w,b_0]_\Omega) \subset S(\Omega) \subset \Lambda_\Omega(\Gamma),
$$
which implies that $S(w)=x$ because $\Gamma$ is projectively visible. Since $w \in \Omega$ was arbitrary, $S(\Omega) = \{x\}$. Since $\Omega$ is open, ${\rm im}\, S = x$. Thus, $S$ is the projectivization of a rank $1$ linear map, so  $\alpha_1(\kappa(\gamma_n))\to\infty$. Since $\{\gamma_n\}$ was arbitrary, $\Gamma$ is $P_{1,d-1}$-divergent.

Next, we prove that $\Gamma$ is $P_{1,d-1}$-transverse. Since each $x \in \Lambda_\Omega(\Gamma)$ has a unique supporting hyperplane, namely $T_x \partial \Omega$, Proposition~\ref{prop: divergent groups preserving properly convex domains} implies that 
\[\Lambda_{1,d-1}(\Gamma)=\{(x,T_x\partial\Omega):x\in\Lambda_\Omega(\Gamma)\}.\]
If $\Gamma$ is not $P_{1,d-1}$-transverse, then there is some $x,y\in\Lambda_\Omega(\Gamma)$ such that $x\in T_y\partial\Omega$. It follows that $[x,y]_\Omega\subset T_y\partial\Omega$, which implies $[x,y]_\Omega\subset\partial\Omega$. This contradicts the visibility of $\Omega$ and hence $\Gamma$ is $P_{1,d-1}$-transverse.

(2): By part (1), $\Gamma$ is $P_{1,d-1}$-divergent, so Proposition \ref{Prop: Islam-Zimmer} part (3) implies that $\gamma_n\to S\in\Pb(\mathrm{End}(\Rb^d))$ given by $S(\Omega)=x$ and $\Pb(\ker(S))=T_y\partial\Omega$, and $\gamma_n^{-1}\to T\in\Pb(\mathrm{End}(\Rb^d))$ given by $T(\Omega)=y$ and $\Pb(\ker T)=T_x\partial\Omega$. Thus, as a sequence in $\PGL(\Rb^{d*})\simeq\PGL(d,\Rb)$, $\gamma_n\to S^* \in\Pb(\mathrm{End}(\Rb^{d*}))$ with the defining property that $S^*(\Omega^*)=T_x\partial\Omega$ and $\Pb(\ker S^*)=y$.

Since $F^1$ does not lie in $T_y\partial\Omega=\ker S$, 
\[\lim_{n\to\infty}\gamma_n(F^1)=S(F^1)=x.\]
Similarly, since $F^{d-1}$ does not contain $y$,
\[\lim_{n\to\infty}\gamma_n(F^{d-1})=S^*(F^{d-1})=T_x\partial\Omega.\]
This proves (2).

(3): This is immediate from (1) and Propositions \ref{prop: limit set convergence action} and \ref{prop: divergent groups preserving properly convex domains}.
\end{proof}

\subsection{Conical limit points}
Let $\Gamma$ act as a convergence group on $M$. Recall that a point $x \in M$ is a \emph{conical limit point} if there exist an 
infinite sequence $\{\gamma_n\}$ of distinct elements in $\Gamma$ and distinct points $a,b \in M$ such that $\gamma_n(x) \rightarrow a$ and $\gamma_n(y) \rightarrow b$ 
for all $y \in  M \setminus \{x\}$. When $\Gamma\subset\PGL(d,\Kb)$ is a $P_\theta$-transverse subgroup, we denote the set of conical limit points of the $\Gamma$ action on $\Lambda_\theta(\Gamma)$ by $\Lambda_{\theta,c}(\Gamma)$. If $k\in\theta$ and 
$\pi_k:\mathcal F_\theta\to\mathrm{Gr}_k(\mathbb K^d)$ is the  projection map, let 
\[\Lambda_k(\Gamma)=\pi_k(\Lambda_\theta(\Gamma))\ \text{ and }\ \Lambda_{k,c}(\Gamma)=\pi_k(\Lambda_{\theta,c}(\Gamma)).\]

Similarly, when $\Gamma\subset{\rm Aut}(\Omega)$ is projectively visible, we denote the set of conical limit points of the $\Gamma$ action on $\Lambda_\Omega(\Gamma)$ by $\Lambda_{\Omega,c}(\Gamma)$. The points in $\Lambda_{\Omega,c}(\Gamma)$ have a characterization very similar to the classical characterization/definition in hyperbolic space.

\begin{lemma} 
\label{geometric conical def}
Suppose that $\Omega \subset \Pb(\Rb^d)$ is a properly convex domain and
$\Gamma \subset \mathrm{Aut}(\Omega)$ is a projectively visible subgroup.  If $x\in\Pb(\Rb^d)$, then
$x \in \Lambda_{\Omega,c}(\Gamma)$ if and only if  there exist $b_0 \in \Omega$ and a sequence $\{\gamma_n\}$ in $\Gamma$ such that $\gamma_n(b_0) \rightarrow x$ and 
$$
\sup_{n \ge 1} \d_\Omega\left(\gamma_n(b_0), [b_0,x)_\Omega\right) < +\infty.
$$
\end{lemma}

\begin{proof} $(\Leftarrow)$: Suppose that  there exist an 
infinite sequence $\{\gamma_n\}$ of distinct elements in $\Gamma$ and distinct points $a,b \in \Lambda_\Omega(\Gamma)$ such that $\gamma_n(x) \rightarrow a$ and $\gamma_n(y) \rightarrow b$ for all $y \in \Lambda_\Omega(\Gamma) \setminus \{x\}$. Fix some $b_0 \in \Omega$. Then Proposition~\ref{prop:visible_implies_conv} implies that $\gamma_n(b_0) \rightarrow b$ and $\gamma_n^{-1}(b_0) \rightarrow x$ (since $a \neq b$). Since $a \neq b$, the visibility property implies that $(a,b)_\Omega \subset \Omega$. Hence 
$$
 \limsup_{n \rightarrow \infty} \d_\Omega\left( \gamma_n^{-1}(b_0), [b_0, x)_\Omega\right) =  \limsup_{n \rightarrow \infty} \d_\Omega\left( b_0, [\gamma_n(b_0), \gamma_n(x)\right)_\Omega) = \d_\Omega\left( b_0, (b,a)_\Omega\right) < +\infty. 
$$

$(\Rightarrow)$: Suppose that $\gamma_n(b_0) \rightarrow x$ and 
$$
\sup_{n \ge 1} \d_\Omega(\gamma_n(b_0), [b_0,x)_\Omega) < +\infty
$$
for some sequence $\{\gamma_n\}$ in $\Gamma$ and $b_0 \in \Omega$. Pick $\{p_n\}$ in $[b_0, x)$ such that $\{\gamma_n^{-1}(p_n)\}$ is relatively compact in $\Omega$. Passing to a subsequence we can suppose that $\gamma_n^{-1}(p_n) \rightarrow p$, $\gamma_n^{-1}(x) \rightarrow a$ and $\gamma_n^{-1}(b_0) \rightarrow b$. Then $a,b \in \Lambda_\Omega(\Gamma)$ and $p \in (a,b)_\Omega$, so $a \neq b$. Further, Proposition~\ref{prop:visible_implies_conv} implies that $\gamma_n^{-1}(y) \rightarrow b$ for all $y \in \Lambda_\Omega(\Gamma) \setminus \{x\}$. So $x \in \Lambda_{\Omega,c}(\Gamma)$. 

\end{proof} 


\section{Transverse representations}


In this section, we develop the basic theory of transverse representations which will be a crucial tool in our work. 
Our main result is that every $P_\theta$-transverse subgroup is the image of $P_\theta$-transverse representation of a projectively visible subgroup.

\begin{definition}\label{def: transverse rep}
Suppose that  $\theta\subset\{1,\dots, d\}$ is a symmetric subset,  $\Omega \subset \mathbb P(\Rb^{d_0})$ is a properly convex domain and 
$\Gamma \subset \mathrm{Aut}(\Omega)$ is a projectively visible subgroup. A representation 
$\rho : \Gamma \rightarrow \PGL(d,\Kb)$ is \emph{$P_\theta$-transverse} if there exists a continuous embedding 
\begin{align*}
\xi : \Lambda_\Omega(\Gamma) \rightarrow \mathcal F_\theta
\end{align*}
with the following properties: 
\begin{enumerate}
\item $\xi$ is $\rho$-equivariant, i.e. $\rho(\gamma) \circ \xi = \xi \circ \gamma$ for all $\gamma \in \Gamma$, 
\item $\xi(\Lambda_\Omega(\Gamma))$ is a transverse subset of $\mathcal F_\theta$, 
\item if $\{\gamma_n\}$ is a sequence in $\Gamma$ so that $\gamma_n(b_0)\to x\in\Lambda_\Omega(\Gamma)$ and 
$\gamma_n^{-1}(b_0)\to y\in\Lambda_\Omega(\Gamma)$ for some (any) $b_0\in\Omega$, then
$\rho(\gamma_n)(F)\to\xi(x)$ if $F\in\mathcal F_\theta$  is transverse to $\xi(y)$.
\end{enumerate}
We refer to $\xi$ as the \emph{limit map} of $\rho$. 
\end{definition}

It follows from Lemma \ref{lem: KAK 2} that if $\rho:\Gamma\to\PGL(d,\Kb)$ is $P_\theta$-transverse, then it is $P_\theta$-divergent, so
it  has finite kernel and $\rho(\Gamma)$ is a $P_\theta$-transverse subgroup.

By Proposition \ref{prop:visible_implies_conv}, if $\Gamma$ is a projectively visible subgroup of ${\rm Aut}(\Omega)$ for some properly convex domain $\Omega\subset\PGL(d,\Rb)$, 
then the inclusion representation $\Gamma\hookrightarrow\PGL(d,\Rb)$ is $P_{1,d-1}$-transverse, and its limit map is given by $x\mapsto (x,T_x\partial\Omega)$.
If, in addition $\Gamma$ acts cocompactly on the convex hull of $\Lambda_\Omega(\Gamma)$ in $\Omega$, then
$\Gamma$ is hyperbolic (see \cite[Thm.\ 1.15]{DGK} or \cite[Thm.\ 5.1]{Zimmer}) and it follows from 
\cite[Thm.\ 1.3]{GGKW} that $P_\theta$-transverse representations in this case coincide  with $P_\theta$-Anosov representations.
Moreover, if  $\Omega$ is the Klein-Beltrami model of the real hyperbolic 2-plane, $\Gamma$ is finitely generated, and $\theta=\{k,d-k\}$, then $P_\theta$-transverse representations 
coincide with the cusped $P_k$-Anosov representations introduced in \cite{CZZ}.

\begin{theorem}
\label{transverse image of visible1}
If $\Gamma \subset \PGL(d,\Kb)$ is $P_\theta$-transverse, then for any $k\in\theta$ such that $k\le d-k$, there exists a properly convex domain $\Omega \subset \Pb(\Rb^{d_0})$ for some $d_0\in\mathbb N$, 
a projectively visible subgroup $\Gamma_0 \subset \mathrm{Aut}(\Omega)$, and a faithful $P_\theta$-transverse representation
$\rho: \Gamma_0 \rightarrow \PGL(d, \Kb)$ with limit map $\xi:\Lambda_{\Omega}(\Gamma_0)\to\mathcal F_\theta$ so that 
\begin{enumerate}
\item $\rho(\Gamma_0)=\Gamma$.
\item  $\xi(\Lambda_{\Omega}(\Gamma_0)) = \Lambda_{\theta}(\Gamma)$. 
\item $\alpha_1(\kappa(\gamma))=\alpha_k(\kappa(\rho(\gamma)))$ for all $\gamma\in\Gamma_0$.
\end{enumerate} 
\end{theorem}

The main content of the proof of Theorem \ref{transverse image of visible1} is the following two propositions, which are motivated by previous work of Zimmer \cite{Zimmer} 
and Danciger--Gueritaud--Kassel \cite{DGK} in the setting of Anosov representations.

The first proposition provides a representation $\phi$, so that $\phi(\Gamma)$ preserves a properly convex domain.
We say that a map $\xi:\Fc_\theta(\Kb^d)\to\Fc_{\theta'}(\Kb^{d_0})$ 
is \emph{transverse} if it sends every transverse pair of flags in $\Fc_\theta(\Kb^d)$ to a transverse pair of flags in $\Fc_{\theta'}(\Kb^{d_0})$.

\begin{proposition}\label{prop:special irreducible repn} If $1 \le k \le d/2$, then there exists a faithful representation 
\[\phi:\PGL(d,\Kb)\to\PGL(d_0,\Rb)\]
for some $d_0\in\mathbb N$, a $\phi$-equivariant, continuous, transverse map
\[\xi_\phi:\Fc_{k,d-k}(\Kb^d)\to\Fc_{1,d_0-1}(\Rb^{d_0}),\] 
and a properly convex domain $\Omega_0 \subset \Pb(\Rb^{d_0})$ such that:
\begin{enumerate}
\item $\phi(\PGL(d,\Kb))\subset{\rm Aut}(\Omega_0)$.
\item $\alpha_1(\kappa(\phi(g)))=\alpha_k(\kappa(g))$ for all $g \in \PGL(d,\Kb)$.
\item If $\Gamma \subset \PGL(d,\Kb)$ is $P_{k,d-k}$-transverse, then $\phi(\Gamma)$ is $P_{1,d_0-1}$-transverse and $\xi_\phi$ induces a homeomorphism 
$\Lambda_{k,d-k}(\Gamma) \rightarrow \Lambda_{1,d_0-1}(\phi(\Gamma))$.
\end{enumerate}
\end{proposition} 

The second proposition shows that one can enlarge the properly convex domain so that $\phi(\Gamma)$ acts as a projectively visible subgroup.

\begin{proposition}\label{prop:transverse+preserves a domain implies visible} If $\Gamma \subset \PGL(d_0, \Rb)$ is $P_{1,d_0-1}$-transverse and preserves a properly convex
domain $\Omega_0 \subset \Pb(\Rb^{d_0})$, then there exists a properly convex domain $\Omega \subset \Pb(\Rb^{d_0})$, containing $\Omega_0$, 
such that $\Gamma$ preserves $\Omega$ and is a 
projectively visible subgroup of $\Aut(\Omega)$. 
\end{proposition}

Assuming these two propositions, we give the proof of Theorem \ref{transverse image of visible1}.

\begin{proof}[Proof of Theorem \ref{transverse image of visible1}]
For any $k\in\theta$, let $\phi$, $\xi_\phi$, $d_0$ and $\Omega_0$ be given by Proposition~\ref{prop:special irreducible repn} and let
\[\Gamma_0=\rho(\Gamma).\]
Then $\Gamma_0$ is $P_{1,d_0-1}$-transverse and preserves $\Omega_0$. So by Proposition~\ref{prop:transverse+preserves a domain implies visible} there exists a properly 
convex domain $\Omega \subset \Pb(\Rb^{d_0})$ where $\Gamma_0\subset{\rm Aut}(\Omega)$ is a projectively visible subgroup. Since $\phi$ is faithful, we may define
\[
\rho=\phi|_{\Gamma}^{-1}:\Gamma_0\to\PGL(d,\Kb)
\]
which is necessarily a faithful representation. 

By Proposition \ref{prop:visible_implies_conv} part (1), there is a $\Gamma_0$-invariant homeomorphism
\[i:\Lambda_\Omega(\Gamma_0)\to\Lambda_{1,d_0-1}(\Gamma_0).\]
Also, since $\rho$ is $P_\theta$-transverse, there is a  $\Gamma$-invariant homeomorphism
\[i':\Lambda_{k,d-k}(\Gamma)\to\Lambda_\theta(\Gamma).\]
Furthermore, by Proposition~\ref{prop:special irreducible repn} part (3) and the fact that $\xi_\phi$ is $\phi$-equivariant, we may define a $\rho$-equivariant homeomorphism 
\[\bar{\xi}=\xi_\phi|_{\Lambda_{k,d-k}(\Gamma)}^{-1}:\Lambda_{1,d_0-1}(\Gamma_0)\to\Lambda_{k,d-k}(\Gamma).\]
Together, these give a $\rho$-equivariant homeomorphism 
\[\xi=i'\circ \bar{\xi}\circ i:\Lambda_\Omega(\Gamma_0)\to\Lambda_\theta(\Gamma).\]
It is immediate that (1) and (2) hold, and (3) is a consequence of Proposition~\ref{prop:special irreducible repn} part (2).

It remains to show that $\rho$ is a $P_\theta$-transverse representation whose limit map is $\xi$. 
To do so, it suffices to prove condition (3) of Definition \ref{def: transverse rep}; the other conditions are clear. 

Let $\{\gamma_n\}$ be a sequence in $\Gamma_0$ so that $\gamma_n(b_0)\to x\in\Lambda_\Omega(\Gamma_0)$ and 
$\gamma_n^{-1}(b_0)\to y\in\Lambda_\Omega(\Gamma_0)$ for any $b_0\in\Omega$. By taking a subsequence, we may assume that 
\[U_\theta(\rho(\gamma_n))\to F^+\quad \text{and} \quad U_\theta(\rho(\gamma_n)^{-1})\to F^-\] 
for some $F^+,F^-\in\Fc_\theta$. By Proposition \ref{Prop: Islam-Zimmer} part (2), $\gamma_n(z) \to x$ for all  $z\in \Lambda_\Omega(\Gamma_0)-\{y\}$, so 
\[\rho(\gamma_n)(\xi(z))=\xi(\gamma_n(z))\to\xi(x)\]
for all $z\in \Lambda_\Omega(\Gamma_0)-\{y\}$. At the same time, Lemma \ref{lem: KAK} implies that $\rho(\gamma_n)(F)\to F^+$ for all $F\in\Fc_\theta$ transverse to $F^-$. 
Thus, $F^+=\xi(x)$, so condition (3) of Definition \ref{def: transverse rep} holds.
\end{proof}
 
 We now prove Propositions \ref{prop:special irreducible repn} and \ref{prop:transverse+preserves a domain implies visible}.

\begin{proof}[Proof of Proposition~\ref{prop:special irreducible repn}:] Using the notation in Section~\ref{sec: special repns}, let
$$
\phi = S_{V} \circ E^k \quad \text{and} \quad \xi_\phi = \xi_{S_{V}} \circ \xi_{E^k}
$$
where $V = \bigwedge^k \Kb^d$. Then $\phi$ is faithful and $\xi_\phi$ is continuous, transverse, and $\phi$-equivariant. Furthermore, by Equations~\eqref{eqn:singular_values_skew_repn}  and~\eqref{eqn:singular_values_symmetric_repn}  in Section \ref{sec: special repns}, we can identify ${\rm Her}(V)$ with $\Rb^{d_0}$ in such a way that 
\begin{equation}\label{eqn:alpha1 of wedge}
\alpha_1(\kappa(\phi(g))) = \alpha_k(\kappa(g))
\end{equation} 
for all $g \in \PGL(d,\Kb)$. Via this identification, 
\[\phi:\PGL(d,\Kb)\to\PGL(d_0,\Rb)\quad\text{and}\quad\xi_\phi:\Fc_{k,d-k}(\Kb^d)\to\Fc_{1,d_0-1}(\Rb^{d_0}),\]
and condition (2) holds. We observed in Section \ref{sec: symmetric tensors}  that
\[
\Omega_0 = \left\{ [X]\in\Pb({\rm Her}(V)) : X\text{ is positive definite}\right\}
\]
is a $S_V(\PGL(V))$-invariant properly convex domain. By the identification of ${\rm Her}(V)$ with $\Rb^{d_0}$, $\Omega_0\subset\Pb(\Rb^{d_0})$ is a $\phi(\PGL(d,\Kb))$-invariant properly convex domain. Thus, condition (1) holds. 

To verify condition (3), it is enough to prove that $\xi_\phi(\Lambda_{k,d-k}(\Gamma))=\Lambda_{1,d_0-1}(\phi(\Gamma))$. This follows from the following lemma.

\begin{lemma} Suppose $\{g_n\}$ is a sequence in $\PGL(d,\Kb)$. If there exist $F^\pm \in \Fc_{k,d-k}(\Kb^d)$ such that 
\begin{align}
\label{eqn:limiting_behavior}
\lim_{n \rightarrow \infty} g_n(F) = F^+
\end{align}
uniformly on compact subsets of $\{ F \in \Fc_{k,d-k}(\Kb^d) : F \text{ transverse to } F^-\}$, then 
\begin{align*}
\lim_{n \rightarrow \infty} \phi(g_n)(F) = \xi_\phi(F^+)
\end{align*}
uniformly on compact subsets of $\{ F \in \Fc_{1,d_0-1}(\Kb^{d_0}) : F \text{ transverse to } \xi_\phi(F^-)\}$.
\end{lemma} 

\begin{proof} Lemma~\ref{lem: KAK} implies that $\alpha_k(\kappa(g_n)) \rightarrow \infty$. Thus by passing to a tail of $\{g_n\}$, we can assume that $\alpha_k(\kappa(g_n)) > 0$ for all $n$. Then Equation~\eqref{eqn:alpha1 of wedge} implies that $\alpha_1(\kappa(\phi(g))) > 0$ for all $n$.

Let $e_1,\dots, e_d$ denote the standard basis for $\Kb^d$. Using the Cartan decomposition we can write $g_n = k_{1,n} a_n k_{2,n}$ where $k_{1,n}, k_{2,n} \in \mathsf{PO}(d)$ and $a_n\in\mathrm{exp}(\mathfrak{a}^+)$. Let 
$$
F_0^+ = \left( \Span(e_1,\dots, e_k), \Span(e_1,\dots, e_{d-k})\right)
$$
and 
$$
F_0^- = \left( \Span(e_{d-k+1},\dots, e_d), \Span(e_{k+1},\dots, e_{d})\right).
$$ 
Then $U_{k,d-k}(g_n)=k_{1,n}(F_0^+)$ and $U_{k,d-k}(g_n^{-1})=k_{2,n}^{-1}(F_0^-)$. Since $\phi(g_n)=\phi(k_{1,n})\phi(a_n)\phi(k_{2,n})$ is a Cartan decomposition of $\phi(g_n)$, the $\phi$-equivariance of $\xi_\phi$ implies
\[U_{1,d_0-1}(\phi(g_n))=\xi_\phi(k_{1,n}(F_0^+))\quad\text{ and }\quad U_{1,d_0-1}(\phi(g_n^{-1}))=\xi_\phi(k_{2,n}^{-1}(F_0^-)).\]

By Lemma \ref{lem: KAK} and Equation~\eqref{eqn:limiting_behavior},
\begin{align*}
\alpha_k(\kappa(g_n)) \to \infty, \quad k_{1,n}(F_0^+)\to F^+ \quad \text{ and } \quad k_{2,n}^{-1}(F_0^-)\to F^-.
\end{align*}
Then by the continuity of $\xi_\phi$,
\[\alpha_1(\kappa(\phi(g_n))) \to \infty, \quad U_{1,d_0-1}(\phi(g_n))\to\xi_\phi(F^+) \quad \text{and} \quad  U_{1,d_0-1}(\phi(g_n^{-1}))\to \xi_\phi(F^-),\] 
so Lemma~\ref{lem: KAK} implies that
\begin{align*}
\lim_{n \rightarrow \infty} \phi(g_n)(F) = \xi_\phi(F^+)
\end{align*}
uniformly on compact subsets of $\{ F \in \Fc_{1,d_0-1}(\Kb^{d_0}) : F \text{ transverse to } \xi_\phi(F^-)\}$.
\end{proof}

\end{proof}

\begin{proof}[Proof of Proposition~\ref{prop:transverse+preserves a domain implies visible}]
Suppose that $\Gamma \subset \PGL(d_0, \Rb)$ is $P_{1,d_0-1}$-transverse and preserves a properly convex domain $\Omega_0 \subset \Pb(\Rb^{d_0})$.
We will enlarge $\Omega_0$  to a properly convex  domain $\Omega$ so that $\Gamma_0\subset\mathrm{Aut}(\Omega)$ is a projectively visible subgroup of $\mathrm{Aut}(\Omega)$.
If $\Gamma_0$ were irreducible it would suffice to consider the convex domain obtained by intersecting all half-spaces containing $\Omega_0$ and bounded
by hyperplanes in $\pi_{d-1}(\Lambda_{1,d-1}(\Gamma))$. In general, we will construct a properly convex domain $D$ in $\Omega_0^*$ and let $\Omega=D^*$.

Let $B \subset \Omega_0^*$ be an open set whose closure is contained in $\Omega_0^*$ and let $D$ denote the convex hull of $\Gamma (B)$ in $\Omega_0^*$. 
Notice that, by construction, $D$ is a non-empty properly convex domain in $\Pb(\Rb^{d_0*})$ such that $\Gamma\subset{\rm Aut}(D)$. Set $\Omega = D^*$. 
By Lemma \ref{lem: dual domain}, $\Omega\subset\Pb(\Rb^{d_0})$ is a properly convex domain and $\Gamma\subset{\rm Aut}(\Omega)$. 

We will show that $\Gamma$ is a visible subgroup of ${\rm Aut}(\Omega)$.
Since $\Gamma$ is $P_{1,d_0-1}$-transverse, Proposition~\ref{prop: divergent groups preserving properly convex domains} implies that
\begin{equation}
\label{eqn:limit set of Gamma in the construction of Omega} 
\Lambda_\Omega(\Gamma)=\pi_{1}(\Lambda_{1,d_0-1}(\Gamma))=\Lambda_{\Omega_0}(\Gamma)
\end{equation}
Applying Proposition~\ref{prop: divergent groups preserving properly convex domains} to the action of $\Gamma$ on $\Omega_0^*$ shows 
\begin{equation}
\label{eqn:closure of B} 
\overline{\Gamma(B)} = \Gamma(\overline{B}) \cup \pi_{d_0-1}(\Lambda_{1,d_0-1}(\Gamma))
\end{equation}
(where we identify $\Gr_{d_0-1}(\Rb^{d_0}) = \Pb(\Rb^{d_0*})$). 

\begin{lemma} \
\begin{enumerate}
\item If $x,y \in \Lambda_\Omega(\Gamma)$ are distinct points, then $(x,y)_\Omega \subset \Omega$. 
\item If $x\in\Lambda_\Omega(\Gamma)$, then $\partial\Omega$ is $C^1$ at $x$.
\end{enumerate}
\end{lemma}

\begin{proof} 
Let $\mathcal C$ be a component of the cone over $\Omega$ in $\mathbb R^{d_0}\setminus\{0\}$. Since $\overline{D}$ is the convex hull of $\overline{\Gamma(B)}$, Equation~\eqref{eqn:closure of B} implies that if $[f]\in\overline{D}$ and $f|_{\mathcal{C}} > 0$, then we may write
$f = \sum_{j=1}^\ell f_j$ where $f_j|_{\mathcal{C}} >0$ and 
$$
[f_j]\in \overline {\Gamma(B)}=\Gamma(\overline{B})\cup\pi_{d_0-1}(\Lambda_{1,d_0-1}(\Gamma))
$$
for all $j=1,\dots,\ell$. Moreover, if $[f] \in \partial D$, then
$$
[f_j]\in \partial D \cap \overline {\Gamma(B)}=\pi_{d_0-1}(\Lambda_{1,d_0-1}(\Gamma))
$$
for all $j=1,\dots,\ell$. 

(1) Fix $p \in (x,y)_\Omega$. Let $\tilde{x}, \tilde{y} \in \partial \mathcal{C}$ be lifts of $x,y$. Then let $\tilde{p}$ be the lift of $p$ contained in the line segment joining $\tilde x$ to $\tilde y$. If $p$ lies
in $\partial \Omega$, then there exists $f\in (\mathbb{R}^{d_0})^*$ so that $f|_{\mathcal{C}} > 0$, $f(\tilde{p})=0$ and $[f]\in \partial D$. Write
$f = \sum_{j=1}^\ell f_j$ where $f_j|_{\mathcal{C}} >0$ and  
$$
[f_j]\in\pi_{d_0-1}(\Lambda_{1,d_0-1}(\Gamma)).
$$
Fix $j\in\{1,\dots,\ell\}$. Since $f_j|_{\mathcal{C}} > 0$, we have $f_j(\tilde{x})\ge0$ and $f_j(\tilde{y})\ge 0$. Since $\Lambda_{1,d-1}(\Gamma)$ is transverse, either $f_j(\tilde{x})>0$ or $f_j(\tilde{y})>0$, which implies that $f_j(\tilde{p})>0$. Since this is true for all $j$, $f(\tilde{p}) > 0$, and we have a contradiction.
 
(2) By Equation~\eqref{eqn:limit set of Gamma in the construction of Omega}, there exists $H_x$ such that   $(x,H_x) \in \Lambda_{1,d_0-1}(\Gamma)$.
Let $H \in \Gr_{d_0-1}(\Rb^{d_0})$ be a supporting hyperplane at $x$, namely $x \in \Pb(H)$ and $\Pb(H) \cap \Omega = \emptyset$. Since $D=\Omega^*$, there exists $f\in (\mathbb{R}^{d_0})^*$ so that $f|_{\mathcal{C}} > 0$,  $[f]\in \partial D$ and   $\mathrm{ker} f = H$. Write
$f = \sum_{j=1}^\ell f_j$ where $f_j|_{\mathcal C}>0$ and 
$$
[f_j]\in \pi_{d_0-1}(\Lambda_{1,d_0-1}(\Gamma)).
$$
Since $f(x) = 0$, we must have $f_j(x) = 0$ for all $j\in\{1,\dots,\ell\}$. Then, by transversality, $\ker f_j = H_x$ for all $j$. Hence, $H=\ker f = H_x$, 
so $x$ is a $C^1$ point of $\partial \Omega$ with $T_x \partial \Omega = H_x$. 
\end{proof} 

It follows from the above lemma that $\Gamma$ is a projectively  visible subgroup of $\Aut(\Omega)$. \end{proof}


\section{Upper bounds on shadows and Hausdorff dimension} 


Suppose $\Omega \subset \Pb(\Rb^{d_0})$ is a properly convex domain and $\Gamma \subset \mathrm{Aut}(\Omega)$ is a projectively visible subgroup. 
Equip $\Omega$ with its Hilbert metric $\d_\Omega$. Given $b, z \in \Omega$ and $r > 0$, we define the {\em shadow}
\begin{align*}
\mathcal O_r(b,z) \subset\partial \Omega
\end{align*}
to be the set of points $x$ where the geodesic ray $[b,x)_\Omega$ intersects the closed ball $\overline{B_\Omega(z,r)}$ of radius $r$ centered at $z$. Then let
$$\widehat{\mathcal O}_r(b,z)=\mathcal O_r(b,z)\cap \Lambda_\Omega(\Gamma)$$
be the intersection of the shadow with the limit set. Notice that both $\mathcal O_r(b,z)$ and $\widehat{\mathcal O}_r(b,z)$ are closed subsets of $\Pb(\Rb^{d_0})$.
Let $B_{\mathrm{Gr}_k(\mathbb K^d)}(z,t)$ denote the open ball of radius $t>0$ about $z\in\mathrm{Gr}_k(\mathbb K^d)$ with respect to the distance defined in Section~\ref{sec: Lie}. 

Our first shadow lemma gives an upper bound on the diameter of the image of a shadow.

\begin{theorem}\label{thm:shadows I} 
Suppose that $\Omega \subset \Pb(\Rb^{d_0})$ is a properly convex domain, $\Gamma \subset \mathrm{Aut}(\Omega)$ is a projectively visible subgroup  and 
$\rho : \Gamma \rightarrow \PGL(d,\Kb)$ is a $P_\theta$-transverse representation with limit map $\xi: \Lambda_\Omega(\Gamma) \rightarrow \mathcal F_\theta$.
For any $k\in\theta$, $b_0 \in \Omega$ and $r > 0$ there exists $C>1$  so that  if $x \in \Lambda_{\Omega}(\Gamma)$,  
$z \in  [b_0,x)_\Omega$, and $\gamma\in\Gamma$ satisfy $\d_{\Omega}(z,\gamma(b_0)) \le r$, then 
\begin{align*}
\xi^k\left(\widehat{\mathcal O}_r(b_0,z)\right)   \subset B_{\mathrm{Gr}_k(\mathbb K^d)}\left(\xi^k(x),C\frac{\sigma_{k+1}(\rho(\gamma))}{\sigma_k(\rho(\gamma))}\right).
\end{align*}
\end{theorem}

\begin{proof} 
We first prove the theorem assuming that $k=1$.
Assume for contradiction that the theorem does not hold. Then there exist  sequences $\{x_n\}$   and $\{y_n\}$ in $\Lambda_\Omega(\Gamma)$,  $\{z_n\}$ in $\Omega$ 
and $\{\gamma_n\}$ in $\Gamma$ so that for all $n$, we have
 $z_n \in [b_0,x_n)_\Omega$, $\d_{\Omega}(z_n,\gamma_n(b_0)) \le r$, $y_n\in\widehat{\mathcal O}_r(b_0,z_n)$ and 
\begin{align}\label{eqn: lower bound hat}
\d_{\Pb(\Kb^d)}(\xi^1(x_n),\xi^1(y_n)) \ge n \frac{\sigma_{2}(\rho(\gamma_n))}{\sigma_{1}(\rho(\gamma_n))}.
\end{align}

We first observe that $\{\gamma_n\}$ is an escaping sequence in $\Gamma$. If not, then $$ \inf_{n \in\mathbb{Z}^+} \frac{\sigma_{2}(\rho(\gamma_n))}{\sigma_{1}(\rho(\gamma_n))} > 0,$$ 
in which case \eqref{eqn: lower bound hat} implies that $\d_{\Pb(\Kb^d)}(\xi^1(x_n),\xi^1(y_n))\to\infty$, which is impossible. 
It  follows that $\d_\Omega(b_0, \gamma_n(b_0))\to\infty$, and hence that $\d_\Omega(b_0, z_n) \rightarrow \infty$. 

For all $n\in\mathbb N$, choose $w_n \in \overline{B_\Omega(z_n,r)} \cap [b_0,y_n)_\Omega$. Since $\{\gamma_n\}$ is an escaping sequence, we may pass to a subsequence so that 
\[\gamma_n^{-1}(x_n) \rightarrow \bar{x}\in\Lambda_\Omega(\Gamma),\quad \gamma_n^{-1}(y_n) \rightarrow \bar{y}\in\Lambda_\Omega(\Gamma),\quad
 \gamma_n^{-1}(z_n) \rightarrow \bar{z}\in\Omega,\]
 \[\gamma_n^{-1}(w_n) \rightarrow \bar{w}\in\Omega\quad\text{and}\quad \gamma_n^{-1}(b_0) \rightarrow \bar{b}\in\Lambda_\Omega(\Gamma).\]
Note that $\bar{z} \in (\bar{b},\bar{x})_\Omega$ and $\bar{w} \in (\bar{b},\bar{y})_\Omega$, which implies that $\bar{x} \neq \bar{b}$ and $\bar{y}\neq \bar{b}$, see Figure \ref{fig: upper bound}.

 \begin{figure}[ht]
    \centering
    \includegraphics[width=0.8\textwidth]{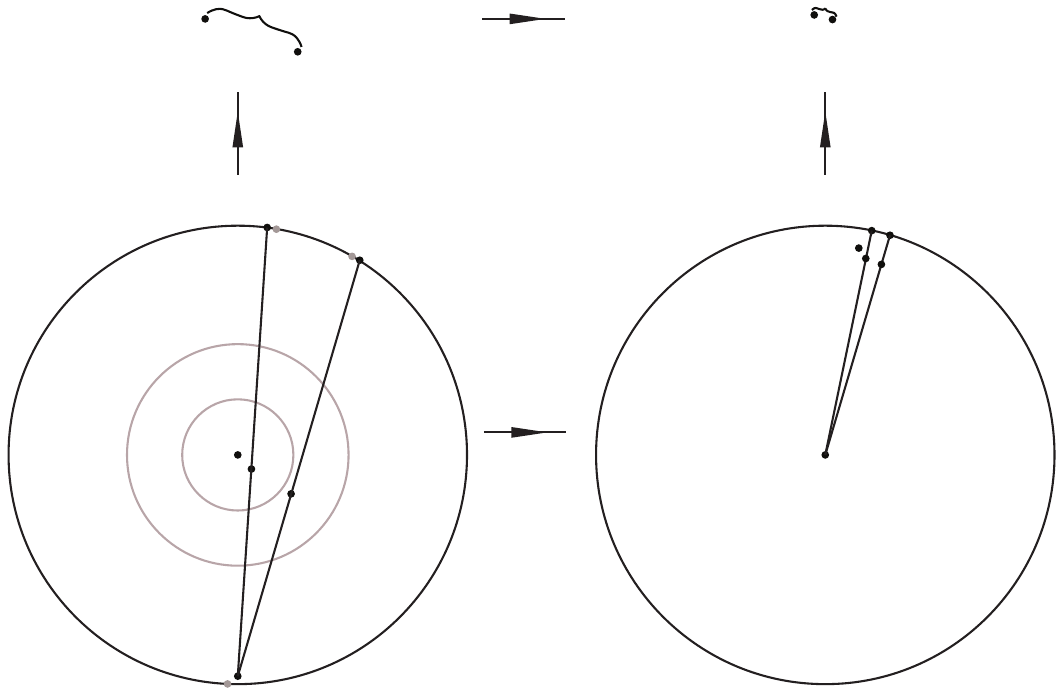}
\put (-190,80){$\gamma_n$}
\put (-197,237){$\rho(\gamma_n)$}
\put (-294,189){$\xi^1$}
\put (-79,189){$\xi^1$}
\tiny
\put (-303,152){$\gamma_n^{-1}(x_n)$}
\put (-308,72){$\gamma_n^{-1}(z_n)$}
\put (-240,149){$\gamma_n^{-1}(y_n)$}
\put (-262,65){$\gamma_n^{-1}(w_n)$}
\put (-280,8){$\gamma_n^{-1}(b_0)$}
\put (-268,160){$\bar{x}$}
\put (-250,145){$\bar{y}$}
\put (-288,-5){$\bar{b}$}
\put (-291, 84){$b_0$}
\put (-85, 75){$b_0$}
\put (-93, 153){$\gamma_n(b_0)$}
\put (-73, 162){$x_n$}
\put (-62, 142){$w_n$}
\put (-77, 146){$z_n$}
\put (-58, 159){$y_n$}
\put (-338,225){$\xi^1(\gamma_n^{-1}(x_n))$}
\put (-262,213){$\xi^1(\gamma_n^{-1}(y_n))$}
\put (-109,227){$\xi^1(x_n)$}
\put (-79, 224){$\xi^1(y_n)$}
\put (-277, 235){$\le 2C$}
\put (-83, 237){$\le 2C\frac{\sigma_2(\rho(\gamma_n))}{\sigma_1(\rho(\gamma_n))}$}
    \caption{Upper bound on the size of shadow.}
    \label{fig: upper bound}
\end{figure}
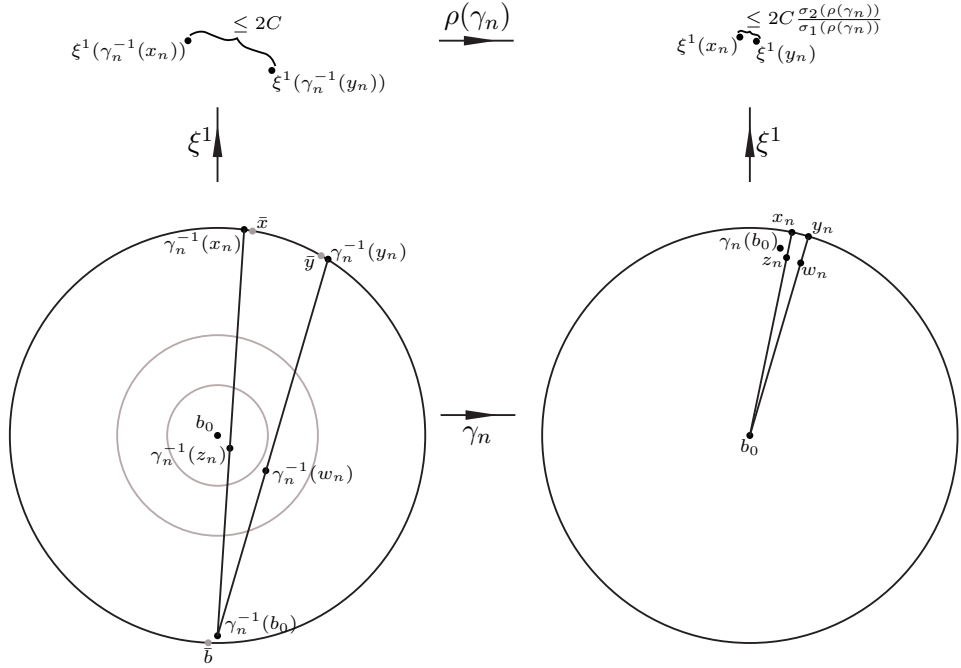

Let $\rho(\gamma_n) = m_n a_n \ell_n$ be a Cartan decomposition of $\rho(\gamma_n)$, where $m_n,\ell_n\in\PU(d,\Kb)$ and $a_n\in\exp(\mathfrak a^+)$. 
For each $n$, let $v_n$ and $u_n$ be unit vectors so that 
\begin{align*}
[v_n]=a_n^{-1}m_n^{-1}(\xi^1(x_n))\ \text{ and } \ [u_n]=a_n^{-1}m_n^{-1}(\xi^1(y_n)).
\end{align*}
Then $\ell_n^{-1}([v_n])=\xi^1(\gamma_n^{-1}(x_n))$ and $\ell_n^{-1}([u_n])=\xi^1(\gamma_n^{-1}(y_n))$. Passing to a subsequence, we can suppose that $\ell_n \rightarrow \ell$, $v_n \rightarrow v$ and 
$u_n \rightarrow u$. Since $\gamma_n^{-1}(b_0) \rightarrow \bar{b}$ and $\rho$ is $P_\theta$-transverse, we may deduce from Lemma  \ref{lem: KAK 2} that
\begin{align*}
\xi^{d-1}(\bar{b}) = \lim_{n\to\infty} U_{d-1}(\rho(\gamma_n^{-1}))=\ell^{-1}(\Span_{\Kb}(e_2,\dots,e_d)).
\end{align*}
Also, by the continuity of $\xi^1$,
\[ \quad \xi^1(\bar{x}) = \ell^{-1} ([  v]) \quad \text{and} \quad \xi^1(\bar{y}) = \ell^{-1} ([ u]).\]

Since $\bar{x} \neq \bar{b}$ and $\bar{y} \neq \bar{b}$,  $\xi^1(\bar{x})$ and $\xi^1(\bar{y})$ are both transverse to $\xi^{d-1}(\bar{b})$, 
or equivalently, $[v]$ and $[u]$ are transverse to $\Span_{\Kb}(e_2,\dots,e_d)$. This implies that there exists $C>0$ so that
\[\tan\angle([v_n],[e_1])\le C\ \text{ and }\ \tan\angle([u_n],[e_1])\le C\] 
for all large enough $n$. Since $\d_{\Pb(\Kb^d)}$ is the angle metric,
\begin{align*}
\d_{\Pb(\Kb^d)}\left( \xi^1(x_n), m_n([e_1]) \right)  =\d_{{\Pb(\Kb^d)}}\left( m_n^{-1}(\xi^1(x_n)), [e_1] \right)
 \le\tan\angle(a_n([v_n]),[e_1])
\le C\frac{\sigma_2(\rho(\gamma_n))}{\sigma_1(\rho(\gamma_n))}.
\end{align*}
Similarly,
\[\d_{\Pb(\Kb^d)}(\xi^1(y_n),m_n([e_1]))\le C\frac{\sigma_2(\rho(\gamma_n))}{\sigma_1(\rho(\gamma_n))},\]
so
\begin{align*}
\d_{\Pb(\Kb^d)}\left( \xi^1(x_n), \xi^1(y_n) \right)  \le 2C\frac{\sigma_2(\rho(\gamma_n))}{\sigma_1(\rho(\gamma_n))},
\end{align*}
which contradicts \eqref{eqn: lower bound hat}. This proves the theorem when $k=1$.

If $k > 1$, consider the exterior power  map $E^k:\PGL(d,\Kb)\to\PGL(D,\Kb)$ and boundary map 
$$
\xi_{E^k}:\Fc_{k,d-k}(\Kb^d) \to\Fc_{1,D-1}(\Kb^D)
$$
defined in Section~\ref{sec: special repns}. Recall that for all $g\in\PGL(d,\Kb)$, 
\[\alpha_k(\kappa(g))=\alpha_1(\kappa(E^k(g))).\]  
Since $\rho$ is a $P_\theta$-transverse representation and $\xi_{E^k}$ is a transverse map,
\[E^k\circ \rho:\Gamma\to\PGL(D,\Kb)\] 
is a $P_{1,D-1}$-tranverse representation with limit map $\xi_{E^k}\circ\xi$.\marginpar{\tiny A: Add Lemma for this?} 
Furthermore, since $\xi_{E^k}$ is 
a smooth embedding, there is some $C'>1$ such that for any $x\in\Lambda_\Omega(\Gamma)$ and $t>0$, 
\[B_{\Pb(\Kb^D)}\left((\xi_{E^k}\circ\xi)^1(x),\frac{t}{C'}\right)\cap \xi_{E^k}^1(\Gr_k(\Kb^d)) \subset \xi_{E^k}^1\left(B_{\Gr_k(\Kb^d)}(\xi^k(x),t)\right)\subset B_{\Pb(\Kb^D)}\left((\xi_{E^k}\circ\xi)^1(x),C't\right).\]
This reduces to the case of $k=1$.
\end{proof}

As a corollary of Theorem \ref{thm:shadows I}, we get the following generalization of the work of Glorieux--Montclair--Tholozan \cite[Thm.\ 4.1]{GMT} and Pozzetti--Sambarino--Wienhard \cite[Prop.\ 4.1]{PSW1}. 

\begin{corollary}\label{prop:upper_bound} 
If $\Gamma \subset \PGL(d,\Kb)$ is $P_{k,d-k }$-transverse, then
\[{\rm dim}_H\left(\Lambda_{k,d-k,c}(\Gamma)\right) \le \delta^{\alpha_k}(\Gamma).\]
In particular, ${\rm dim}_H\left(\Lambda_{k,c}(\Gamma)\right) \le \delta^{\alpha_k}(\Gamma)$.
\end{corollary} 

\begin{proof}
By Theorem \ref{transverse image of visible1}, there is a projectively visible subgroup $\Gamma_0\subset\mathrm{Aut}(\Omega)$ for some properly convex domain $\Omega\subset\Pb(\Rb^{d_0})$, and a $P_\theta$-transverse representation $\rho:\Gamma_0\to\PGL(d,\mathbb K)$ with limit map $\xi$, such that $\rho(\Gamma_0)=\Gamma$ and $\xi(\Lambda_\Omega(\Gamma_0))=\Lambda_{k,d-k}(\Gamma)$ for all $\gamma\in\Gamma_0$. Since $\xi$ is $\rho$-equivariant, injective, and continuous, we have 
\[\Lambda_{k,d-k,c}(\Gamma)=\xi(\Lambda_{\Omega,c}(\Gamma_0)).\]

Fix a base point $b_0\in\Omega$. For $r > 0$, let $L_r\subset \Lambda_{\Omega,c}(\Gamma_0)$ denote the set of conical limit points $x$ with the property that there exists a sequence $\{\gamma_n\}$ in $\Gamma_0$ such that $\gamma_n(b_0)\to x$ and $\d_\Omega(\gamma_n(b_0),[b_0,x)_\Omega)<r$ for all $n$. Then by definition 
\begin{align*}
\Lambda_{\Omega,c}(\Gamma_0) = \bigcup_{r=1}^\infty L_r. 
\end{align*}

Fix, for the moment, $r>0$. For $\gamma \in \Gamma$ define 
\[a(\gamma)=\min\{\alpha_k(\kappa(\rho(\gamma))),\alpha_{d-k}(\kappa(\rho(\gamma)))\}.\]
Then Theorem \ref{thm:shadows I} guarantees that there exists $C=C(r) > 1$ so that 
\begin{align}\label{eqn: mesh size}
{\rm diam}\,  \xi\left(\widehat{\mathcal O}_r\big(b_0, \gamma(b_0)\big)\right)\le Ce^{-a(\gamma)}.
\end{align}
Further, for any $N$,
$$
L_r\subset \bigcup_{\{\gamma\in\Gamma_0\ |\ a(\gamma)>N\} }\widehat{\mathcal O}_r\big(b_0, \gamma(b_0)\big).
$$
By \eqref{eqn: mesh size} the diameter of each element of this covering goes to 0 as $N \rightarrow \infty$ and 
\begin{align*}
\sum_{\{\gamma\in\Gamma_0\ |\ a(\gamma)>N\}} & \left(\mathrm{diam}\,\xi\left(\widehat{\mathcal O}_r\big(b_0, \gamma(b_0)\big)\right)\right)^s
\le C^s\sum_{\gamma\in\Gamma} e^{-sa(\gamma)}\\
&\le C^s\left( \sum_{\gamma\in\Gamma}e^{-s\alpha_k(\kappa(\rho(\gamma)))}+\sum_{\gamma\in\Gamma} e^{-s\alpha_{d-k}(\kappa(\rho(\gamma)))}\right)
\end{align*}
 is finite when $s>\max\{\delta^{\alpha_k}(\Gamma),\delta^{\alpha_{d-k}}(\Gamma)\}$. Note that $\delta^{\alpha_k}(\Gamma)=\delta^{\alpha_{d-k}}(\Gamma)$ because $\alpha_k(\kappa(\rho(\gamma)))=\alpha_{d-k}(\kappa(\rho(\gamma^{-1})))$ for all $\gamma\in\Gamma$. Therefore,
$${\rm dim}_H(\xi(L_r))\le\delta^{\alpha_k}(\Gamma).$$
Since $r$ was arbitrary, we see that
$${\rm dim}_H\left(\Lambda_{k,d-k,c}(\Gamma)\right)={\rm dim}_H\left(\xi(\Lambda_{\Omega,c}(\Gamma_0))\right) =\sup_{r \in\mathbb N} {\rm dim}_H \left( \xi(L_r)  \right) \le \delta^{\alpha_k}(\Gamma).$$
For the second statement, observe that ${\rm dim}_H\left(\Lambda_{k,c}(\Gamma)\right)\le {\rm dim}_H\left(\Lambda_{k,d-k,c}(\Gamma)\right)\le \delta^{\alpha_k}(\Gamma).$
\end{proof}

We use Corollary \ref{prop:upper_bound} and results from the theory of Kleinian groups to prove Theorem \ref{KMO general}.

\begin{theorem}\label{thm:KMO general in body} 
Let $\rho_1:\Gamma\to \mathsf{SO}(d_1-1,1)$ and $\rho_2:\Gamma\to \mathsf{SO}(d_2-1,1)$ be faithful geometrically finite representations so that
$\rho_1(\alpha)$ is parabolic if and only if $\rho_2(\alpha)$ is parabolic. If we regard $\rho=\rho_1\oplus\rho_2$ as a representation into
$\mathsf{PSL}(d_1+d_2,\mathbb R)$, then
$$
{\rm dim}_H\left(\Lambda_{2}(\rho(\Gamma))\right) =\max\left\{{\rm dim}_H\left(\Lambda_{1}(\rho_1(\Gamma))\right), {\rm dim}_H\left(\Lambda_{1}(\rho_2(\Gamma))\right)\right\}.
$$
\end{theorem}

\begin{proof} Let $\mathcal P$ denote the collection of maximal subgroups of $\Gamma$ whose images under each $\rho_i$ consist entirely of parabolic (or trivial) elements
of $\mathsf{SO}(d_i-1,1)$. 
Then $(\Gamma,\mathcal P)$ is relatively hyperbolic and there exists a $\rho_i$-equivariant homeomorphism $\nu_i:\partial(\Gamma,\mathcal P)\to\Lambda_1(\rho_i(\Gamma))$
where $\partial(\Gamma,\mathcal P)$ is the Bowditch boundary of $(\Gamma,\mathcal P)$.
Moreover, $\Lambda_{1,c}(\rho_i(\Gamma))$ is the complement of the images of  fixed points of elements of $\mathcal P$
(See Bowditch \cite[Sec.\ 9]{bowditch-relhyp} and Tukia \cite{tukia-isomorphism}). One may then easily check that
$$\Lambda_2(\rho(\Gamma))=\{\langle \nu_1(x),\nu_2(x)\rangle\ :\ x\in\partial(\Gamma,\mathcal P)\}\subset \mathrm{Gr}_2(\mathbb R^{d_1+d_2}),$$
that $\Lambda_{2,c}(\rho(\Gamma))$ is the complement of the images of fixed points of elements of $\mathcal P$, and that $\rho$ is $P_2$-transverse.
Then, Corollary \ref{prop:upper_bound} implies that 
$${\rm dim}_H\left(\Lambda_{2}(\rho(\Gamma))\right) =
{\rm dim}_H\left(\Lambda_{2,c}(\rho(\Gamma))\right)\le \delta^{\alpha_2}(\rho).$$
Since $\sigma_2(\rho(\gamma))=\min\{\sigma_1(\rho_1(\gamma),\sigma_1(\rho_2(\gamma))\}$ for all $\gamma\in\Gamma$  
$$\#\{\gamma\in\Gamma\ :\ \alpha_2(\rho(\gamma))\le T\} 
\le \#\{\gamma\in\Gamma\ :\ \alpha_1(\rho_1(\gamma))\le T\} +  \#\{\gamma\in\Gamma\ :\ \alpha_1(\rho_2(\gamma))\le T\}$$
which implies that
$$\delta^{\alpha_2}(\rho)\le \max\left\{\delta^{\alpha_1}(\rho_1),\delta^{\alpha_1}(\rho_2)\right\}.$$
Sullivan \cite{sullivan-gf}  showed that ${\rm dim}_H\left(\Lambda_{1}(\rho_i(\Gamma))\right)=\delta^{\alpha_1}(\rho_i)$, so  we may
combine these results to see that 
$${\rm dim}_H\left(\Lambda_{2}(\rho(\Gamma))\right) \le\max\left\{{\rm dim}_H\left(\Lambda_{1}(\rho_1(\Gamma)\right), {\rm dim}_H\left(\Lambda_{1}(\rho_2(\Gamma)\right)\right\}.$$
Since there is a smooth projection of $\Lambda_2(\rho(\Gamma))$ onto $\Lambda_1(\rho_i(\Gamma))$, for both $i=1,2$, we see that equality must hold.
\end{proof}

\begin{remark*} \ \begin{enumerate} 
\item It is possible to give an elementary proof of Corollary \ref{prop:upper_bound} without developing the theory of
projectively visible subgroups. However, the current approach also yields additional structural informations about shadows which may be more
broadly applicable.

\item In the proof of Theorem \ref{KMO general} one may alternatively check that the assumptions of Corollary \ref{prop:upper_bound} hold
by verifying that $\rho$ is $P_2$-relatively dominated in the sense of Zhu \cite{feng}.
\end{enumerate}
\end{remark*}


\section{Singular values and radial projection} 


In this section, we first show that singular values of products of elements in the image of a  $P_\theta$-transverse representation satisfy a coarsely multiplicative lower bound
if the orbit of the elements in the domain  proceed towards infinity ``without backtracking.'' To quantify ``without backtracking'' we use the radial projection map introduced in Section~\ref{sec: properly convex domains}. We will use this result to show that simple root functionals satisfy a coarsely additive lower bound in the absence of backtracking.

Recall that given a properly convex domain $\Omega \subset \mathbb P(\mathbb R^{d_0})$ and $b_0\in\Omega$, the radial projection map based at $b_0 \in \Omega$ is the map 
\[\iota_{b_0}:\overline{\Omega}-\{b_0\}\to\partial\Omega\] 
where $\iota_{b_0}(z) \in \partial \Omega$ is the unique boundary point satisfying $z \in (b_0, \iota_{b_0}(z)]_\Omega$.

\begin{lemma}\label{lem:lower_bd_on_sk} 
Suppose that $\Omega\subset\mathbb P(\mathbb R^{d_0})$ is properly convex, $\Gamma\subset\Aut(\Omega)$ is projectively visible and 
$\rho : \Gamma \rightarrow \PGL(d,\Kb)$ is a $P_\theta$-transverse representation. For any $b_0 \in \Omega$ and $\epsilon> 0$ there exists $C> 0$ so that:  if 
$\gamma, \eta \in \Gamma$,
\begin{align*}
\d_{\Pb(\Rb^{d_0})}\Big(\iota_{b_0}(\gamma^{-1}(b_0)),\iota_{b_0}(\eta(b_0))\Big)\ge \epsilon
\end{align*}
and $k \in \theta$, then 
\begin{align*}
\sigma_k (\rho(\gamma\eta))\ge C\sigma_k(\rho(\gamma))\sigma_k(\rho(\eta)).
\end{align*}
\end{lemma} 

\begin{proof} 
If not, there exist sequences $\{\gamma_n\}$ and $\{\eta_n\}$ in  $\Gamma$  and $k\in\theta$, such that
\begin{align*}
\d_{\Pb(\Rb^{d_0})}\Big(\iota_{b_0}(\gamma_n^{-1}(b_0)),\iota_{b_0}(\eta_n(b_0))\Big) \ge \epsilon
\end{align*}
for all $n \ge 1$ and 
\begin{align*}
\lim_{n \rightarrow \infty} \frac{\sigma_k (\rho(\gamma_n \eta_n))}{\sigma_k(\rho(\gamma_n))\sigma_k(\rho(\eta_n))} = 0.
\end{align*}
Since $\sigma_k(\rho(\gamma_n\eta_n))\ge\max\{\sigma_k(\rho(\gamma_n))\sigma_d(\rho(\eta_n)),\sigma_d(\rho(\gamma_n))\sigma_k(\rho(\eta_n))\}$, this is only possible if $\{\gamma_n\}$ and $\{\eta_n\}$ are escaping sequences. So we can pass to subsequences so that 
$\gamma_n^{-1}(b_0) \rightarrow x \in \Lambda_\Omega(\Gamma)$ 
and $\eta_n(b_0) \rightarrow y \in \Lambda_\Omega(\Gamma)$, in which case $\d_{\Pb(\mathbb R^{d_0})}(x,y)\ge \epsilon$.

Consider Cartan decompositions 
\[\rho(\gamma_n) = m_{n}a_{n} \ell_{n}\ \text{ and }\ \rho(\eta_n) = \hat m_{n} \hat a_{n} \hat\ell_{n},\] 
where $m_n,\ell_n,\hat m_n, \hat\ell_n\in \PU(d,\Kb)$, and $a_n,\hat a_n\in\exp(\mathfrak a^+)$. Let $W=\Span_{\Kb}(e_1,\dots, e_k)$ and $W^\perp=\Span_{\Kb}(e_{k+1},\dots, e_d)$. Then let 
\[\pi_W:\Kb^d\to W\ \text{ and }\ \pi_{W^\perp}:\Kb^d\to W^\perp\] 
denote the orthogonal projections. 

Let $\xi : \Lambda_\Omega(\Gamma) \rightarrow \mathcal F_\theta$ be the limit map of $\rho$. Since $\rho$ is $P_\theta$-transverse, Lemma \ref{lem: KAK 2} implies that 
\[
\lim_{n \rightarrow \infty} \ell_n^{-1} (W^\perp)=\lim_{n \rightarrow \infty} U_{d-k}(\rho(\gamma_n)^{-1}) = \xi^{d-k}(x)
\]
and 
\[
\lim_{n \rightarrow \infty} \hat m_n (W)=\lim_{n \rightarrow \infty} U_k(\rho(\eta_n))= \xi^k(y).
\] 
Since $x\neq y$, $\xi(x)$ is transverse to $\xi(y)$, so there is some $c>0$ such that for large enough $n$, all $u\in W^\perp-\{0\}$, and all $v\in W-\{0\}$, we have
\[ \angle(u, \ell_n \hat m_n(v))= \angle(\ell_n^{-1}(u),  \hat m_n(v))\ge c.\]
Hence, 
there is some $C>0$ such that for large enough $n$ and all $v\in W$, 
\[\norm{\pi_{W}(\ell_n \hat m_n(v))}\ge C\norm{v}.\]
Then for all $v \in W$ we have 
\begin{align*}
\norm{\rho(\gamma_n\eta_n) \hat \ell_n^{-1} (v)} & =  \norm{a_n \ell_n \hat m_n \hat a_n (v)} \ge \norm{\pi_W(a_n \ell_n \hat m_n \hat a_n (v))}   =  \norm{a_n \pi_W(\ell_n \hat m_n \hat a_n (v))}\\
& \ge \sigma_k(\rho(\gamma_n)) \norm{ \pi_W(\ell_n \hat m_n \hat a_n (v))} \ge C  \sigma_k(\rho(\gamma_n))\norm{\hat a_n (v)} \ge C \sigma_k(\rho(\gamma_n))\sigma_k(\rho(\eta_n)) \norm{v}. 
\end{align*}
So by the max-min characterization of singular values 
\begin{align*}
\frac{\sigma_k (\rho(\gamma_n \eta_n))}{\sigma_k(\rho(\gamma_n))\sigma_k(\rho(\eta_n))} \ge \min_{v \in W,\norm{v}=1 }\frac{\norm{\rho(\gamma_n\eta_n) \hat \ell_n^{-1} (v)}}{\sigma_k(\rho(\gamma_n))\sigma_k(\rho(\eta_n))}
\ge C
\end{align*}
and we have a contradiction. 
\end{proof}

As an immediate consequence we see that the first fundamental weight is coarsely additive in the absence of backtracking.

\begin{lemma}\label{lem:lower_bd_on_prod_sk} 
Suppose that $\Omega\subset\mathbb P(\mathbb R^{d_0})$ is properly convex, $\Gamma\subset\Aut(\Omega)$ is projectively visible and 
$\rho : \Gamma \rightarrow \PGL(d,\Kb)$ is a $P_\theta$-transverse representation. For any $b_0 \in \Omega$ and $\epsilon> 0$ there exists $C> 0$ so that:  if 
$\gamma, \eta \in \Gamma$, 
\begin{align*}
\d_{\Pb(\Rb^{d_0})}\Big(\iota_{b_0}(\gamma^{-1}(b_0)),\iota_{b_0}(\eta(b_0))\Big)\ge \epsilon
\end{align*}
and $k \in \theta$, then 
\begin{align*}
(\sigma_1 \cdots \sigma_k) (\rho(\gamma\eta))\ge C(\sigma_1 \cdots \sigma_k) (\rho(\gamma))(\sigma_1 \cdots \sigma_k) (\rho(\eta)).
\end{align*}
\end{lemma} 

\begin{proof} Let $E^k : \PGL(d,\Kb) \rightarrow \PGL(\bigwedge^k \Kb^d)$ be the exterior power representation, see Section~\ref{sec: special repns}. 
Since $k \in \theta$, the representation $E^k \circ \rho$ is $P_{1,d_1-1}$-transverse where  $d_1 = \dim_{\Kb} \bigwedge^k \Kb^d$. 
Further, if we fix the standard inner product on $\bigwedge^k \Kb^d$, then 
$$
(\sigma_1 \cdots \sigma_k)(g) = \sigma_1\left(E^k(g)\right)
$$
for all $g \in \PGL(d,\Kb)$. So Lemma~\ref{lem:lower_bd_on_sk} immediately implies the result. 
\end{proof} 

The proofs of the next two results will use the following well known estimate. 

\begin{observation}\label{obs:trivial_bd_on_products_of_s} If $g,h \in \PGL(d,\Kb)$, then 
$$
(\sigma_1 \cdots \sigma_k)(gh) \le (\sigma_1 \cdots \sigma_k)(g)(\sigma_1 \cdots \sigma_k)(h). 
$$
\end{observation} 

\begin{proof} By definition $\sigma_1(gh) \le \sigma_1(g)\sigma_1(h)$. For $k > 1$,
\[
(\sigma_1 \cdots \sigma_k)(gh) = \sigma_1\left(E^k(gh)\right) \le \sigma_1\left(E^k(g)\right)\sigma_1\left(E^k(h)\right)= (\sigma_1 \cdots \sigma_k)(g)(\sigma_1 \cdots \sigma_k)(h). \qedhere
\]

\end{proof}

We next show that if $k\in\theta$, then the $k^{\rm th}$ simple root of the Cartan projection has a coarsely additive  lower bound if there is no backtracking.

\begin{lemma}\label{lem:lower_bd_on_sk_2} 
Suppose that $\Omega\subset\mathbb P(\mathbb R^{d_0})$ is properly convex, $\Gamma\subset\Aut(\Omega)$ is projectively visible 
and $\rho : \Gamma \rightarrow \PGL(d,\Kb)$ is a $P_\theta$-transverse representation. For any $b_0 \in \Omega$ and $\epsilon> 0$ there exists $C> 0$ so that:  if 
$\gamma, \eta \in \Gamma$,
\begin{align*}
\d_{\Pb(\Rb^{d_0})}\Big(\iota_{b_0}(\gamma^{-1}(b_0)),\iota_{b_0}(\eta(b_0))\Big)\ge \epsilon
\end{align*}
and $k \in \theta$, then 
\begin{align*}
\alpha_k (\kappa(\rho(\gamma\eta)))\ge \alpha_k(\kappa(\rho(\gamma)))+\alpha_k(\kappa(\rho(\eta)))-C.
\end{align*}
\end{lemma} 

\begin{proof} By Lemmas~\ref{lem:lower_bd_on_sk} and~\ref{lem:lower_bd_on_prod_sk}  there exist $C_1,C_2 > 0$, which depend on $\epsilon$, 
but are independent of $\gamma$ and $\eta$, such that 
\begin{equation*}
\sigma_k(\rho(\gamma\eta)) \ge C_1 \sigma_k(\rho(\gamma))\sigma_k(\rho(\eta))
\end{equation*}
and 
\begin{equation*}
(\sigma_1 \cdots \sigma_k) (\rho(\gamma\eta))\ge C_2(\sigma_1 \cdots \sigma_k) (\rho(\gamma))(\sigma_1 \cdots \sigma_k) (\rho(\eta)).
\end{equation*}
Combining these facts with Observation~\ref{obs:trivial_bd_on_products_of_s}, we see that
\begin{align*}
\alpha_k (\kappa(\rho(\gamma\eta)))
&=\log\left( \frac{\sigma_k(\rho(\gamma\eta))}{\sigma_{k+1}(\rho(\gamma\eta))} \right)=\log \left(\frac{(\sigma_1 \cdots \sigma_k)(\rho(\gamma\eta)) }{(\sigma_1 \cdots \sigma_{k+1})(\rho(\gamma\eta)) }\sigma_k(\rho(\gamma\eta))\right) \\
& \ge \log \left( C_1C_2 \frac{(\sigma_1 \cdots \sigma_k)(\rho(\gamma)) (\sigma_1 \cdots \sigma_k)(\rho(\eta))}{ (\sigma_1 \cdots \sigma_{k+1})(\rho(\gamma)) (\sigma_1\cdots \sigma_{k+1})(\rho(\eta))}\sigma_k(\rho(\gamma))\sigma_k(\rho(\eta)) \right)\\
& = \log(C_1C_2)+ \alpha_k(\kappa(\rho(\gamma)))+\alpha_k(\kappa(\rho(\eta))).
\end{align*}
So, Lemma \ref{lem:lower_bd_on_sk_2} holds with $C=-\log(C_1C_2)$.
\end{proof} 

A very similar argument shows that simple roots of the Cartan projection sometimes admit a coarsely additive  upper bound if there is no backtracking.
We will apply this result to the special case of Hitchin representations.

\begin{lemma}\label{lem:upper_bd_on_sk} 
Suppose that $\Omega\subset\mathbb P(\mathbb R^{d_0})$ is properly convex, $\Gamma\subset\Aut(\Omega)$ is projectively visible and  $\rho : \Gamma \rightarrow \PGL(d,\Kb)$ is a $P_\theta$-transverse representation. For any $b_0 \in \Omega$ and $\epsilon> 0$ there exists $C> 0$ so that:  if 
$\gamma, \eta \in \Gamma$,
\begin{align*}
\d_{\Pb(\Rb^{d_0})}\Big(\iota_{b_0}(\gamma^{-1}(b_0)),\iota_{b_0}(\eta(b_0))\Big)\ge \epsilon
\end{align*}
and $k-1, k,k+1\in\theta \cup \{0,d\}$, then
\begin{align*}
\alpha_k (\kappa(\rho(\gamma\eta)))\le \alpha_k(\kappa(\rho(\gamma)))+\alpha_k(\kappa(\rho(\eta)))+C.
\end{align*}
\end{lemma} 

\begin{proof} In the following argument, if $k=1$, then we use the convention that  $(\sigma_1 \cdots \sigma_{k-1})(g) = 1$ for all $g \in \PGL(d,\Kb)$. 

By Lemmas~\ref{lem:lower_bd_on_sk} and~\ref{lem:lower_bd_on_prod_sk}  there exist $C_1 ,C_2 > 0$, which depend only on $\epsilon$, such that 
\begin{equation*}
\sigma_{k+1}(\rho(\gamma\eta)) \ge C_1 \sigma_{k+1}(\rho(\gamma))\sigma_{k+1}(\rho(\eta))
\end{equation*}
(this is obvious when $k+1=d$) and 
\begin{align*}
(\sigma_1 \cdots \sigma_{k-1})(\rho(\gamma\eta)) 
\ge C_2(\sigma_1\cdots \sigma_{k-1})(\rho(\gamma)) (\sigma_1\cdots \sigma_{k-1})(\rho(\eta)).
\end{align*}
Combining these facts with Observation~\ref{obs:trivial_bd_on_products_of_s} we see that
\begin{align*}
\alpha_k (\kappa(\rho(\gamma\eta)))
&=\log\left( \frac{\sigma_k(\rho(\gamma\eta))}{\sigma_{k+1}(\rho(\gamma\eta))} \right)=\log \left(\frac{(\sigma_1 \cdots \sigma_{k})(\rho(\gamma\eta)) }{(\sigma_1 \cdots \sigma_{k-1})(\rho(\gamma\eta)) \sigma_{k+1}(\rho(\gamma\eta))}\right) \\
& \le \log \left(  \frac{\sigma_k(\rho(\gamma)) \sigma_k(\rho(\eta))}{C_1C_2\sigma_{k+1}(\rho(\gamma))\sigma_{k+1}(\rho(\eta))} \right)\\
& = -\log(C_1C_2)+ \alpha_k(\kappa(\rho(\gamma)))+\alpha_k(\kappa(\rho(\eta))),
\end{align*}
which completes the proof. 
\end{proof} 

Finally, we obtain a lower bound which applies when $\gamma(b_0)$ lies within a bounded distance of  the projective line segment joining $b_0$ to $\eta(b_0)$.

\begin{lemma}\label{lem:lower_bd_on_sk_3} 
Suppose that $\Omega\subset\mathbb P(\mathbb R^{d_0})$ is properly convex, $\Gamma\subset\Aut(\Omega)$ is projectively visible and 
$\rho : \Gamma \rightarrow \PGL(d,\Kb)$ is a $P_\theta$-transverse representation. For any $b_0 \in \Omega$ and $r> 0$ there exists $C> 0$ so that:  if $\gamma, \eta \in \Gamma$,
\begin{align*}
\d_\Omega\left( \gamma(b_0), [b_0, \eta(b_0)]_\Omega \right) \le r
\end{align*}
and $k\in\theta$, then
\begin{align*}
\alpha_k (\kappa(\rho(\eta)))\ge \alpha_k(\kappa(\rho(\gamma)))+\alpha_k(\kappa(\rho(\gamma^{-1}\eta)))-C.
\end{align*}
\end{lemma} 

\begin{proof} 
If not, there exist sequences $\{\gamma_n\}$ and $\{\eta_n\}$ in  $\Gamma$  where
\begin{align*}
\d_\Omega\left( \gamma_n(b_0), [b_0, \eta_n(b_0)]_\Omega \right) \le r
\end{align*}
and 
\begin{align*}
\alpha_k (\kappa(\rho(\eta_n)))\le \alpha_k(\kappa(\rho(\gamma_n)))+\alpha_k(\kappa(\rho(\gamma_n^{-1}\eta_n)))-n.
\end{align*}

We claim that $\{\gamma_n^{-1} \eta_n\}$ and $\{\gamma_n\}$ are escaping sequences. Notice that 
\begin{align*}
\alpha_k(\kappa(\rho(\gamma_n)))& \ge n + \alpha_k (\kappa(\rho(\eta_n)))-\alpha_k(\kappa(\rho(\gamma_n^{-1}\eta_n))) \\
&\ge n -  \log \frac{\sigma_1( \rho(\gamma_n))}{\sigma_d( \rho(\gamma_n))}.
\end{align*}
So $\{\gamma_n\}$ must be escaping. A similar argument shows that  $\{\gamma_n^{-1}\eta_n\}$ is escaping. Then
\begin{align*}
\d_\Omega\left( b_0, [\gamma_n^{-1}(b_0), \gamma_n^{-1}\eta_n(b_0)]_\Omega \right) =\d_\Omega\left( \gamma_n(b_0), [b_0, \eta_n(b_0)]_\Omega \right)  \le r,
\end{align*}
so 
$$
\liminf_{n \rightarrow \infty} \d_{\Pb(\Rb^d)}\Big(\iota_{b_0}(\gamma_n^{-1}(b_0)),\iota_{b_0}(\gamma_n^{-1}\eta_n(b_0))\Big) > 0. 
$$
By Lemma~\ref{lem:lower_bd_on_sk_2} there exists some $C > 0$ such that 
\begin{align*}
\alpha_k (\kappa(\rho(\eta_n)))\ge \alpha_k(\kappa(\rho(\gamma_n)))+\alpha_k(\kappa(\rho(\gamma_n^{-1}\eta_n)))-C
\end{align*}
for $n$ sufficiently large. So we have a contradiction. 
\end{proof} 


\section{Lower bounds on shadows}


We will need to restrict to a special class of transverse representations to obtain lower bounds on the inradii of images of shadows. Recall, from the introduction, that a discrete subgroup $\Gamma_0 \subset \PGL(d,\Kb)$ is $(1,1,q)$-hypertransverse if $\Gamma_0$ is $P_\theta$-transverse for some $\theta$ containing $1$ and $q$, and 
$$
F^1+G^1+H^{d-q}
$$
is a direct sum for all pairwise distinct $F,G,H\in\Lambda_\theta(\Gamma_0)$. Then we say that a $P_\theta$-transverse representation $\rho : \Gamma \rightarrow \PGL(d,\Kb)$ is {\em $(1,1,q)$-hypertransverse} if  $1,q \in \theta$ and $\rho(\Gamma)$ is a $(1,1,q)$-hypertransverse subgroup of $\mathsf{PGL}(d,\mathbb K)$. Theorem~\ref{transverse image of visible1} implies that every $(1,1,q)$-hypertransverse subgroup is the image of a $(1,1,q)$-hypertransverse representation. 

We obtain a bound on the inradius of limit sets of hypertransverse representations at uniformly conical limit points, which
generalizes work of Pozzetti-Sambarino-Wienhard \cite{PSW1} from the Anosov setting.
For any $b_0 \in \Omega$ and $R > 0$, let
$$\Lambda_{b_0,R}(\Gamma) \subset \Lambda_\Omega(\Gamma)$$
denote the set of points $x \in \Lambda_\Omega(\Gamma)$ such that the geodesic ray $[b_0,x)_\Omega$ lies in a closed neighborhood of 
radius $R$ of the orbit of $b_0$, i.e. 
$$
[b_0, x)_\Omega \subset \Gamma\left(\overline{B_\Omega(b_0, R)}\right).
$$
We say that such limit points are \emph{$R$-uniformly conical} from $b_0$.

\begin{theorem}\label{thm:size_of_shadows_inside}
Suppose that $\Omega \subset \Pb(\Rb^{d_0})$ is a properly convex domain, $\Gamma \subset \mathrm{Aut}(\Omega)$ is a projectively visible subgroup,
$\rho : \Gamma \rightarrow \PGL(d,\Kb)$ is $(1,1,q)$-hypertransverse  with limit map $\xi:\Lambda_\Omega(\Gamma)\to \mathcal F_{1,q,d-q,d-1}$ and 
\begin{align*}
\sigma_2(\rho(\gamma)) = \sigma_q(\rho(\gamma)) 
\end{align*}
for all $\gamma \in \Gamma$. 
For any $b_0 \in \Omega$ and $r,R > 0$, there exists $C>1$ so that: if $x \in \Lambda_{b_0,R}(\Gamma)$,  
$z \in [b_0,x)_\Omega$ and $\gamma\in\Gamma$ satisfy $\d_{\Omega}(z,\gamma(b_0)) < r$, then 
\begin{align*}
B_{\mathbb P(\mathbb K^d)}\left(x, \frac{1}{C} \frac{\sigma_{2}(\rho(\gamma))}{\sigma_{1}(\rho(\gamma))}\right)\cap \Lambda_1(\rho(\Gamma))\subset\xi^1\left(\widehat{\mathcal O}_r(b_0,z)\right).
\end{align*}
\end{theorem}

\subsection{Projection onto a line} 
One major tool in the proof will be a projection of $\partial \Omega$ onto the projective line through a point $x\in\Lambda_\Omega(\Gamma)$ and
the basepoint $b_0$.

If $x \in \Lambda_\Omega(\Gamma)$, let $x_{\rm opp} \in \partial \Omega\setminus \{x\}$ be the other point of intersection of the projective line through
$b_0$ and $x$ with $\partial\Omega$.
For every $x \in \partial \Omega$, fix a supporting hyperplane $H_x$ to $\Omega$ at $x$. Notice that $H_x$ may intersect $\partial\Omega-\{x\}$,  but it does not contain $x_{\rm opp}$ since $(x,x_{\rm opp}) \subset \Omega$. Moreover, if $x\in\Lambda_\Omega(\Gamma)$, then $H_x=T_x\partial\Omega$ and 
$H_x$ cannot intersect $\Lambda_\Omega(\Gamma)-\{x\}$. Thus, for each $x\in\Lambda_\Omega(\Gamma)$, the codimension $2$ projective subspace
\[W_x=H_x\cap H_{x_{\rm opp}}\] 
does not intersect $\Lambda_\Omega(\Gamma)$, so we may define
\[\pi_x:\Lambda_\Omega(\Gamma)\to [x_{\rm opp}, x]_\Omega\]
by $\pi_x(y)=[x_{\rm opp}, x]_\Omega\cap \left(y\oplus W_x\right)$. Note that $\pi_x(y)=x$ if and only if $y=x$.

\begin{remark} When $\Omega$ is the Klein-Beltrami model of real hyperbolic $d$-space, then $\pi_x(y)$ is the orthogonal projection of $y$ onto the geodesic $(x_{\rm opp},x)_\Omega$.
\end{remark}

The construction of the map $\pi_x$ involves choosing a supporting hyperplane at each boundary point. So, when 
$\partial \Omega$ is not $C^1$, there is no reason to expect that $\pi_x(y)$ or $W_x$ varies continuously with $x$. 
However, the following weak continuity property will suffice for our purposes. 

\begin{lemma}\label{lem: continuity of pi}
Let $\{x_n\}$ and $\{y_n\}$ be sequences in $\Lambda_\Omega(\Gamma)$ such that $x_n\to x$ and $y_n\to y$. Then $x=y$ if and only if $\pi_{x_n}(y_n)\to x$.
\end{lemma}

\begin{proof} Since $\overline{\Omega}$ is compact, it suffices to consider the case when $\lim_{n\to\infty}\pi_{x_n}(y_n)$ exists. 

By taking a subsequence, we may assume that $H_{x_{n,\rm opp}}\to H$ for some supporting hyperplane $H$ to $\Omega$ at $x_{\rm opp}$. Since $x_n\in\Lambda_\Omega(\Gamma)$ for all $n$ and $\Gamma$ is projectively visible, we have 
\[
H_{x_n}=T_{x_n}\partial\Omega\to T_x\partial\Omega= H_{x},
\] 
so $W_{x_n}\to W=H_{x}\cap H$. (It is possible that $H\neq H_{x_{\rm opp}}$ and $W\neq W_{x}$.) 
 
 First, suppose that $x=y$. By definition, $\pi_{x_n}(y_n)\in y_n \oplus W_{x_n}$, which implies
\[
\lim_{n\to\infty}\pi_{x_n}(y_n)\in y\oplus W=x\oplus W=H_{x}.
\] 
Since $\pi_{x_n}(y_n)\in [x_{n,{\rm opp}},x_n]_\Omega$ for all $n$, we have
\[\lim_{n\to\infty}\pi_{x_n}(y_n)\in [x_{\rm opp},x]_\Omega.\] 
It follows that $\pi_{x_n}(y_n)\to H_x\cap[x_{\rm opp},x]_\Omega=x$.

Conversely, suppose that $\pi_{x_n}(y_n)\to x$. Since $y_n\in W_{x_n}\oplus \pi_{x_n}(y_n)$, this implies that $y\in W\oplus x=H_{x}$, which is possible only if $x=y$.
\end{proof}

As a consequence, we prove that if $y$ is close enough to $x$, then $\pi_x(y)\in (b_0,x)$.

\begin{lemma} 
\label{right side}
There exists $\delta > 0$ so that  if $x,y \in \Lambda_\Omega(\Gamma)$ and $0 < \d_{\Pb(\Rb^{d_0})}(x,y) \le \delta$, 
then $\pi_x(y) \in (b_0, x)_\Omega$. 
\end{lemma} 

\begin{proof} 
If not, then there exist sequences $\{x_n\}$ and $\{y_n\}$ in $\Lambda_\Omega(\Gamma)$ such that
\[\pi_{x_n}(y_n) \in [x_{n, {\rm opp}}, b_0]_\Omega \quad \text{and} \quad 0 < \d_{\Pb(\Rb^{d_0})}(x_n,y_n) < 1/n\] 
for all $n$. Passing to a subsequence, we can suppose that 
\[x_n \rightarrow x\in\Lambda_\Omega(\Gamma),\quad y_n \rightarrow y\in\Lambda_\Omega(\Gamma) \quad \text{and} \quad  \pi_{x_n}(y_n) \rightarrow p\in [x_{\rm opp},b_0]_\Omega.\]
Since $\d_{\Pb(\Rb^{d_0})}(x_n,y_n)\to 0$, we have that $x=y$, so Lemma \ref{lem: continuity of pi} implies that $p=x$, which is a contradiction.
\end{proof} 

The following lemma shows that if $z\in [b_0,x)_\Omega$, $z$ is near enough to the orbit $\Gamma(b_0)$ and $\pi_x(y)$ lies between $z$ and $x$ and far enough from $z$, 
then $y$ lies in the shadow of $z$ from $b_0$.

\begin{lemma}
\label{placement} 
Given $r,r'>0$, there exists $T>0$ such that if $x, y \in \Lambda_\Omega(\Gamma)$, $z \in [b_0, x)_\Omega$, 
$\d_\Omega(z, \Gamma(b_0)) \le r'$, $\pi_x(y) \in (z,x)_\Omega$ and $d_\Omega(\pi_{x}(y),z) \ge T$, then $y \in \mathcal O_r(b_0,z)$. 
\end{lemma}

\begin{proof} 
If not, there exist sequences $\{x_n\}$, $\{y_n\}$ in $\Lambda_\Omega(\Gamma)$ and $\{z_n\}$ in $\Omega$ so that 
\begin{align*} z_n & \in [b_0,x_n)_\Omega,  \quad \d_\Omega(z_n, \Gamma(b_0)) \le r', \quad \pi_{x_n}(y_n) \in (z_n,x_n)_\Omega, \\
&  \d_\Omega\left(\pi_{x_n}(y_n),z_n\right) \ge n\quad \text{and}\quad y_n \notin \mathcal O_r(b_0,z_n).
\end{align*}

For all $n\in\mathbb N$, choose $\gamma_n \in \Gamma$ so that $\d_\Omega(\gamma_n^{-1}(z_n),b_0) \le r'$ and pass to a subsequence so that
\[
\gamma_n^{-1}(z_n) \rightarrow \bar{z}\in\Omega, \quad \gamma_n^{-1}(b_0) \rightarrow \bar{b}\in\overline{\Omega}, \quad \gamma_n^{-1}(x_n) \rightarrow \bar{x}\in\Lambda_\Omega(\Gamma),\quad \gamma_n^{-1}(x_{n,{\rm opp}}) \rightarrow \hat{x}\in\partial\Omega, \]
\[\gamma_n^{-1}(y_n) \rightarrow \bar{y}\in\Lambda_\Omega(\Gamma) \quad\text{ and }\quad\gamma_n^{-1}(H_{x_n,{\rm opp}})\to\bar H.
\] 
Let $\bar{W}=\bar H\cap H_{\bar x}$, and note that $\gamma_n^{-1}(W_{x_n})\to\bar{W}$.

Notice that $\bar{z} \in [\bar{b},\bar{x}]_\Omega \cap \Omega$, hence $\bar{b} \neq \bar{x}$, which implies that $\bar{b}\in[\hat{x},\bar{x})_\Omega$. If $\bar{b}=\hat{x}$, then $\hat{x}\in\Lambda_\Omega(\Gamma)$, and the visibility of $\Gamma$ implies that $(\hat{x},\bar{x})_\Omega \subset \Omega$. If $\bar{b}\neq\hat{x}$, then $\bar{b}\in\Omega$, and the convexity of $\Omega$ implies that $(\hat{x},\bar{x})_\Omega \subset \Omega$. In either case, since $\bar H$ is a supporting hyperplane to $\Omega$ at $\hat{x}$, it follows that $\bar x\notin\bar H$, which implies that $\bar{x}\notin \bar{W}$.

Since $\gamma_n^{-1}(\pi_{x_n}(y_n))\in(\gamma_n^{-1}(z_n),\gamma_n^{-1}(x_n))_\Omega$, and 
$$
\lim_{n \rightarrow \infty} \d_\Omega\left(\gamma_n^{-1}(\pi_{x_n}(y_n)),\gamma_n^{-1}(z_n) \right) = \lim_{n \rightarrow \infty} \d_\Omega\left(\pi_{x_n}(y_n),z_n\right)=\infty,
$$
we have $\gamma_n^{-1}(\pi_{x_n}(y_n)) \rightarrow \bar{x}$. Also, since $\bar{x} \notin \bar W$, 
$$
\bar{y} = \lim_{n \rightarrow \infty} \gamma_n^{-1}(y_n) \in \lim_{n \rightarrow \infty} \gamma_n^{-1}(W_{x_n}\oplus\pi_{x_n}(y_n))=\bar{W}\oplus\bar x
= \lim_{n \rightarrow \infty} \gamma_n^{-1}(H_{x_n})= H_{\bar x}.
$$
As such, $\bar{y}=\bar{x}$ because $\bar{x},\bar{y} \in \Lambda_\Omega(\Gamma)$ and $\Gamma$ is a projectively visible subgroup. Then 
\begin{align*}
\lim_{n \rightarrow \infty} \d_\Omega(z_n, (b_0, y_n)_\Omega) =\lim_{n \rightarrow \infty} \d_\Omega(\gamma_n^{-1}(z_n), \gamma_n^{-1}(b_0, y_n)_\Omega) =\d_\Omega(\bar{z}, (\bar{b}, \bar{x})_\Omega) = 0
\end{align*}
and so $y_n \in \mathcal{O}_r(b_0,z_n)$ for $n$ sufficiently large. This is a contradiction.  
\end{proof}

\subsection{Proof of Theorem \ref{thm:size_of_shadows_inside}}
Fix $b_0\in\Omega$ and $r,R>0$.  Let $\delta>0$ be the constant given by Lemma~\ref{right side}.

The following lemma is the crucial estimate in the proof. It shows that if $\gamma(b_0)$ is near $\pi_x(y)$, then one  obtains a lower bound on 
the distance between $\xi^1(x)$ and $\xi^1(y)$ in terms of $\alpha_1(\rho(\gamma))$.

\begin{lemma}\label{lem:lower_bd_on_distance}
There exists $C_1>1$ so that  if $x,y \in \Lambda_\Omega(\Gamma)$, $\gamma\in\Gamma$, $0<\d_{\Pb(\Rb^{d_0})}(x,y) \le \delta$ and 
$\d_{\Omega}(\gamma(b_0),\pi_x(y))\le R$, then
\begin{align*}
\d_{\Pb(\Kb^d)}\left( \xi^1(x), \xi^1(y) \right) \ge \frac{1}{C_1} \frac{\sigma_{2}(\rho(\gamma))}{\sigma_{1}(\rho(\gamma))}.
\end{align*}
\end{lemma} 

\begin{proof} 
If not, there exist sequences $\{x_n\}$ and $\{y_n\}$ in $\Lambda_\Omega(\Gamma)$ and $\{\gamma_n\}$ in $\Gamma$ such that for all $n$, we have
$0 < \d_{\Pb(\Rb^{d_0})}(x_n,y_n) \le \delta$, $\d_{\Omega}(\gamma_n(b_0),\pi_{x_n}(y_n)) \le R$, and
\begin{align}\label{eqn: lower_bd_on_distance}
 \d_{\Pb(\Kb^d)}\left( \xi^1(x_n), \xi^1(y_n) \right) \le \frac{1}{n} \frac{\sigma_{2}(\rho(\gamma_n))}{\sigma_{1}(\rho(\gamma_n))}.
\end{align}
Passing to a subsequence, we may assume that 
\[x_n\to x\in\Lambda_\Omega(\Gamma),\quad y_n\to y\in\Lambda_\Omega(\Gamma),\quad\gamma_n^{-1}(\pi_{x_n}(y_n)) \rightarrow \bar{z}\in\Omega,\]
\[\gamma_n^{-1}(b_0) \rightarrow \bar{b}\in\Omega\cup\Lambda_\Omega(\Gamma),\quad\gamma_n(b_0) \rightarrow b\in\Omega\cup\Lambda_\Omega(\Gamma),\quad \gamma_n^{-1}(x_n) \rightarrow \bar{x}\in\Lambda_\Omega(\Gamma),\]
\[\gamma_n^{-1}(y_n) \rightarrow \bar{y}\in\Lambda_\Omega(\Gamma),\quad\gamma_n^{-1}(x_{n,{\rm opp}})\to \hat x\in\partial\Omega\quad  \text{ and }\ \gamma_n^{-1}(H_{x_{n, {\rm opp}}}) \rightarrow \bar{H}.\]

\begin{figure}[h]
    \centering
    \includegraphics[width=0.7\textwidth]{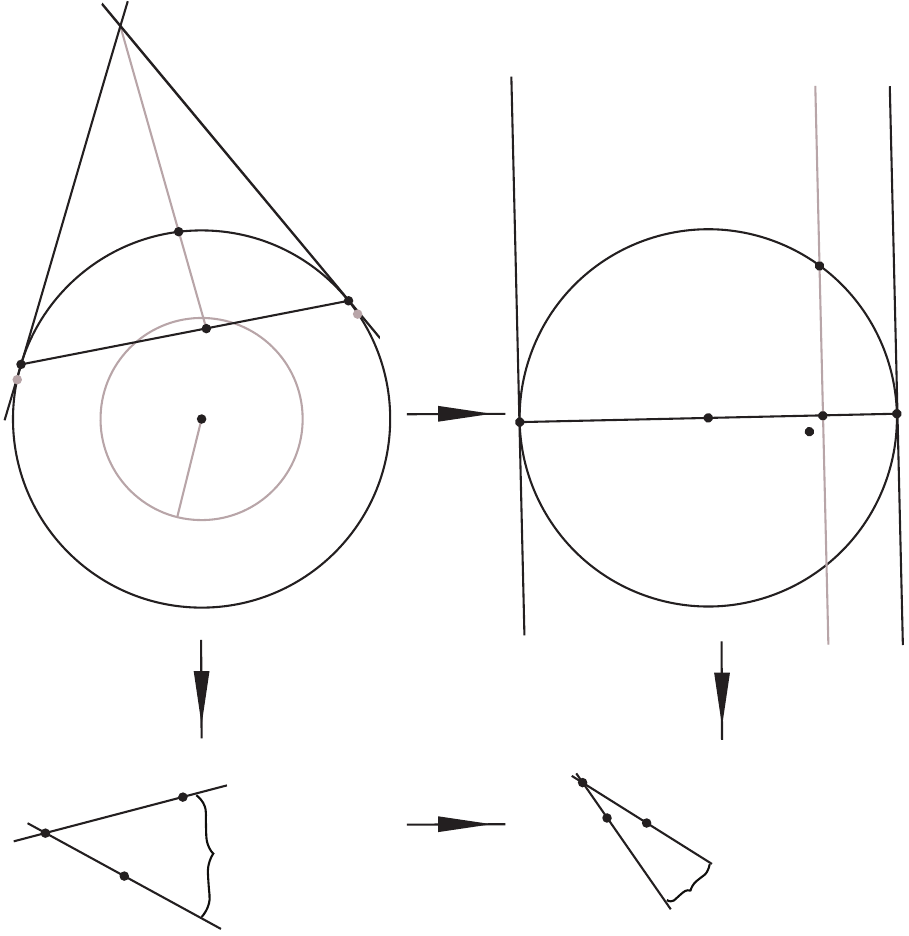}
\put (-171,27){$\rho(\gamma_n)$}
\put (-165,189){$\gamma_n$}
\put (-261,84){$\xi^1$}
\put (-62,84){$\xi^1$}
\tiny
\put (-78,184){$b_0$}
\put (-58,170){$\gamma_n(b_0)$}
\put (-49,186){$\pi_{x_n}(y_n)$}
\put (-30,237){$y_n$}
\put (-2,182){$x_n$}
\put (-134,175){$x_{n,{\rm opp}}$}
\put (-288,217){$\gamma_n^{-1}(\pi_{x_n}(y_n))$}
\put (-198,227){$\gamma_n^{-1}(x_n)$}
\put (-197,211){$\bar{x}$}
\put (-286,249){$\gamma_n^{-1}(y_n)$}
\put (-354,203){$\gamma_n^{-1}(x_{n,{\rm opp}})$}
\put (-311,190){$\hat{x}$}
\put (-248,184){$b_0$}
\put (-252,162){$R$}
\put (-363,28){$\xi^{d-q}(\gamma_n^{-1}(x_{n,{\rm opp}}))$}
\put (-313,13){$\xi^1(\gamma_n^{-1}(x_n))$}
\put (-292,52){$\xi^1(\gamma_n^{-1}(y_n))$}
\put (-243,27){$\ge\frac{1}{C}\abs{ \frac{v_n^j}{v_n^1} -  \frac{w_n^j}{w_n^1}}$}
\put (-115,57){$\xi^{d-q}(x_{n,{\rm opp}})$}
\put (-128,34){$\xi^1(x_n)$}
\put (-91,42){$\xi^1(y_n)$}
\put (-75,10){$\ge\frac{1}{C} \frac{\sigma_j(\rho(\gamma_n))}{\sigma_1(\rho(\gamma_n))}\abs{ \frac{v_n^j}{v_n^1} -  \frac{w_n^j}{w_n^1}}$}
    \caption{Lower bound on the size of shadow.}
    \label{fig: lower bound}
\end{figure}

We first show that $\{\gamma_n\}$ is an escaping sequence, $b=x=y$, $\bar{b}=\hat x$ and $\hat x\ne \bar x$. 
Since $\xi^1$ is injective and continuous, and $\d_{\Pb(\Kb^d)}(\xi^1(x_n),\xi^1(y_n))\to 0$, we know that $x=y$. By Lemma \ref{lem: continuity of pi}, 
we have $\pi_{x_n}(y_n)\to x\in\Lambda_\Omega(\Gamma)$. Then, since $\{\gamma_n^{-1}(\pi_{x_n}(y_n)):n\ge 1\}$
is a relatively compact subset of $\Omega$, it follows that $\{\gamma_n\}$ is an escaping sequence, so $b,\bar{b}\in\Lambda_\Omega(\Gamma)$. 
By Proposition~\ref{Prop: Islam-Zimmer}(2), 
\[b=\lim_{n\to\infty}\gamma_n\big(\gamma_n^{-1}(\pi_{x_n}(y_n))\big)=\lim_{n\to\infty}\pi_{x_n}(y_n)=x.\] 
Since $\Gamma$ acts as a convergence group on $\Lambda_\Omega(\Gamma)$ 
(see Proposition~\ref{prop:visible_implies_conv}(3)),
$\gamma_n^{-1}(b_0)$ converges to $\bar b$ uniformly on compacta in $\Lambda_\Omega(\Gamma)-\{ x\}$. Since $\{x_{n,{\rm opp}}:n\ge 1\}$ is relatively compact in $\Lambda_\Omega(\Gamma)-\{x\}$, 
we have $\hat x=\bar b$. Since $d\left(b_0,\gamma_n^{-1}\left([x_n,x_{n,opp}]_\Omega\right)\right)\le R$, we see that $\hat x\ne\bar x$.

Note that $\bar{H}$ is a supporting hyperplane to $\Omega$ at $\hat x$. Since $\hat x=\bar{b}\in\Lambda_\Omega(\Gamma)$, the supporting hyperplane to $\Omega$ at $\hat x$ is unique, so $\bar{H}=H_{\hat x}$. Thus, if we denote $\bar W=H_{\bar x}\cap H_{\hat x}$, then
$$\bar z=\Big(\bar W\oplus \bar y\Big)\cap [\hat x,\bar x]_\Omega.$$
In particular, $\hat x$, $\bar x$ and $\bar y$ are pairwise distinct points in $\Lambda_\Omega(\Gamma)$ (see Figure \ref{fig: lower bound}), so the transversality of $\rho$ gives 
\begin{align}\label{eqn: transverse}
\xi^1(\bar{x})+\xi^{d-1}(\hat x)=\Kb^d=\xi^1(\bar{y})+\xi^{d-1}(\hat x),
\end{align}
and the assumption that $\rho$ is $(1,1,q)$-hypertransverse implies
\begin{align}\label{eqn: hypertransverse}
 \xi^1(\bar{x}) + \xi^1(\bar{y}) + \xi^{d-q}(\hat x)
 \end{align}
is a direct sum.

For each $n$, let $\rho(\gamma_n) = m_n a_n \ell_n$ be a Cartan decomposition of $\rho(\gamma_n)$, where $m_n,\ell_n\in\PU(d,\Kb)$ and $a_n\in\exp(\mathfrak a^+)$. Also let $v_n=(v_n^1,\dots,v_n^d)$ and $w_n=(w_n^1,\dots,w_n^d)$ be unit vectors so that 
\begin{align*}
[v_n]=a_n^{-1}m_n^{-1}(\xi^1(x_n))\quad \text{and} \quad [w_n]=a_n^{-1}m_n^{-1}(\xi^1(y_n)).
\end{align*}
Then $\ell_n^{-1}([v_n])=\xi^1(\gamma_n^{-1}x_n)$ and $\ell_n^{-1}([w_n])=\xi^1(\gamma_n^{-1}y_n)$. Passing to a subsequence, we can suppose that 
\[\ell_n \rightarrow \ell,\quad v_n \rightarrow v=(v^1,\dots,v^d) \quad \text{and}\quad 
w_n \rightarrow w=(w^1,\dots,w^d).\] 

Since $\rho$ is a $P_{1,q,d-q,d-1}$-transverse representation, and $\{\gamma_n\}$ is an escaping sequence such that $\gamma_n^{-1}(b_0) \rightarrow \hat x$ and $\gamma_n(b_0) \rightarrow x$, Lemma~\ref{lem: KAK 2} implies that
\begin{align}\label{eqn: transverse2}
\xi^{d-i}(\hat{x}) =\lim_{n\to\infty} U_{d-i}(\rho(\gamma_n^{-1}))= \ell^{-1}(\Span_{\Kb}(e_{i+1}, \dots, e_d))
\end{align}
for $i=1$ and $i=q$.
Also, since $\xi^1$ is continuous, 
\[ \quad \xi^1(\bar{x}) = \ell^{-1} (\left[  v\right])\quad\text{and}\quad \xi^1(\bar{y}) = \ell^{-1} (\left[  w \right]).\]
Hence,  Equations \eqref{eqn: transverse} and \eqref{eqn: transverse2} (with $i=1$) imply that $v^1\ne 0\ne w^1$, so $v_n^1\ne 0\ne w_n^1$ for sufficiently large $n$.
Equations \eqref{eqn: hypertransverse} and \eqref{eqn: transverse2} (with $i=q$) imply that the 
collection of vectors $\{v,w,e_{q+1},\ldots,e_d\}$ are linearly independent over $\Kb$. In particular, there is some $j\in\{2,\dots,q\}$ such that 
$$\frac{v^j}{v^1}-\frac{w^j}{w^1}\ne 0$$

Consider the affine chart 
\[
\mathbb{A}^{d-1}=\left\{ \left[ (u^1, \dots, u^d) \right] \in \Pb(\Kb^d) : u^1\neq 0\right\}
\] 
of $\Pb(\Kb^d)$. Let $\d_{\mathbb{A}^{d-1}}$ denote the pullback to $\mathbb{A}^{d-1}$ of the standard metric on $\Kb^{d-1}$ via the identification $\Ab^{d-1}\simeq\Kb^{d-1}$ given by
\[
\left[ (u^1, \dots, u^d) \right]  \mapsto \left(\frac{u^2}{u^1},\dots,\frac{u^d}{u^1}\right).
\]
Since $v_n^1\ne 0\ne w_n^1$ for sufficiently large $n$ and $\alpha_1(\kappa(\rho(\gamma_n)))\to\infty$, we have
\[\lim_{n\to\infty}a_n([v_n])=[e_1]= \lim_{n\to\infty}a_n([w_n]),\]
so there is a compact $\Kc\subset\mathbb{A}^{d-1}$ that contains $a_n([v_n])$ and $a_n([w_n])$ for sufficiently large $n$. Choose $C>1$ so that $\d_{\Pb(\Kb^{d})}$ and 
$\d_{\mathbb{A}^{d-1}}$ are $C$-bilipschitz on $\Kc$. Then
\begin{align*}
\d_{\Pb(\Kb^{d})}\left( \xi^1(x_n), \xi^1(y_n) \right) & =\d_{\Pb(\Kb^{d})}\left( \rho(\gamma_n)\xi^1(\gamma_n^{-1}(x_n)), \rho(\gamma_n)\xi^1(\gamma_n^{-1}(y_n)) \right)  = \d_{\Pb(\Kb^{d})}\left( a_n([v_n]), a_n ([w_n]) \right) \\
& \ge\frac{1}{C} \d_{\mathbb{A}^{d-1}}\left( a_n([v_n]), a_n ([w_n]) \right)  =\frac{1}{C} \sqrt{\sum_{i=2}^d \frac{\sigma_i(\rho(\gamma_n))^2}{\sigma_1(\rho(\gamma_n))^2}\abs{ \frac{v_n^i}{v_n^1} -  \frac{w_n^i}{w_n^1}}^2}\\
&\ge\frac{1}{C} \frac{\sigma_j(\rho(\gamma_n))}{\sigma_1(\rho(\gamma_n))}\abs{ \frac{v_n^j}{v_n^1} -  \frac{w_n^j}{w_n^1}}.
\end{align*}
Thus, by \eqref{eqn: lower_bd_on_distance} and the hypothesis that $\sigma_q(\rho(\gamma_n))=\sigma_2(\rho(\gamma_n))$,
\begin{align*}
0 &= \lim_{n \rightarrow \infty} \frac{\sigma_1(\rho(\gamma_n))}{\sigma_2(\rho(\gamma_n))} \d_{\Pb(\Kb^{d})}\left( \xi^1(x_n), \xi^1(y_n) \right)  \ge \frac{1}{C} \lim_{n \rightarrow \infty}\frac{\sigma_1(\rho(\gamma_n))}{\sigma_2(\rho(\gamma_n))} \left(\abs{ \frac{v_n^j}{v_n^1} -  \frac{w_n^j}{w_n^1}}  \frac{\sigma_j(\rho(\gamma_n))}{\sigma_1(\rho(\gamma_n))}\right)\\
& = \frac{1}{C}\abs{\frac{v^j}{v^1}-\frac{w^j}{w^1}} \neq 0 
\end{align*}
and we have  achieved a contradiction.
\end{proof} 

If Theorem \ref{thm:size_of_shadows_inside} does not hold, then for every $n$ there exist $x_n \in \Lambda_{b_0,R}(\Gamma)$, $y_n\in\Lambda_\Omega(\Gamma)$, 
$z_n \in [b_0,x_n)_\Omega$  and $\gamma_n \in \Gamma$ such that
$$\d_{\Omega}(z_n,\gamma_n(b_0)) \le r,\quad   y_n\notin \mathcal O_r(b_0,z_n)\quad \text{and} \quad 
\d_{\Pb(\Kb^d)}(\xi^1(x_n),\xi^1(y_n))\le\frac{1}{n} \frac{\sigma_{2}(\rho(\gamma_n))}{\sigma_{1}(\rho(\gamma_n))}.$$
In particular, $x_n\neq y_n$ for all $n$. Also, by taking subsequences, we may assume that 
\[x_n\to x\in\Lambda_\Omega(\Gamma),\quad y_n\to y\in\Lambda_\Omega(\Gamma)\quad \text{and}\quad  \gamma_n(b_0)\to b\in\Omega\cup\Lambda_\Omega(\Gamma).\]

Since $\xi$ is a continuous embedding, $x=y$. So by passing to a tail of our sequences, 
we can assume that $0<\d_{\Pb(\Rb^{d_0})}(x_n,y_n) \le \delta$ for all $n$. Also, since $x_n \in \Lambda_{b_0,R}(\Gamma)$,
there exists $\beta_n \in \Gamma$ with $\d_{\Omega}(\pi_{x_n}(y_n),\beta_n(b_0)) \le R$. 
Then, by Lemma~\ref{lem:lower_bd_on_distance} (for the first inequality) and assumption (for the second), there exists $C_1>1$ that does not depend on $n$, so that
\begin{align}\label{eqn: upper and lower bounds 1}
\frac{\sigma_{2}(\rho(\beta_n))}{\sigma_{1}(\rho(\beta_n))} \le C_1 \d_{\Pb(\Kb^d)}\left( \xi^1(x_n), \xi^1(y_n) \right) \le \frac{C_1}{n}\frac{\sigma_{2}(\rho(\gamma_n))}{\sigma_{1}(\rho(\gamma_n))}.
\end{align}

Notice that $\sigma_j( \rho(\beta_n)) = \sigma_j(\rho(\gamma_n) \rho(\gamma_n^{-1} \beta_n))$ and 
\begin{align*}
\d_\Omega( \gamma_n^{-1} \beta_n(b_0), b_0) =\d_\Omega( \beta_n(b_0), \gamma_n(b_0)) \le \d_{\Omega}(\pi_{x_n}(y_n), z_n)+R+r.
\end{align*}
So if we set
\begin{align*}
S_n = \max\left\{ \sigma_1(\rho(\eta)) : \d_\Omega( \eta(b_0), b_0) \le d_{\Omega}(\pi_{x_n}(y_n), z_n)+R+r\right\},
\end{align*}
then for all $j\in\{1,\dots,d\}$, we have
\begin{align*}
\frac{1}{S_n} \sigma_j( \rho(\gamma_n)) \le \sigma_j(\rho(\beta_n)) \le S_n  \sigma_j( \rho(\gamma_n)).
\end{align*}
Thus, \eqref{eqn: upper and lower bounds 1} implies that $S_n\ge\sqrt{\frac{n}{C_1}}\to\infty$, so
\begin{align}\label{eqn: upper and lower bounds 2}
\lim_{n \rightarrow \infty} \d_{\Omega}(\pi_{x_n}(y_n), z_n) = \infty.
\end{align}

Since $y_n \notin \mathcal O_r(b_0,z_n)$ and $\d_{\Omega}(\pi_{x_n}(y_n), z_n) \to \infty$, Lemma \ref{placement} (with $r'=r$) implies that by taking the tail end of the sequence, 
we may assume that $\pi_{x_n}(y_n)\notin (z_n,x_n)_\Omega$ for all $n$. At the same time, since $0<\d_{\Pb(\Rb^{d_0})}(x_n,y_n) \le \delta$ for all $n$, 
Lemma \ref{right side} implies that $\pi_{x_n}(y_n)\in (b_0,x_n)_\Omega$. Thus $\pi_{x_n}(y_n)\in [b_0,z_n]_\Omega$ for all $n$.

Equation~\eqref{eqn:Hilbert metric Hausdorff distance} implies that
$$\d_\Omega^{\rm Haus}\Big([b_0,z_n]_\Omega,[b_0,\gamma_n(b_0)]_\Omega\Big)\le \d_\Omega(z_n,\gamma_n(b_0))\le r.$$
Then, since $\pi_{x_n}(y_n)\in [b_0,z_n]_\Omega$ and $d(\beta_n(b_0),\pi_{x_n}(y_n))\le R$, we see that
$$\d_\Omega\left( \beta_n(b_0), [b_0, \gamma_n(b_0)]_\Omega\right) \le R+r.
$$
So by Lemma~\ref{lem:lower_bd_on_sk_3} there exists $C > 0$, depending only on $R+r$, such that 
$$
\alpha_1(\kappa(\rho(\gamma_n))) \ge  \alpha_1(\kappa(\rho(\beta_n)))+\alpha_1(\kappa(\rho(\beta_n^{-1}\gamma_n))) -C\ge  \alpha_1(\kappa(\rho(\beta_n)))-C
$$
for all $n$. But this contradicts Equation~\eqref{eqn: upper and lower bounds 1}.
\qed


\section{Lower bounds on Hausdorff dimension}


In this section we complete the proof of Theorem~\ref{thm:BishopJones} which we restate here.

\begin{theorem}\label{thm:lower_bd} 
Suppose that $\Gamma\subset\PGL(d,\Kb)$ is  $(1,1,q)$-hypertransverse and $\sigma_2(\gamma) = \sigma_q(\gamma)$ for all $\gamma \in \Gamma$. Then 
\begin{align*}
{\rm dim}_H(\Lambda_{1,c}(\Gamma)) = \delta^{\alpha_1}(\Gamma).
\end{align*} 
\end{theorem}

As mentioned in the introduction, Theorem~\ref{thm:lower_bd} generalizes earlier results of Pozzetti--Sambarino--Wienhard \cite{PSW1} for $(1,1,2)$-hyperconvex Anosov representations and  Bishop--Jones \cite{bishop-jones} for Kleinian groups. 

\subsection{Proof of Theorem \ref{thm:lower_bd}}
Let $\theta=\{1,q,d-q,d-1\}$. By Theorem \ref{transverse image of visible1}, there exist a properly convex domain $\Omega\subset\Pb(\Rb^{d_0})$, a projectively visible subgroup $\Gamma_0\subset\mathrm{Aut}(\Omega)$ 
 and a $P_\theta$-transverse representation $\rho:\Gamma_0\to\PGL(d,\mathbb K)$ with limit map 
$\xi:\Lambda_\Omega(\Gamma_0)\to\Fc_\theta$, such that $\rho(\Gamma_0)=\Gamma$ and $\xi(\Lambda_{\Omega}(\Gamma_0)) = \Lambda_{\theta}(\Gamma)$. Then, by definition, $\rho$ is a $(1,1,q)$-hypertransverse representation. 

Choose $b_0 \in \Omega$ so that ${\rm Stab}_{\Gamma_0}(b_0)=\id$. Then  the orbit map $\gamma \mapsto \gamma(b_0) \in\Omega$ is injective and 
$\Gamma_0(b_0)$ is a closed discrete subset of $\Omega$. Further, the function $c : \Gamma_0(b_0) \rightarrow (0,1]$ given by 
\begin{align*}
c(\gamma(b_0)) := \frac{\sigma_2(\rho(\gamma))}{\sigma_1(\rho(\gamma))}
\end{align*}
is well defined. 

The following technical result places us in the situation where we may apply the argument of Bishop--Jones. 
 
\begin{proposition}\label{prop:finding_children} For any $0 < \delta < \delta^{\alpha_1}(\Gamma)$, there exist $r_0, D_0 > 0$ such that if $z \in \Gamma_0(b_0)$, 
then there exists a finite subset $\mathcal C(z) $ of $\Gamma_0(b_0)-\{z\}$ 
with the following properties
\begin{enumerate}
\item if $w\in \mathcal C(z)$, then $\mathcal O_{2r_0}(b_0,w) \subset \mathcal O_{r_0}(b_0,z)$,
\item the sets $\left\{\mathcal O_{2r_0}(b_0,w)\right\}_{w\in \mathcal C(z)}$ are pairwise disjoint, 
\item  if $w\in \mathcal C(z)$, then $\d_{\Omega}(z,w) \le D_0$, 
\item  $\sum_{w\in \mathcal C(z)}  c(w)^{\delta} \ge c(z)^{\delta}$.
\end{enumerate}
\end{proposition}

Assuming Proposition~\ref{prop:finding_children} for the moment we prove Theorem~\ref{thm:lower_bd}. 

\medskip\noindent
{\bf Outline of proof:}
We will construct, for each $0 < \delta < \delta^{\alpha_1}(\Gamma)$, a set $E_\delta\subset\Lambda_{\Omega,c}(\Gamma_0)$, 
a measure $\mu_\delta$ on $\partial \Omega$  supported on $E_\delta$ and constants $C,t_0 > 0$ such that:  
$$\xi^1_*\mu_\delta\left(B_{\Pb(\Kb^d)}(p,t)\right)\le C t^\delta$$
for all $p\in\Pb(\Kb^d)$ and $0 < t < t_0$. 
Once we do so, then we may apply simple covering arguments in the spirit of Frostman's Lemma, see for instance~\cite[Thm.\ 1.2.8 and Lem.\ 3.1.1]{bishop-peres},
to show that 
\[\dim_{H}(\Lambda_{1,c}(\Gamma))\ge\dim_{H}(\xi^1(E_\delta))\ge\delta.\]
Taking the limits as $\delta\to \delta^{\alpha_1}(\Gamma)$ will yield that ${\rm dim}_H(\Lambda_{1,c}(\Gamma)) \ge \delta^{\alpha_1}(\Gamma)$.
Since Corollary \ref{prop:upper_bound} implies that 
${\rm dim}_H(\Lambda_{1,c}(\Gamma)) \le \delta^{\alpha_1}(\Gamma)$, our theorem will follow. 

\medskip

Fix $0 < \delta < \delta^{\alpha_1}(\Gamma)$, and let $r_0,D_0>0$ be as given by Proposition~\ref{prop:finding_children}. 

We inductively construct a tree $\mathcal T=\mathcal T_\delta\subset\Omega$ with root $b_0$,  
whose vertices are a subset of $\Gamma_0(b_0)$,
and with the property that for every vertex $z \in \mathcal T$, its children $\mathcal C(z)$ satisfy the conditions in Proposition~\ref{prop:finding_children}. 
(Notice that, by definition, $\mathcal O_{r_0}(b_0,z)=\partial\Omega$ if $d_\Omega(b_0,z)<r_0$.)
Properties (1) and (2) in Proposition~\ref{prop:finding_children} guarantee that our inductive construction does indeed produce a tree. 

Let $E=E_\delta \subset \Lambda_\Omega(\Gamma_0)$ be the set of accumulation points of the vertices of $\mathcal T$. 
Since $\Gamma_0(b_0)$ is discrete and $\mathcal T$ is infinite, $E$ is non-empty.  We now observe that $E$ is uniformly conical.

\begin{lemma} \label{lem: bishopjones1}\
\begin{enumerate}
\item For any $x \in E$ there exists a (discrete) geodesic ray $\{x_n\}$  in $\mathcal T$ such that $x_0 = b_0$, $x_n\to x$ and
$$x\in\bigcap_{n\in\mathbb N} \mathcal O_{r_0}(b_0,x_n).$$
\item $E \subset \Lambda_{b_0,R_0}(\Gamma_0)\subset \Lambda_{\Omega,c}(\Gamma_0)$ where $R_0 = \frac{1}{2}D_0 + 2r_0$. 
\end{enumerate}
\end{lemma}

\begin{proof} 
(1): Fix $x \in E$. Then there exist a sequence of vertices $\{w_m \}$ in $\mathcal T$ with 
$w_m \rightarrow x$. 
We may pass to a subsequence so that $w_m\to w\in\partial_\infty\mathcal T$, the abstract visual boundary of $\mathcal T$.
Let $\sigma: \mathbb Z_{\ge 0} \rightarrow \mathcal T$ be a geodesic ray that starts at $b_0$ and limits to $w$. Let $x_n = \sigma(n)$ and let 
$$k_m=\max\{n\ :\ w_m \ \text{ is a descendent of }\ x_n\}.$$
Observe that $k_m\to\infty$, so we may pass to a subsequence so that the sequence $\{k_m\}$ is strictly increasing.

Since $x_{n+1}$ is a child of $x_n$ for all $n$ and $w_m$ is a descendent of $x_{k_m}$ for all $m$, Proposition~\ref{prop:finding_children} part (1) gives
\begin{align*}
\mathcal O_{2r_0}(b_0,x_{n+1}) \subset \mathcal O_{r_0}(b_0,x_n)\ \text{ and }\ \mathcal O_{2r_0}(b_0,w_m)\subset\mathcal O_{r_0}(b_0,x_{k_m})
\end{align*}
for all $n$ and $m$. Hence,
$$\iota_{b_0}(w_m)\in\mathcal O_{r_0}(b_0,x_{k_m})=\bigcap_{n\le k_m}\mathcal O_{r_0}(b_0,x_n),$$
where $\iota_{b_0}$ is the radial projection. Since $\{k_m\}$ is strictly increasing and the shadows are closed, 
\begin{align}\label{eqn: nested intersection to x}
x=\lim_{m\to\infty} w_m=\lim_{m\to\infty} \iota_{b_0}(w_m)\in\bigcap_{n\in\mathbb N} \mathcal O_{r_0}(b_0,x_n).
\end{align}

Finally, we prove that $x_n \rightarrow x$. By \eqref{eqn: nested intersection to x}, there exists a sequence $\{y_n\}$ along the geodesic ray $[b_0,x)_\Omega$ in $\Omega$ such 
that $\d_{\Omega}(y_n, x_n) \le r_0$ for all $n$. Suppose for contradiction that the sequence $\{x_n\}$ does not converge to $x$. Then there exists a subsequence $\{x_{n_j}\}$
with $x_{n_j} \to y\in E$ for some $y \neq x$. Let $\{\gamma_j\}$ be the sequence in $\Gamma_0$ so that $\gamma_j(b_0) = x_{n_j}$ for all $j$. 
Then Proposition~\ref{Prop: Islam-Zimmer} part (2) implies that $\gamma_j(b) \rightarrow y$ for all $b \in \Omega$ 
and the convergence is locally uniform. Since $\{\gamma_j^{-1}(y_{n_j}) \}$ is relatively compact in $\Omega$, 
we have $y_{n_j} = \gamma_n(\gamma_n^{-1}(y_n)) \rightarrow y$. However, by construction $y_n\in  [b_0, x)$ for all $n$, so this is not possible.

(2): Fix $x\in E$. By (1), there is a geodesic ray $\{x_n\}$ in $\mathcal T$ such that $x_n\to x$, and there exists a sequence $\{y_n\}$ along the geodesic ray 
$[b_0,x)_\Omega$ in $\Omega$ such that $\d_{\Omega}(y_n, x_n) \le r_0$ for all $n$. Then by Proposition~\ref{prop:finding_children} part (3),
\begin{align*}
\d_\Omega(y_n, y_{n+1}) \le \d_\Omega(x_n, x_{n+1}) + 2r_0\le D_0+2r_0. 
\end{align*}
It then follows that
\begin{align*}
[b_0,x)_\Omega \subset \bigcup_{n=1}^\infty \overline{B_\Omega\left(  y_n, \frac{1}{2}D_0 + r_0\right)} \subset \bigcup_{n=1}^\infty \overline{B_\Omega\left( x_n, \frac{1}{2}D_0 + 2r_0\right)}.
\end{align*}
Since $x_n\in\Gamma_0(b_0)$ for all $n$, we have $x \in \Lambda_{b_0,R_0}(\Gamma_0)$, where $R_0 = \frac{1}{2}D_0 + 2r_0$.
Therefore, $E\subset  \Lambda_{b_0,R_0}(\Gamma_0)$.
Lemma \ref{geometric conical def} implies that $ \Lambda_{b_0,R_0}(\Gamma_0)\subset\Lambda_{\Omega,c}(\Gamma_0)$, which completes the proof.
\end{proof} 

We now construct a well-behaved probability measure $\mu$ on $E$. 

\begin{lemma} \label{lem: measure}
There exists a Borel probability measure $\mu=\mu_\delta$ on $\partial\Omega$ such that $\mu(E) = 1$ and
\begin{enumerate}
\item
$\mu\left( \mathcal O_{2r_0}(b_0,z)\right) \le c(z)^{\delta}$ for every vertex $z$ of $\mathcal T$ and
\item
there exists $C>0$ and  $t_0>0$ so that if $0<t<t_0$ and $p\in \mathbb P(\mathbb K^d)$,
then 
$$\xi^1_*\mu(B(p,t))\le C t^\delta,$$
where $B(p,t)$ is the ball of radius $t$ in $\Pb(\Kb^d)$ centered at $p$.
\end{enumerate}
\end{lemma}

\begin{proof} 
Let $\mathcal V_n$ denote the set of vertices of $\mathcal T$ at distance $n$ from $b_0$ (with respect to the integer-valued metric on the vertices of $\mathcal T$). 
We inductively define a sequence of Borel probability measures 
supported on $\mathcal V_n$. First let $\mu_0 $ be the Dirac measure $\delta_{b_0}$ at $b_0$. Then, inductively define 
\begin{align*}
\mu_n = \sum_{z \in \mathcal V_{n-1}} \Bigg(\frac{\mu_{n-1}(z)}{\sum_{w \in \mathcal C(z)}  c(w)^{\delta}} \sum_{w \in \mathcal C(z)} c(w)^{\delta} \delta_{w}\Bigg).
\end{align*}

We may view $\{\mu_n\}$ as a sequence of probability measures on the compact space $\mathcal T\cup E$, so it has a weak-$*$ 
subsequential limit $\mu$, which is a probability measure on $\mathcal T\cup E$. Further, since 
$\d_{\Omega}(b_0,\mathcal V_n)\to\infty$, the support of $\mu$ lies in $E$. 

(1): Since the shadows are closed, we introduce open neighborhoods to apply Portmanteau's Theorem to estimate the measures of shadows. 

By Proposition~\ref{prop:finding_children} part (2) and induction, 
if $n \in \Nb$ and $z,w \in \mathcal V_n$ are distinct, then 
\begin{align*}
\mathcal O_{2r_0}(b_0,z) \cap \mathcal O_{2r_0}(b_0,w) = \emptyset.
\end{align*}
Then, since each shadow is closed and each $\mathcal V_n$ is finite, we can construct open sets $\{\mathcal{U}_z\}_{z\in \mathcal T}$ in $\partial \Omega$ such that 
\begin{itemize}
\item  $\mathcal O_{2r_0}(b_0,z) \subset \mathcal{U}_z$ for all $z \in \mathcal T$ and 
\item if  $n \in \Nb$ and $z,w \in \mathcal V_n$ are distinct, then 
\begin{align*}
\overline{\mathcal{U}_z} \cap \overline{\mathcal{U}_w} = \emptyset.
\end{align*}
\end{itemize}
Next let $\mathcal K_z$ be the interior (in $\overline\Omega$) of the closed convex cone of $\overline{\mathcal U_z}$ based at $b_0$.

Now fix a vertex $z$ of $\mathcal T$, and let $n$ be the integer such that $z\in\mathcal V_n$. By Proposition~\ref{prop:finding_children} part (1), the set of vertices in $\cup_{j\ge n}\mathcal V_j$ that lie in  $\mathcal K_z$ are precisely the set of descendants of $z$. Then the definition of $\mu_n$ implies that
\begin{align*}
\mu_m(\mathcal K_z)=\mu_n(\mathcal K_z)
\end{align*}
for all $m\ge n$. Furthermore, if $\{b_0=x_0,x_1,\dots,x_{n-1},x_n=z\}$ is the geodesic in $\mathcal{T}$ between $b_0$ and $z$, then by induction and 
Proposition~\ref{prop:finding_children} part (4), we have 
\[\mu_n(z)=c(z)^{\delta}\prod_{i=0}^{n-1}\frac{c(x_i)^{\delta}}{\sum_{w \in \mathcal C(x_i)}  c(w)^{\delta}}\le c(z)^{\delta}.\]
Since $\mathcal K_z\cap\mathcal V_n=\{z\}$, it now follows that for all $m\ge n$
\[\mu_m(\mathcal K_z)=\mu_n(\mathcal K_z)=\mu_n(z) \le c(z)^{\delta}.\]
Thus, by Portmanteau's Theorem, 
$$
\mu(\mathcal O_{2r_0}(b_0,z))\le \liminf_{m\to\infty}\mu_m(\mathcal K_z) \le c(z)^{\delta}, 
$$
so (1) holds.

(2): Let $z,w \in \mathcal T$ be any two adjacent vertices. By Proposition~\ref{prop:finding_children} part (3), $\d_{\Omega}(z,w) \le D_0$. So if we set
\begin{align*}
S = \max\left\{ \frac{\sigma_1(\rho(\gamma))}{\sigma_d(\rho(\gamma))} : \d_{\Omega}(b_0,\gamma(b_0)) \le D_0\right\},
\end{align*}
then 
\begin{align}
\label{eqref:singular_value_gap}
\frac{1}{S}c(z) \le c(w) \le S c(z).
\end{align}

By Lemma \ref{lem: bishopjones1} part (1), for any $x\in E$  there exists a geodesic ray $\{x_n\}$  in $\mathcal T$ 
starting at $x_0 = b_0$ such that $x_n \rightarrow x$ in $\overline{\Omega}$ and
\begin{align*}
x \in \bigcap_{n \in\mathbb N} \mathcal O_{r_0}(b_0,x_n).
\end{align*}
So there exists a sequence $\{y_n\}$ along the geodesic ray $[b_0,x)$ in $\Omega$ 
with $\d_{\Omega}(y_n, x_n) \le r_0$ for all $n$. 
By Theorem~\ref{thm:size_of_shadows_inside} and Lemma \ref{lem: bishopjones1} part (2), there exists $C_1> 1$ (that depends on $R_0$ but not $x$) such that 
\begin{align}\label{eqn: inradius}
B\left(\xi^1(x), \frac{c(x_n)}{C_1} \right)\cap\Lambda_1(\Gamma)\subset\xi^1\left(\widehat{\mathcal O}_{r_0}(b_0, y_n)\right)\subset\xi^1\left(\widehat{\mathcal O}_{2r_0}(b_0, x_n)\right)
\end{align}
for all $n$. Then by (1),
\begin{align}\label{eqn:push1}
\xi^1_*\mu\left(B\left(\xi^1(x),\frac{c(x_n)}{C_1}\right)\right)\le c(x_n)^{\delta}.
\end{align}

Since $\Gamma_0$ is $P_1$-divergent, $c(x_n)\to 0$. Let
$$t_1=\frac{1}{C_1}\min \{c(z) : z\in \mathcal V_1\}.$$
If $0<t<t_1$, there is some positive integer $n$ such that $c(x_{n+1}) \le tC_1 \le c(x_n)$. 
By \eqref{eqref:singular_value_gap}, 
\begin{align*}
\frac{1}{S}\frac{c(x_n)}{C_1} \le t \le \frac{c(x_n)}{C_1},
\end{align*}
so, Equation \eqref{eqn:push1} implies that
$$\xi^1_*\mu\left(B\left(\xi^1(x),t\right)\right)\le \xi^1_*\mu\left(B\left(\xi^1(x),\frac{c(x_n)}{C_1}\right)\right)\le c(x_n)^{\delta} \le(C_1S)^{\delta} t^{\delta}.$$

Let $t_0=t_1/2$. Suppose that $p\in\mathbb P(\mathbb K^d)$ and $0<t<t_0$.
Since $\mu$ is supported on $E$, either $B(p,t)\cap \xi^1(E)$ is empty, in which case $(\xi^1)^*\mu(B(p,t))=0$ or there exists $x\in E$ so that 
$\xi^1(x)\in B(p,t)\cap \xi^1(E)$, in which case
$B(p,t)\subset B(\xi^1(x),2t)$, so 
$$\xi^1_*\mu\left(B(p,t)\right) \le \xi^1_*\mu\left(B(\xi^1(x),2t)\right) \le(C_1S)^{\delta} (2t)^{\delta}\le (2C_1S)^{\delta} t^{\delta}. $$
Therefore, (2) holds with $C= (2C_1S)^{\delta}$.
\end{proof}

Let $\{B(p_i,r_i)\}$ be a countable covering of $\xi^1(E)$ by open balls each with radii strictly less that $t_0$. Lemma \ref{lem: measure}  implies that
$$\sum_i r_i^\delta\ge \frac{1}{C} \sum_i \xi^1_*\mu(B(p_i,r_i))\ge\frac{1}{C},$$
so $\dim_{H}(\xi^1(E))\ge\delta.$
Since $\xi^1:\Lambda_\Omega(\Gamma_0)\to \Lambda_{1}(\Gamma)$ is a $\rho$-equivariant homeomorphism and $E\subset\Lambda_{\Omega,c}(\Gamma)$,
we see that $\xi^1(E)\subset\Lambda_{1,c}(\Gamma)$.
Therefore, 
\[\dim_{H}(\Lambda_{1,c}(\Gamma))\ge\dim_{H}(\xi^1(E))\ge\delta.\]
Since $\delta$ can be chosen to be any positive number less than $\delta^{\alpha_1}(\Gamma)$, it follows
that
$$\dim_{H}(\Lambda_{1,c}(\Gamma))\ge\delta^{\alpha_1}(\Gamma).$$
By Corollary \ref{prop:upper_bound}, $\dim_{H}(\Lambda_{1,c}(\Gamma))\le\delta^{\alpha_1}(\Gamma),$ and this completes the proof of Theorem \ref{thm:lower_bd}.

\subsection{Proof of Proposition  \ref{prop:finding_children}}

It remains to prove Proposition \ref{prop:finding_children}. We start with some observations about shadows. For any integer  $n\ge 0$, let 
\[\mathcal A_n = \left\{ z \in \Gamma_0(b_0) : e^{-(n+1)} < c(z) \le e^{-n} \right\}.\] 
Notice that, since $\Gamma$ is $P_1$-divergent, $\mathcal A_n$ is finite.
We first show that there is a uniform lower bound on distance which implies that two points in $\mathcal A_n$ have disjoint shadows.

\begin{lemma} \label{lem: close shadows}
For any $r > 0$ there exists $C_0=C_0(r) > 0$ such that if $n \ge 0$, $z,w \in \mathcal A_n$ and $\d_{\Omega}(z,w) > C_0$, then 
$$ \mathcal O_{r}(b_0,z) \cap \mathcal O_{r}(b_0,w)=\emptyset.$$
\end{lemma}

\begin{proof} We prove the contrapositive. Suppose 
$$
x \in  \mathcal O_{r}(b_0,\gamma(b_0)) \cap \mathcal O_{r}(b_0,\eta(b_0))
$$ 
and $\gamma(b_0), \eta(b_0) \in \Ac_n$. Then there exist
$z^\prime, w^\prime \in [b_0,x)_\Omega$ such that $\d_{\Omega}(\gamma(b_0),z^\prime)\le r$ and $\d_{\Omega}(\eta(b_0),w^\prime) \le r$. 
Without loss of generality we can assume that $z^\prime \in [b_0, w^\prime]_\Omega$. Then Equation~\eqref{eqn:Hilbert metric Hausdorff distance} implies that
 $$
 \d_\Omega( \gamma(b_0), [b_0, \eta(b_0)]_\Omega) \le  r+\d_\Omega( z^\prime, [b_0, \eta(b_0)]_\Omega) \le 2r.
 $$
 So, by Lemma~\ref{lem:lower_bd_on_sk_3}, there exists $C > 0$, which only depends on $2r$, such that 
 $$
 \alpha_1(\kappa(\rho(\eta)) )\ge  \alpha_1(\kappa(\rho(\gamma)))+\alpha_1(\kappa(\rho(\gamma^{-1}\eta)))-C.
 $$
 Since $\gamma(b_0), \eta(b_0) \in \Ac_n$,  we have that 
 \[\alpha_1(\kappa(\rho(\eta)))-\alpha_1(\kappa(\rho(\gamma)))=-\log c(\eta(b_0))+\log c(\gamma(b_0))\le 1,\]
 and so
 $$
 \alpha_1(\kappa(\rho(\gamma^{-1}\eta)))\le C+\alpha_1(\kappa(\rho(\eta)))-\alpha_1(\kappa(\rho(\gamma)))  \le C+1.
 $$
 
Thus, if we set
 $$
 C_0 := \max\left\{ \d_\Omega( b_0, \beta(b_0)) :\ \beta\in\Gamma_0\ \text{and}\  \alpha_1(\kappa(\rho(\beta))) \le C+1 \right\},
 $$
 then 
 $$
 \d_\Omega( \gamma(b_0), \eta(b_0)) = \d_\Omega( b_0, \gamma^{-1}\eta(b_0)) \le C_0.
 $$
 Since $C_0$ only depends on $r$ this completes the proof. 
\end{proof} 

The next lemma is needed to define the constants $r_0$ and $D_0$ in Proposition \ref{prop:finding_children}.

 \begin{lemma}\label{obs:shadows} Given $t>0$, there exist $r_0=r_0(t)>0$ and $N_0=N_0(t) > 0$ such that if 
 $z,w \in \Gamma_0(b_0)-\{b_0\}$, $\d_{\Omega}(b_0,w) \ge N_0$, 
 and $\d_{\Pb(\Rb^{d_0})}\big(\iota_{b_0}(z),\iota_{b_0}(w)\big)\ge t$, then 
 \begin{align*}
 \mathcal O_{2r_0}(z,w) \subset \mathcal O_{r_0}(z,b_0) \cap \mathcal O_{2r_0+1}(b_0,w).
 \end{align*}
 \end{lemma} 
 
 \begin{proof}
 Since $\Gamma_0\subset \Aut(\Omega)$ is projectively visible, there exists $r_0 > 0$ such that if $p,q \in \Lambda_\Omega(\Gamma_0)$ and 
 $\d_{\Pb(\Rb^{d_0})}(p,q) \ge t$, then 
 \begin{align*}
 \d_\Omega( b_0, (p,q)_\Omega) \le \frac{r_0}{2}. 
 \end{align*}
 
Suppose, for contradiction, that the lemma fails for this value of $r_0$. Then for every positive integer $n$, there exist 
$z_n, w_n \in \Gamma_0(b_0)-\{b_0\}$ and
  \begin{align*}
x_n \in  \mathcal O_{2r_0}(z_n,w_n) - \left(\mathcal O_{r_0}(z_n,b_0) \cap \mathcal O_{2r_0+1}(b_0,w_n)\right)
 \end{align*}
such that $\d_{\Omega}(b_0,w_n) \ge n$ and $\d_{\Pb(\Rb^{d_0})}(\iota_{b_0}(z_n),\iota_{b_0}(w_n))\ge t$. For each $n$, let $\gamma_n \in \Gamma_0$ be the element with $w_n = \gamma_n(b_0)$, and choose $u_n \in [z_n,x_n)_\Omega$ so that $\d_\Omega(w_n,u_n)\le 2r_0$. By passing to a subsequence, we can suppose that 
\[w_n \rightarrow w\in\Lambda_\Omega(\Gamma_0),\quad z_n \rightarrow z\in\overline{\Omega},\quad x_n \rightarrow x\in\partial\Omega,\quad\gamma_n^{-1}(b_0) \rightarrow \bar{b}\in\Lambda_\Omega(\Gamma_0),\]
\[\gamma_n^{-1}(u_n) \rightarrow \bar{u}\in\Omega\quad \text{and}\quad\gamma_n^{-1}(x_n) \rightarrow \bar{x}\in\partial\Omega.\] 
Then $\d_{\Pb(\Rb^{d_0})}(\iota_{b_0}(z),w)\ge t$ and so (see Figure \ref{fig: radial})
 \begin{align}
 \label{eqn:distance to bo}
 \d_\Omega( b_0, (z,w)_\Omega)\le \d_\Omega( b_0, (\iota_{b_0}(z),w)_\Omega) \le \frac{r_0}{2}. 
 \end{align}
 We will obtain a contradiction by showing that $x_n\in \mathcal O_{r_0}(z_n,b_0) \cap \mathcal O_{2r_0+1}(b_0,w_n)$ for $n$ sufficiently large. 
 
  \begin{figure}[ht]
    \centering
    \includegraphics[width=0.4\textwidth]{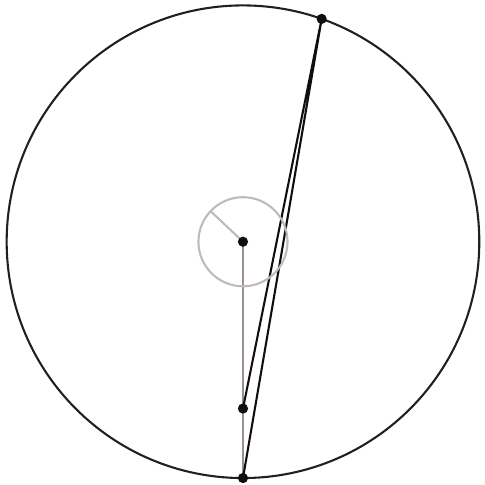}
\tiny
\put (-88,90){$b_0$}
\put (-97,99){$\frac{r_0}{2}$}
\put (-98,28){$z$}
\put (-65,179){$w=x$}
\put (-100,-4){$\iota_{b_0}(z)$}
    \caption{$x\in\mathcal O_{r_0}(z,b_0)$.}
    \label{fig: radial}
\end{figure}
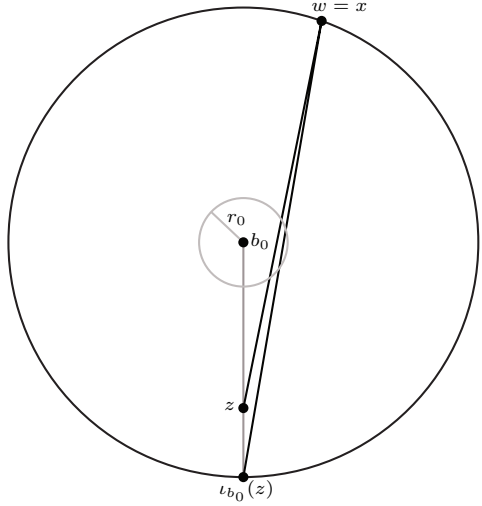

Proposition~\ref{Prop: Islam-Zimmer} part (2) implies that 
\begin{equation}
\label{eqn:forward orbits}
\gamma_n(v) \rightarrow w
\end{equation}
locally uniformly over all $v \in \Pb(\Rb^{d_0}) - T_{\bar{b}} \partial \Omega$ and 
\begin{equation}
\label{eqn:backward orbits}
\gamma_n^{-1}(v) \rightarrow \bar{b}
\end{equation}
locally uniformly over all $v \in \Pb(\Rb^{d_0}) - T_{w} \partial \Omega$. Since $\{\gamma_n^{-1}(u_n)\}$ is relatively compact in $\Omega$, \eqref{eqn:forward orbits} implies that
\[\lim_{n\to\infty}u_n = \lim_{n\to\infty}\gamma_n(\gamma_n^{-1}(u_n))=w.\]
Since $u_n \in [z_n, x_n)_\Omega$ for all $n$, it follows that $w \in [z,x]_\Omega$. Since  $\d_{\Pb(\Rb^{d_0})}(\iota_{b_0}(z),w)\ge t$  and $w\in\partial\Omega$, we must
have $x=w$. Then \eqref{eqn:distance to bo} implies that 
$$
\lim_{n \rightarrow \infty} \d_\Omega( b_0, (z_n,x_n)_\Omega) = \d_\Omega( b_0, (z,x)_\Omega) = \d_\Omega( b_0, (z,w)_\Omega) \le \frac{r_0}{2}. 
$$
Hence $x_n\in\mathcal O_{r_0}(z_n,b_0)$ for $n$ sufficiently large.

 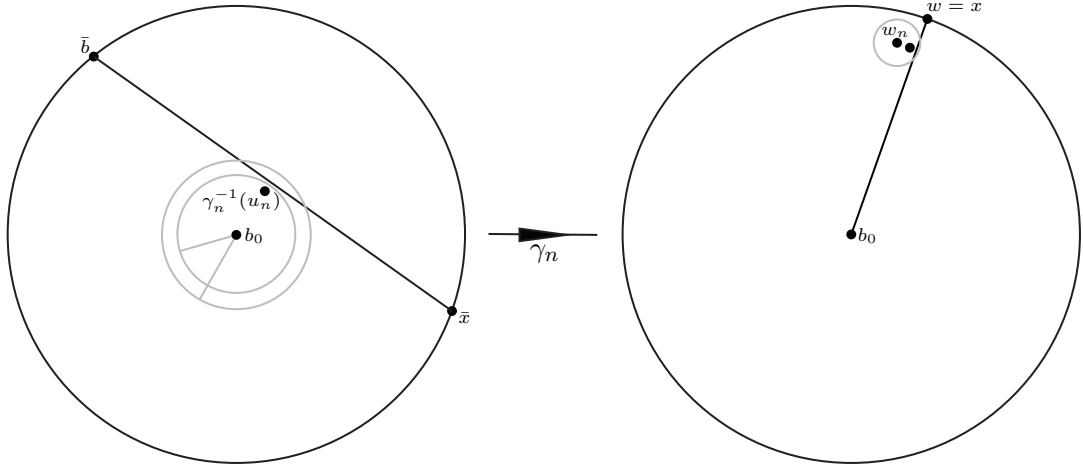
\begin{figure}[ht]
    \centering
    \includegraphics[width=0.9\textwidth]{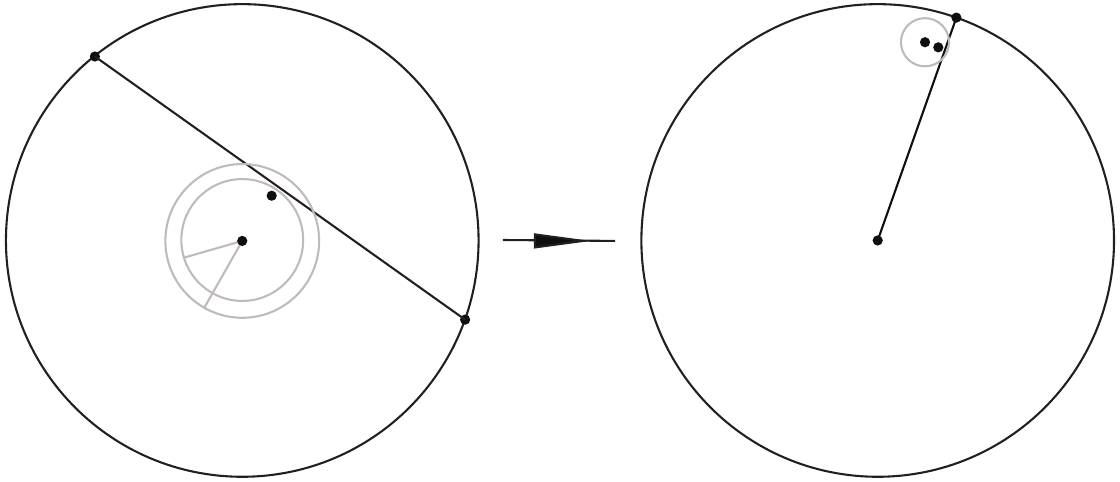}
\put (-210,80){$\gamma_n$}
\tiny
\put (-318,86){$b_0$}
\put (-380,157){$\bar{b}$}
\put (-237,55){$\bar{x}$}
\put (-334,99){$\gamma_n^{-1}(u_n)$}
\put (-86,86){$b_0$}
\put (-60,173){$w=x$}
\put (-77,164){$w_n$}
    \caption{$x\in \mathcal O_{2r_0+1}(b_0,w)$.}
    \label{fig: radial2}
\end{figure}

Since  $(z,w)_\Omega \subset \Omega$,  we see that  $z \notin T_w \partial \Omega$. Hence, $\{z_n\}$ is relatively compact in $\Pb(\Rb^{d_0}) - T_{w} \partial \Omega$, so \eqref{eqn:backward orbits} implies that $\gamma_n^{-1}(z_n) \rightarrow \bar{b}$. Then $\bar{u} \in (\bar{b}, \bar{x})_\Omega$, which gives
$$
\lim_{n \rightarrow \infty} \d_{\Omega}(u_n,(b_0,x_n)_\Omega) = \lim_{n \rightarrow \infty} \d_{\Omega}\left(\gamma_n^{-1}(u_n),\gamma_n^{-1}(b_0,x_n)_\Omega\right) = \d_{\Omega}(\bar{u},(\bar{b},\bar{x})_\Omega)=0.
$$
So for sufficiently large $n$, 
\begin{align*}
\emptyset\neq \overline{B_\Omega(u_n, 1)} \cap (b_0, x_n)_\Omega \subset \overline{B_\Omega(w_n,2r_0+1)} \cap (b_0, x_n)_\Omega,
\end{align*}
or equivalently, $x_n\in \mathcal O_{2r_0+1}(b_0,w_n)$, see Figure \ref{fig: radial2}. Thus $x_n\in \mathcal O_{r_0}(z_n,b_0) \cap \mathcal O_{2r_0+1}(b_0,w_n)$ for $n$ sufficiently large, which is a contradiction. 
 \end{proof}

For any $y\in\partial\Omega$ and any $t>0$, set
\[
\mathcal{B}(y,t)=\left\{w\in\Gamma_0(b_0)-\{b_0\}: \d_{\Pb(\Rb^{d_0})}(\iota_{b_0}(w),y) <t\right\}.
\]
Recall that $0 < \delta < \delta^{\alpha_1}(\Gamma)$. Fix $\epsilon > 0$ such that $0 < \delta < \delta+\epsilon < \delta^{\alpha_1}(\Gamma)$. Then
 \begin{align*}
 \sum_{w\in\Gamma_0(b_0)} c(w)^{\delta + \epsilon} = \infty.
 \end{align*}
Thus, for each $n>0$, there is some $x_n\in\partial\Omega$ such that 
  \begin{align*}
 \sum_{w\in\mathcal B(x_n,1/n)} c(w)^{\delta + \epsilon} = \infty.
 \end{align*}
By taking a subsequence, we may assume that $x_n\to x \in \partial \Omega$. Then for every $t>0$,
  \begin{align}\label{eqn: x infinity}
 \sum_{w\in\mathcal B(x,t)} c(w)^{\delta + \epsilon} = \infty
 \end{align}
 Since $\mathcal B(x,t)$ is infinite for every $t > 0$, we must have $x \in \Lambda_\Omega(\Gamma_0)$. 
 
We prove the following refinement of Equation \eqref{eqn: x infinity}.

\begin{lemma}\label{lem: massive band}
If $y\in\Gamma_0(x)$ and $t>0$, then 
\begin{align*}
\limsup_{n \rightarrow \infty} \sum_{w\in\mathcal A_n\cap \mathcal B(y,t)} c(w)^{\delta} = \infty.
\end{align*}
\end{lemma}

\begin{proof} 
Choose $\gamma\in\Gamma_0$ so that $\gamma(x)=y$. Since $\gamma$ is a diffeomorphism and  $\Pb(\Rb^{d_0})$ is compact,  there exists $D=D(\gamma) > 1$ so that 
\[\d_{\Pb(\Rb^{d_0})}(\gamma(a), \gamma(b)) \le D\d_{\Pb(\Rb^{d_0})}(a,b)\] 
for all $a,b \in \Pb(\Rb^{d_0})$. Also, since
\[\sigma_d(\rho(\gamma))\sigma_j(\rho(\eta))\le \sigma_j(\rho(\gamma\eta))\le \sigma_1(\rho(\gamma))\sigma_j(\rho(\eta))\]
for all $j\in\{1,\dots,d\}$, we may enlarge $D$ if necessary to ensure that 
\[c(\gamma w) \ge \frac{1}{D}c(w)\] 
for all $w \in \Gamma_0(b_0)$. Then
\begin{align*}
\sum_{w\in\mathcal B(y,t)} c(w)^{\delta + \epsilon}&=\sum_{w\in\mathcal B(\gamma(x),t)} c(w)^{\delta + \epsilon} \ge\frac{1}{D^{\delta+\epsilon}} \sum_{w\in\mathcal B(x,t/D)} c(w)^{\delta + \epsilon}= \infty. 
\end{align*}

If $\limsup_{n \rightarrow \infty} \sum_{w\in\mathcal A_n\cap \mathcal B(y,t)} c(w)^{\delta} < \infty$,
then there exists $C > 0$ such that
  \begin{align*}
 \sum_{w\in \mathcal A_n\cap\mathcal B(y,t)} c(w)^{\delta} \le C.
 \end{align*}
for all $n\in\mathbb N$.
Then 
 \begin{align*}
\sum_{w\in\mathcal B(y,t)} c(w)^{\delta+\epsilon} \le\sum_{n=0}^{\infty} \bigg( e^{-\epsilon n}
\sum_{w\in \mathcal A_n\cap\mathcal B(y,t)} c(w)^{\delta}\bigg) \le C  \sum_{n=0}^{\infty} e^{-\epsilon n} < \infty,
 \end{align*}
which is a contradiction.
 \end{proof}

We are now ready to finish the proof of Proposition \ref{prop:finding_children}.  
Choose $x'\in\Gamma(x)-\{x\}$ and fix
$$0<t_0< \frac{1}{4}\d_{\Pb(\Rb^{d_0})}(x,x').$$ 
 
Note that $\sigma_1(\rho(\gamma\eta))\le \sigma_1(\rho(\gamma))\sigma_1(\rho(\eta))$ for all $\gamma, \eta \in \Gamma$. Thus, by  Lemma~\ref{lem:lower_bd_on_sk}, there exists $C_1=C_1(t_0)>0$ so that if $\gamma,\eta\in\Gamma_0$ satisfy $\d_{\Pb(\Rb^{d_0})}(\iota_{b_0}(\eta(b_0)),\iota_{b_0}(\gamma^{-1}(b_0)))\ge t_0$, then
\begin{align}\label{eqn: BPS2} 
\frac{\sigma_2(\rho(\gamma\eta))}{\sigma_1(\rho(\gamma\eta))}\ge C_1\frac{\sigma_2(\rho(\gamma))}{\sigma_1(\rho(\gamma))}\frac{\sigma_2(\rho(\eta))}{\sigma_1(\rho(\eta))}.
\end{align}
Also, let $r_0=r_0(t_0)>0$ and $N_0=N_0(t_0)>0$ be the constants 
given by Lemma \ref{obs:shadows}. Then let $C_0 = C_0(2r_0+1) > 0$ be the constant from Lemma~\ref{lem: close shadows}.
Since $\Gamma_0$ is a discrete group, there exist a finite partition 
\[\Gamma_0(b_0)=P_1\cup\dots\cup P_L\] 
such that each $P_i$ is $C_0$-separated. By definition
$$
\lim_{n \rightarrow \infty} \min_{z \in \mathcal A_n} \d_\Omega(b_0, z) = \infty.
$$
So by Lemma \ref{lem: massive band}, there exist $n, n^\prime \ge 1$ such that 
\begin{align}
   \min_{z \in \mathcal A_n \cup \mathcal A_{n^\prime}} &\d_\Omega(z, b_0) \ge N_0,  \label{eqn: distance bigger than N0} \\
    \sum_{w\in \mathcal A_{n}\cap\mathcal B(x,t_0)} &c(w)^{\delta} \ge \frac{L}{C_1^{\delta}},\text{ and }\nonumber\\
    \sum_{w\in \mathcal A_{n^\prime}\cap\mathcal B(x',t_0)}& c(w)^{\delta} \ge \frac{L}{C_1^{\delta}}. \nonumber
\end{align}
Thus, there exist $i_x,i_{x^\prime}\in\{1,\dots,L\}$ such that 
   \begin{align}
   \label{eqn:sum is large}
 \sum_{w\in S_x} c(w)^{\delta} \ge \frac{1}{C_1^{\delta}}\quad\text{and}\quad \sum_{w\in S_{x^\prime}} c(w)^{\delta} \ge \frac{1}{C_1^{\delta}},
 \end{align}
 where $S_x=P_{i_x}\cap \mathcal A_n \cap \mathcal B(x,t_0)$ and $S_{x'}=P_{i_{x^\prime}}\cap \mathcal A_{n'} \cap \mathcal B(x^\prime,t_0)$. 
 Set 
 $$
 D_0=\max_{z \in \Ac_n \cup \Ac_{n^\prime}} \d_\Omega(b_0, z).
 $$

Fix $z=\gamma(b_0) \in\Gamma(b_0)$. Since $\d_{\Pb(\Rb^{d_0})}(x,x')> 4t_0$, there exists $y\in\{x,x'\}$ such that 
\begin{equation}\label{eqn: 2t0}
\d_{\Pb(\Rb^{d_0})}\left(y,\iota_{b_0}(\gamma^{-1}(b_0)) \right)> 2t_0.
\end{equation} 
Let $\mathcal C(z)= \gamma(S_y)\subset\Gamma_0(b_0)-\{z\}$. We check that $\mathcal C(z)$ satisfies parts (1)--(4) of Proposition \ref{prop:finding_children}. 
 
 Since $S_y\subset\mathcal B(y,t_0)$, \eqref{eqn: 2t0} implies that
 \begin{equation}\label{eqn: t0}
\d_{\Pb(\Rb^{d_0})}\left(\iota_{b_0}(\gamma^{-1}(w)),\iota_{b_0}(\gamma^{-1}(b_0))\right)> t_0
 \end{equation} 
 for all $w\in\mathcal C(z)$. Since $S_y\subset\mathcal A_n \cup \mathcal A_{n^\prime}$, \eqref{eqn: distance bigger than N0} implies that
 \begin{equation}\label{eqn: progress}
 \d_{\Omega}(b_0,\gamma^{-1}(w))  \ge N_0
 \end{equation}
for all $w\in\mathcal C(z)$. Since $S_y\subset P_{i_y}$, $S_y$ is $C_0$-separated, so by Lemma \ref{lem: close shadows}, 
\begin{align}
\label{eqn:disjoint_shadows}
\mathcal O_{2r_0+1}(b_0,\gamma^{-1}(w)) \cap \mathcal O_{2r_0+1}(b_0,\gamma^{-1}(w')) = \emptyset
\end{align}
for all distinct $w,w'\in\mathcal C(z)$. 

Lemma~\ref{obs:shadows}, \eqref{eqn: t0}, and \eqref{eqn: progress} imply that
\begin{align*}
\mathcal O_{2r_0}(b_0,w) = \gamma\Big(\mathcal O_{2r_0}(\gamma^{-1}(b_0), \gamma^{-1}(w))\Big) \subset \gamma\Big(\mathcal O_{r_0}(\gamma^{-1}(b_0),b_0)\Big) = \mathcal O_{r_0}(b_0,z)
\end{align*}
for all $w\in\mathcal C(z)$, so part (1) holds. Lemma~\ref{obs:shadows} and \eqref{eqn:disjoint_shadows} imply that
\begin{align*}
\gamma^{-1} \Big( \mathcal O_{2r_0}(b_0,w) \cap \mathcal O_{2r_0}(b_0,w') \Big)&
= \mathcal O_{2r_0}(\gamma^{-1}(b_0), \gamma^{-1}(w)) \cap \mathcal O_{2r_0}(\gamma^{-1}(b_0), \gamma^{-1}(w')) \\
& \subset \mathcal O_{2r_0+1}(b_0,\gamma^{-1}(w)) \cap \mathcal O_{2r_0+1}(b_0,\gamma^{-1}(w') ) = \emptyset
\end{align*}
for all distinct $w,w'\in\mathcal C(z)$, so part (2) holds. Since $S_y\subset\mathcal A_n \cup \mathcal A_{n^\prime}$, we know that 
\[\d_{\Omega}(z,w)=\d_{\Omega}(b_0,\gamma^{-1}(w)) \le D_0,\] 
for all $w\in\mathcal C(z)$, so part (3) holds. Finally, for each $w\in\mathcal C(z)$, choose $\eta_w \in \Gamma_0$ so  that 
\[\eta_w(b_0) = \gamma^{-1}(w) \in S_y.\] 
By \eqref{eqn: t0}, we have $\d_{\Pb(\Rb^{d_0})}(\iota_{b_0}(\eta_w(b_0)),\iota_{b_0}(\gamma^{-1}(b_0))> t_0$, 
so by \eqref{eqn: BPS2},
\begin{align*}
c(w) =\frac{\sigma_2(\rho(\gamma\eta_w))}{\sigma_1(\rho(\gamma\eta_w))} \ge C_1 \frac{\sigma_2(\rho(\gamma))}{\sigma_1(\rho(\gamma))}\frac{\sigma_2(\rho(\eta_w))}{\sigma_1(\rho(\eta_w))} = C_1c(z) c( \gamma^{-1} (w)) .
\end{align*}
Then, by \eqref{eqn:sum is large},
\[\sum_{w\in\mathcal C(z)} c(w)^{\delta} \ge C_1^{\delta} c(z)^{\delta} \sum_{w\in S_y} c(w)^{\delta} \ge  c(z)^{\delta}.\qedhere\]
This completes the proof of  Proposition \ref{prop:finding_children} and hence the proof of Theorem \ref{thm:lower_bd}.


\section{Critical exponents and entropies}


In this section, we show that critical exponents and entropies agree for cusped $P_\theta$-Anosov representations of geometrically finite Fuchsian groups. This generalizes results of 
Glorieux--Montclair--Tholozan \cite[Thm.\ 3.1]{GMT}  and Pozzetti--Sambarino--Weinhard \cite[Prop.\ 4.1]{PSW1} from the Anosov setting. 
Notice that Proposition \ref{critandent} from the introduction is a special case of this result.

\begin{proposition}
\label{entropy and critical exponent}
If $\Gamma\subset\mathsf{PSL}(2,\mathbb R)$ is geometrically finite, $\rho:\Gamma\to \mathsf{PGL}(d,\mathbb K)$ is $P_\theta$-Anosov, and
$\phi\in \Bc_\theta^+(\rho)$, then $h^{\phi}(\rho)= \delta^{\phi}(\rho)$.
\end{proposition}

\begin{remark} After this paper was submitted for publication, the above proposition was established for all relatively Anosov groups, see~\cite[Cor.\ 1.7]{BCZZ}.
\end{remark}

\begin{proof} We will make crucial use of the fact that if there is a definite angle between the attracting $k$-plane and the repelling $(d-k)$-plane of a $P_k$-proximal
element, then the $k^{th}$ fundamental weights of the Cartan and Jordan projections are uniformly close.

\begin{lemma}\cite[Lem.\ 2.36]{GMT}
\label{angle weights}
Given $\eta>0$, there exists $C>0$ so that if $g\in\PGL(d,\Kb)$ is $P_k$-proximal, $V$ is its attracting $k$-plane and $W$ is its repelling $(d-k)$-plane,
and $\angle(V,W)\ge\eta$, then
$$\big|\omega_k(\nu(g))-\omega_k(\kappa(g))\big|\le C.$$
\end{lemma}

We use this to control the difference  between the Jordan and Cartan projections of images of elements whose axes in $\mathbb H^2$ pass
through a compact  subset of $\mathbb H^2$.

\begin{lemma}
\label{compact weights}
Given a compact subset $K\subset\Hb^2$, there exists $C>0$ so that if $\gamma\in\Gamma_{hyp}$ and its axis passes through $K$, then 
$$\big|\omega_k(\kappa(\rho(\gamma)))-\omega_k(\nu(\rho(\gamma)))\big|\le C$$
for all $k\in\theta$.
\end{lemma}

\begin{proof}
By compactness and the transversality of $\xi$, there exists $\eta>0$ so that if the geodesic with endpoints $z,w\in\Lambda(\Gamma)$ intersect $K$, then $\angle(\xi^k(z),\xi^{d-k}(w))\ge\eta$ for all $k\in\theta$. We may now apply Lemma~\ref{angle weights} to deduce the lemma.
\end{proof}

We first prove that $h^{\phi}(\rho)\ge \delta^{\phi}(\rho)$. Denote $a(\phi)=\sum_{k\in\theta} |a_k|$, where 
\[\phi=\sum_{k\in\theta} a_k\omega_k.\]
Since $\Gamma$ is geometrically finite, there is a compact subset $K\subset\Hb^2$ so that if $\gamma\in\Gamma$ is hyperbolic, 
then there exists a conjugate of $\gamma$ whose axis intersects $K$. Moreover, by Lemma~\ref{compact weights}, there exists $C>0$ such that
$$\big|\phi(\kappa(\rho(\gamma)))-\phi(\nu(\rho(\gamma)))\big|\le a(\phi)  C$$
for all $\gamma\in\Gamma_{hyp}$ whose axis intersects $K$. It follows that
\begin{align*}
\#\left\{[\gamma]\in[\Gamma_{hyp}]\ \Big|\ \ell^\phi(\gamma)\le T\right\}&\le
\#\left\{\gamma\in\Gamma_{hyp}\ \Big|\ \text{the axis of }\gamma\text{ intersects }K\text{ and } \phi(\nu(\rho(\gamma)))\le T\right\}\\
&\le
\#\left\{\gamma\in\Gamma\ \Big|\ \phi(\kappa(\rho(\gamma)))\le T+ a(\phi)C\right\}
\end{align*}
for all $T>0.$ Since 
$$\delta^\phi(\rho) = \limsup_{T \rightarrow \infty} \frac{1}{T} \log\#\left\{\gamma\in\Gamma\ \Big|\ \phi(\kappa(\rho(\gamma)))\le T\right\},$$
this implies that $h^{\phi}(\rho)\le \delta^{\phi}(\rho).$

We now prove the opposite inequality. A special case of a result of Abels--Margulis--Soifer \cite[Thm.\ 4.1]{AMS}
implies that there exist $\mu>0$ and a finite subset $\mathcal A$ of $\Gamma$ such that if
$\gamma\in\Gamma$, then there exists $\alpha\in\mathcal A$ so that $\gamma\alpha\in\Gamma_{hyp}$ and 
$\d_{\partial\Hb^2}((\gamma\alpha)^+,(\gamma\alpha)^-)\ge\mu$.

Let $K'\subset\Hb^2$ be a compact set so that any bi-infinite geodesic whose
endpoints are at least $\mu$ apart pass through $K' $. By Lemma \ref{compact weights}, there is some $C'>0$ such that
$$\big|\phi(\kappa(\rho(\gamma)))-\phi(\nu(\rho(\gamma)))\big|\le a(\phi)  C'$$
for all $\gamma\in\Gamma_{hyp}$ whose axis intersects $K' $. For all $\gamma\in\Gamma$, set
\[\Ac_\gamma=\{\eta\in(\gamma\mathcal{A})\cap\Gamma_{hyp}:\text{ the axis of }\eta\text{ intersects }K' \}.\]
Since $\Ac$ is finite, there exists $E > 0$ so that if $\gamma\in\Gamma$ and $\eta\in\Ac_\gamma$, then
\[\phi(\kappa(\rho(\eta)))\le E+\phi(\kappa(\rho(\gamma))).\]

For any $\eta\in\Gamma_{hyp}$, let $\ell(\eta)$ denote the minimal translation distance of $\eta$ in $\Hb^2$. 
Since $K' $ is compact,  there exists a positive integer $D$ such that for any $\eta\in\Gamma_{hyp}$, at most $D \ell(\eta)$ of the $\Gamma$-translates of the axis 
of $\eta$ intersect $K' $, 
or equivalently, at most $D\ell(\eta)$ conjugates of $\eta$ have axes that intersect $K' $. Also, since $\rho$ is $P_\theta$-Anosov, there are constants $C>0$ and $c>1$ such that 
\[\ell(\eta)\le \d_{\Hb^2}(\eta(b_0),b_0)\le c\phi(\kappa(\rho(\eta)))+C\] 
for all $\eta\in\Gamma_{hyp}$. Furthermore, note that any element $\eta\in\Gamma$ can lie in $\mathcal A_\gamma$ for at most $|\mathcal  A|$  different $\gamma\in\Gamma$. Thus, for all $T>0$, 
{\small\begin{align*}
\#\left\{\gamma\in\Gamma\ \Big|\ \phi(\kappa(\rho(\gamma))) \le T\right\}&\le |\mathcal  A|\#\left\{\eta\in\Gamma\ \Big|\ \eta\in \mathcal A_\gamma\text{ for some }\gamma\in\Gamma\text{ such that } \phi(\kappa(\rho(\gamma)))\le T\right\}\\
&\le |\mathcal  A|\#\left\{\eta\in\Gamma_{hyp}\ \Big|\ \text{the axis of }\eta\text{ intersects }K' \text{ and } \phi(\kappa(\rho(\eta)))\le E+T\right\}\\
&\le |\mathcal  A|D(c(E+T)+C)\#\left\{[\eta]\in[\Gamma_{hyp}]\ \Big|\ \phi(\nu(\rho(\gamma)))\le E+T+a(\phi)C'\right\}.
\end{align*}
}
This implies that $\delta^{\alpha_1}(\rho)\le h^{\alpha_1}(\rho)$.
\end{proof}

We use a similar argument to show that the $\phi$-Poincar\'e series $Q_{\rho(\Gamma)}^\phi$ diverges at its critical exponent.

\begin{proposition}
\label{crit exp divergence}
If $\Gamma\subset\mathsf{PSL}(2,\mathbb R)$ is geometrically finite, $\rho:\Gamma\to \mathsf{PGL}(d,\mathbb K)$ is $P_\theta$-Anosov, and
$\phi\in \Bc_\theta^+(\rho)$, then $Q_{\rho(\Gamma)}^\phi(\delta^\phi(\rho))=+\infty$.
\end{proposition}

\begin{remark} After this paper was submitted for publication, the above proposition was established for all relatively Anosov groups, see~\cite[Thm.\ 1.1]{CZZ4}
\end{remark} 

\begin{proof}
It follows immediately from Theorem \ref{BCKM counting} that there exists $N\in\mathbb N$ and $A>0$ so that if
$n\ge N$ and
$$ S_n^\phi(\rho)=\left\{ [\gamma]\in[\Gamma_{hyp}]\ :\ n\le \phi(\nu(\rho(\gamma)))< n+1\right\},\qquad\text{then}\qquad \# S_n^\phi(\rho)\ge \frac{Ae^{n\delta^\phi(\rho)}}{n}.$$
The argument above shows that  if $[\gamma]\in S_n^\phi(\rho)$, then there exists $\hat\gamma\in[\gamma]$ so that
$$\phi(\kappa(\rho(\hat\gamma))\le n+1+a(\phi)C.$$ 
Therefore, if $n\ge N$, 
$$\sum_{\{\gamma\ : [\gamma]\in S_n^\phi(\rho)\} } e^{-\delta^\phi(\rho)\phi( \kappa(\rho(\gamma)))}
 \ge \frac{Ae^{n\delta^\phi(\rho)}}{n}e^{-\delta^\phi(\rho)(n+1+a(\phi)C)}\ge \frac{Ae^{-\delta^\phi(\rho)(1+a(\phi)C)}}{n}$$
which implies that $Q_{\rho(\Gamma)}^\phi(\delta^\phi(\rho))=+\infty$.
\end{proof}

\begin{remark}
Propositions \ref{entropy and critical exponent} and \ref{crit exp divergence} do not hold for all $P_\theta$-transverse groups and respresentations.
If $\Gamma$ is a geometrically infinite surface subgroup of $\mathsf{PO}(3,1)$ such that $\Lambda(\Gamma)$ is not all of $\partial\mathbb H^3$,
then $\Gamma$ is a $P_1$-transverse subgroup of $\mathsf{PSL}(4,\mathbb R)$. However,
$R^{\omega_1}_T(\Gamma)$ is infinite for all sufficiently large $T$ (see Bonahon \cite{bonahon-ends}), so $h^{\omega_1}(\Gamma)=+\infty$ with our definition.
Moreover, $\delta^{\omega_1}(\Gamma)=2$ (see Canary \cite[Cor.\ 4.2]{canary-laplace}) 
and  the $\alpha_1$-Poincar\'e series  of $\Gamma$ converges at $s=2$ (see Sullivan \cite[Thm.\ II]{sullivan-ergodic}). 
Also notice that the inclusion map of $\Gamma$ into $\mathsf{PSL}(4,\mathbb R)$ is a $P_1$-transverse representation.
\end{remark}


\section{Quint's indicator set}


In this section, we develop analyticity and convexity properties for the set of linear functionals with entropy 1.

Suppose that $\Gamma$ is a geometrically finite Fuchsian group and
$\rho:\Gamma\to\mathsf{PGL}(d,\Kb)$ is $P_\theta$-Anosov.
Following Sambarino \cite{sambarino-hyperconvex}, we define the {\em Quint indicator set}
$$\mathcal{Q}_\theta(\rho)=\{\phi\in\mathcal B_\theta^+(\rho(\Gamma)) \ |\ h^\phi(\rho)=1\}$$
which arises as the boundary of
$$\mathcal{R}_\theta(\rho)=\{\phi\in \mathcal{B}_\theta^+(\rho(\Gamma)) \ |\  h^\phi(\rho)\le 1\}.$$

The following result generalizes a result of Sambarino \cite[Prop.\ 4.7]{sambarino-hyperconvex} for  Hitchin representations, see also \cite[Prop.\ 4.11]{potrie-sambarino} and
 \cite{sambarino-dichotomy}. (Sambarino \cite{sambarino-hyperconvex} was inspired by work of Quint \cite{quint-ps} who studied the behavior of the entropy on the space of linear 
functionals in the case that the representation is Zariski dense.) Theorem \ref{quintintro} from the introduction is a special case of the result below.

\begin{theorem}
\label{quint indicator}
Suppose that $\Gamma\subset\mathsf{PSL}(2,\mathbb R)$ is geometrically finite and
$\rho:\Gamma\to\mathsf{PGL}(d,\mathbb R)$ is a $P_\theta$-Anosov representation. Then
\begin{enumerate}
\item
$\mathcal{Q}_\theta(\rho)$ intersects each ray (based at  ${\bf 0}$) in $\mathcal{B}_\theta^+(\rho(\Gamma))$ exactly once.
\item
$\mathcal{Q}_\theta(\rho)$ is an analytic submanifold of $\mathfrak{a}_\theta^*$.
\item
$\mathcal{R}_\theta(\rho)$ is a convex subset of $\mathfrak{a}_\theta^*$.
\item
If $\phi_1,\phi_2\in\mathcal{Q}_\theta(\rho)$, then the line segment in $\mathfrak{a}_\theta^*$ between $\phi_1$ and $\phi_2$ lies in $\mathcal{Q}_\theta(\rho)$
if and only if 
$$\phi_1(\nu(\rho(\gamma)))=\phi_2(\nu(\rho(\gamma)))$$
for all $\gamma\in\Gamma$.
\end{enumerate}
\end{theorem}

\begin{remark} After this paper was submitted for publication, part (4) of the above theorem was established for all relatively Anosov groups, see~\cite[Thm.\ 1.8]{CZZ4}
\end{remark} 

In this proof, we will use the technology developed in Bray--Canary--Kao--Martone \cite{BCKM}. We note that the results of that paper were stated for representations into $\SL(d,\Rb)$, but a careful reading verifies that the same arguments taken verbatim work for representations into $\PGL(d,\Rb)$. 

One first associates to any torsion-free, geometrically finite Fuchsian group $\Gamma$, 
a topologically mixing, countable, one-sided Markov shift $(\Sigma^+,\sigma)$ with the big images and pre-images property and countable alphabet $\mathcal A$.
If $\Gamma$ is convex cocompact, 
we use the Bowen-Series coding \cite{bowen-series}, if $\Gamma$ is a non-cocompact lattice, we use the
Stadlbauer--Ledrappier--Sarig coding \cite{ledrappier-sarig,stadlbauer}, and in the remaining cases, we use the coding of Dal'bo--Peign\'e \cite{dalbo-peigne}. 
In all cases, such a coding gives a pair of maps 
\[
G:\mathcal A\to \Gamma\quad\text{and}\quad\omega:\Sigma^+\to\Lambda(\Gamma)
\]
with the property that if $x=\overline{x_1\cdots x_n}\in {\rm Fix}^n(\Sigma^+)$, then $\omega(x)$ is the attracting fixed point of $G(x_1)\dots G(x_n)$. 
If $\rho:\Gamma\to \mathsf{PGL}(d,\mathbb R)$ is $P_\theta$-Anosov, then one can use Quint's Iwasawa cocycle \cite{quint-ps} to construct a 
locally H\"older continuous, vector-valued function 
$\tau_\rho:\Sigma^+\to\mathfrak{a}_\theta$ so that if $x=\overline{x_1\cdots x_n}\in {\rm Fix}^n(\Sigma^+)$, then
$$S_n\tau_\rho(x)=\sum_{i=0}^{n-1} \tau_\rho(\sigma^i(x))=p_\theta(\nu(G(x_1)\cdots G(x_n)))$$
see \cite[Thm.\ ${\rm D}^*$]{BCKM}. See \cite{BCKM} for careful definitions and more detailed statements.

Given a locally H\"older continuous function $g:\Sigma^+\to\mathbb R$ one may define its \emph{Gurevich pressure}
$$P(g)=\sup_{\mu\in\mathcal M} h_\sigma(\mu)+\int_{\Sigma^+} g d\mu$$
where $\mathcal M$ is the space of $\sigma$-invariant probability measures on $\Sigma^+$ and $h_\sigma(\mu)$ is the measure-theoretic entopy of $\sigma$ with respect to $\mu$.
Notice that this pressure need not be finite. However, it is analytic and convex on the space of finite pressure,  locally $\alpha$-H\"older continuous functions for any $\alpha>0$ 
(see Mauldin--Urbanski \cite[Thm.\ 2.6.12]{MU} or Sarig \cite[Cor.\ 4]{sarig-2003}-\cite[Prop.\ 4.4]{sarig-2009}). 
A measure $\mu_g\in\mathcal M$ is said to be an {\em equilibrium measure} for a locally H\"older continuous
function $g:\Sigma^+\to\mathbb R$ if 
$$P(g)=h_\sigma(\mu_g)+\int_{\Sigma^+} g d\mu_g.$$

\begin{proof}[Proof of Theorem \ref{quint indicator}]
First, notice that we may assume throughout that $\Gamma$ is torsion-free, since any geometrically finite Fuchsian group $\Gamma$ has a finite
index torsion-free
subgroup $\Gamma_0$ and $h^\phi(\rho(\Gamma))=h^\phi(\rho(\Gamma_0))$ for all $\phi\in\mathcal B^+_\theta(\rho(\Gamma))=\mathcal B_\theta^+(\rho(\Gamma_0))$.

If $\alpha$ is a parabolic element of $\Gamma$, then, by \cite[Cor.\ 4.2]{CZZ}, the quantity
\[c_k(\rho,\alpha)=\lim_{s\to\infty}\frac{ \omega_k(\kappa(\rho(\alpha^s)))}{\log s}\] 
is well-defined positive integer. If  $\phi=\sum_{k\in\theta} a_k\omega_k\in\mathcal B_\theta^+(\rho(\Gamma))$, let
\[c(\rho,\phi)=\inf\left\{\sum_{k\in\theta} a_kc_k(\rho,\alpha): \alpha\in\Gamma\text{ parabolic}\right\}.\]
We recall the following results from Bray--Canary--Kao--Martone \cite{BCKM}.

\begin{proposition}\label{prop: BCKM}
Let $\phi\in \mathcal B_\theta^+(\rho(\Gamma))$. 
\begin{enumerate}
\item \cite[Cor.\ 1.2, Lem.\ 3.3, Thm.\ ${\rm D}^*$]{BCKM} $0<h^\phi(\rho)<+\infty$ and $P(-t\phi\circ\tau_\rho)=0$ if and only if $t=h^\phi(\rho)$. 
\item \cite[Lem.\ 3.3, Thm.\ ${\rm D}^*$]{BCKM} $P(-\phi\circ\tau_\rho)<+\infty$ if and only if $c(\rho,\phi)>1$. In particular,
$$\mathcal E(\rho)=\{\phi\in \mathcal B_\theta^+(\rho(\Gamma))\ |\ P(-\phi\circ\tau_\rho)<+\infty\}$$
is an open subset of $\mathfrak{a}_\theta^*$.
\item  \cite[Lem.\ 3.4, Thm.\ ${\rm D}^*$]{BCKM} If $P(-\phi\circ\tau_\rho)<+\infty$, then there is a unique equilibrium measure
$dm_{-\phi\circ\tau_\rho}$ for $-\phi\circ\tau_\rho$.
\end{enumerate}
\end{proposition}

Notice that (1) follows immediately from Proposition \ref{prop: BCKM} part (1) and the observation that $h^{k\phi}=\frac{1}{k}h^\phi$ for all $k>0$.

The function $\hat P:\mathcal E(\rho)\to \mathbb R$ given by 
\[\hat P(\phi)=P(-\phi\circ\tau_\rho)\] 
is convex and analytic since $P$ is convex and analytic. By Proposition \ref{prop: BCKM} part (1), $\hat P(\phi)=0$ if and only if $h^\phi(\rho)=1$. Also, note that $\hat P(k\phi)<\hat P(\phi)$ for all $1<k<\infty$. Thus,
\[\mathcal R_\theta(\rho)=\hat P^{-1}((-\infty,0])\,\text{ and }\,\mathcal Q_\theta(\rho)=\hat P^{-1}(0).\] 
This proves (2) and (3).

Let $\phi_1,\phi_2\in\mathcal{Q}_\theta(\rho)$ be distinct. Suppose that the line segment in $\mathfrak{a}_\theta^*$ between $\phi_1$ and $\phi_2$ also lies in $\mathcal{Q}_\theta(\rho)$, or equivalently, that
$$f(t)=P( (t-1)\phi_1\circ\tau_\rho-t\phi_2\circ\tau_\rho)=0$$
for all $t\in[0,1]$. For any $t\in [0,1]$, Proposition \ref{prop: BCKM} part (3) states that there is a unique equilibrium measure $m_t$ for the function $(t-1)\phi_1\circ\tau_\rho-t\phi_2\circ\tau_\rho$, so Mauldin and Urbanski's formula for the derivative of the pressure function \cite[Prop.\ 2.6.13]{MU} (see also
the discussion in \cite[Sec.\ 2]{BCKM}) implies that
$$f'(t)=\int_{\Sigma^+} \big(\phi_1\circ\tau_\rho-\phi_2\circ\tau_\rho\big)\ dm_t=0.$$
(Notice that by Proposition \ref{prop: BCKM} part (2), $\Ec(\rho)$ is open, so $f$ can be defined on some open interval containing $[0,1]$.) In particular, we can write the equality $f'(1)=0$ as
$$\int_{\Sigma^+} \phi_1\circ\tau_\rho\ dm_1=\int_{\Sigma^+}\phi_2\circ\tau_\rho\ dm_1.$$

Since, $P(-\phi_2\circ\tau_\rho)=0$ and $m_1$ is an equilibrium measure for $-\phi_2\circ \tau_\rho$,
$$h_\sigma(m_1)=\int_{\Sigma^+}\phi_2\circ \tau_\rho\  dm_1$$
so we see that
$$h_\sigma(m_1)-\int_{\Sigma^+} \phi_1\circ\tau_\rho\ dm_1
=\int_{\Sigma^+} \phi_2\circ\tau_\rho \ dm_1-\int_{\Sigma^+} \phi_1\circ\tau_\rho\ dm_1=0.$$
So, $m_1$ is also an equilibrium measure for
$-\phi_1\circ\tau_\rho$. Since Proposition \ref{prop: BCKM} part (3) implies that $-\phi_1\circ\tau_\rho$ has a unique equilibrium measure,
$m_0=m_1$. Sarig  \cite[Thm.\ 4.8]{sarig-2009} 
showed that this only happens
when $\phi_1\circ \tau_\rho$ and $\phi_2\circ\tau_\rho$ are cohomologous, 
so the Livsic Theorem
\cite[Thm.\ 1.1]{sarig-2009} (see also Mauldin-Urbanski \cite[Thm.\ 2.2.7]{MU}) implies that
$\phi_1(\nu(\rho(\gamma)))=\phi_2(\nu(\rho(\gamma)))$
for all $\gamma\in\Gamma$.
We have completed the proof.
\end{proof}

Since domain groups of (traditional) Anosov representations admit topologically mixing, finite Markov codings (see
Bridgeman--Canary--Labourie--Sambarino \cite{BCLS} or Constantine--LaFont--Thompson \cite{CLT}), one may apply
the exact same argument to obtain the analogous result in this setting, see also Sambarino \cite{sambarino-dichotomy}.

\begin{theorem}
Suppose that $\Gamma$ is a hyperbolic group and
$\rho:\Gamma\to\mathsf{PGL}(d,\mathbb R)$ is a $P_\theta$-Anosov representation. Then
\begin{enumerate}
\item
$\mathcal{Q}_\theta(\rho)$ intersects each ray (based at  ${\bf 0}$) in $\mathcal{B}_\theta^+(\rho(\Gamma))$ exactly once.
\item
$\mathcal{Q}_\theta(\rho)$ is a analytic submanifold of $\mathfrak{a}_\theta^*$.
\item
$\mathcal{R}_\theta(\rho)$ is a convex subset of $\mathfrak{a}_\theta^*$.
\item
If $\phi_1,\phi_2\in\mathcal{Q}_\theta(\rho)$, then the line segment in $\mathfrak{a}_\theta^*$ between $\phi_1$ and $\phi_2$ lies in $\mathcal{Q}_\theta(\rho)$
if and only if 
$$\phi_1(\nu(\rho(\gamma)))=\phi_2(\nu(\rho(\gamma)))$$
for all $\gamma\in\Gamma$.
\end{enumerate}
\end{theorem}


\section{Entropy rigidity} 


In this section we establish our main result. It remains to prove Proposition \ref{hitchinlimitset}, which bounds the Hausdorff dimension of the limit set of 
a Hitchin representation, check that all exterior powers of Hitchin representations are $(1,1,2)$-hypertransverse and carefully apply Sambarino's classification
of Zariski closures of Hitchin representations.

\subsection{Hausdorff dimension of positive curves} 
We prove that the image of a continuous, positive map has Hausdorff dimension at most $1$.
The proof below is a mild generalization of the standard proof that the graph of a monotonic function $f:[a,b]\to \mathbb R$ has Hausdorff dimension 1. 
Notice that Proposition \ref{hitchinlimitset} is an immediate consequence.

\begin{proposition}\label{pairwise transverse}
Let $X\subset\mathbb{S}^1$ be a closed subset. If $\xi:X\to\Fc$ is a continuous, positive map, then $\dim_H(\xi(X))\le 1$.
\end{proposition}

\begin{proof} If $X$ is a finite set, then $\dim_H(\xi(X))=0$. Hence, we may assume that $X$ is infinite.

Let $x,y,z\in X$ be mutually distinct points and let 
\[I:=\{w\in X:x\le w\le y<z<x\}.\]
Since any positive $n$-tuple of flags consists of mutually transverse flags, see Lemma \ref{lem: basic positivity}, $\xi(I)$ is a compact subset of the affine chart
\[\mathbb{A}_{\xi(z)}:=\{F\in\Fc:F\text{ is transverse to }\xi(z)\}.\] 
Furthermore, if we let $U_{\xi(z)}\subset\SL(d,\Rb)$ denote the subgroup of unipotent elements that fix $\xi(z)$, then we have a real analytic diffeomorphism 
\[L:U_{\xi(z)}\to \mathbb{A}_{\xi(z)}\] 
given by $L(w)=w(\xi(x))$. Let
\[u:=L^{-1}\circ\xi:I\to U_{\xi(z)}.\]
To prove the lemma, it is sufficient to show that $\dim_H(u(I))\le 1$ (with respect to any Riemannian metric on $U_{\xi(z)}$).

Observe that any ordered basis compatible with $(\xi(z),\xi(x))$ defines an identification of $U_{\xi(z)}$ with the group of unipotent, upper triangular matrices in $\SL(d,\Rb)$. Since $\xi$ is a positive map, there exists an ordered basis $\Bc$ compatible with $(\xi(z),\xi(x))$ so that if $s,t\in I$ such that $x \le s<t \le  y$, then 
\[u(s)^{-1}u(t)\in U_{>0}(\Bc)\subset U_{\xi(z)}.\]
In particular, every upper triangular entry of $u(t)$ is strictly larger than the corresponding upper triangular entry of $u(s)$. Thus, it is now sufficient to prove that if $I\subset [0,1]$ is a closed subset and 
 \[f=(f_1,\dots,f_D):I\to\Rb^D\] 
 is a continuous map such that $f_i:I\to\Rb$ is an increasing function for all $i$, then $\dim_H(f(I))\le 1$ (in the Euclidean metric) in $\Rb^\ell$. We will verify the stronger condition that the $1$-dimensional Hausdorff measure of $f(I)$, denoted $\mathcal{H}^1(f(I))$, is finite.

 For every $n\in\Zb^+$ and $k\in\{0,\dots,n-1\}$, let 
 \[I_{n,k}=\left[\frac{k}{n},\frac{k+1}{n}\right]\cap I,\]
 and let
 \[R_{n,k}=[f_1(\min I_{n,k}),f_1(\max I_{n,k})]\times\dots\times [f_D(\min I_{n,k}),f_D(\max I_{n,k})]\subset\Rb^D.\]
Observe that for all $n\in\Zb^+$, $\{R_{n,k}:k=0,\dots,n-1\}$ covers $f(I)$ because $f_i$ is increasing for all $i$. Furthermore,
 \[\sum_{k=0}^{n-1}{\rm diam}(R_{n,k})\le\sum_{k=0}^{n-1}\left(\sum_{i=1}^D f_i\left(\max I_{n,k}\right)-f_i\left(\min I_{n,k}\right)\right)\le\sum_{i=1}^D f_i\left(\max I\right)-f_i\left(\min I\right).\]
Since $f$ is continuous,
\[\lim_{n\to\infty}\sup\{{\rm diam}(R_{n,k}):k=0,\dots,n-1\}=0.\]
It follows that $\mathcal{H}^1(f(I))\le \sum_{i=1}^D f_i\left(\max I\right)-f_i\left(\min I\right)<\infty$.
\end{proof}

\subsection{Simple root critical exponents of Hitchin representations} 
We first observe that all exterior powers of Hitchin representations are $(1,1,2)$-hypertransverse. Let $E^k : \PGL(d,\Rb) \rightarrow \PGL(\bigwedge^k\Rb^d)$ be the representation defined in Section~\ref{sec: skew-symmetric tensors}.

\begin{proposition}\label{prop: Hitchin transverse}
If $\rho:\Gamma\to\PSL(d,\Rb)$ is a Hitchin representation, then it is $P_\Delta$-transverse. Furthermore, $E^k\circ\rho(\Gamma)$ is $(1,1,2)$-hypertransverse for all $k\in\{1,\ldots,d-1\}$.
\end{proposition}

In the case where $\Gamma$ is a cocompact lattice, the hypertransversality follows from the work of Labourie \cite{labourie-invent}, see \cite[Prop.\ 9.6]{PSW1}. We give a more direct proof using positivity, which holds for any discrete $\Gamma\subset\PSL(2,\Rb)$.

\begin{proof}
In \cite[Thm.\ 7.1]{CZZ} we proved that if $\Gamma$ is finitely generated, then $\rho$ is $P_\Delta$-Anosov, hence $P_\Delta$-transverse. 
If $\Gamma$ is infinitely generated, the proof given there shows that $\rho$ is $P_\Delta$-transverse. Let $\xi : \Lambda(\Gamma) \rightarrow \Fc$ be the boundary map (which is positive). 

Fix $k \in \{1,\ldots,d-1\}$. We begin by defining the boundary map of $E^k\circ\rho$. For all $F\in\Fc$, choose a basis $(v_1,\dots,v_d)$ of $\Rb^d$ that is compatible with $F$, 
i.e. $F^i=\Span_{\Rb}(v_1,\dots,v_i)$ for all $i=1,\dots,d$. This induces a basis 
\[(w_1,\dots,w_\ell)=(v_{i_1}\wedge\dots\wedge v_{i_k})_{1\le i_1<\dots<i_k\le d}\quad\text{ of }\quad\bigwedge^k\Rb^d,\] 
enumerated according to the dictionary ordering in the subscripts. Let $D = \dim \bigwedge^k\Rb^d$ and $\hat{\Fc} = \Fc_{1,2,D-2,D-1}(\bigwedge^k\Rb^d)$. 
One can verify that we have a well-defined map $W_k:\Fc\to \hat{\Fc}$ given by
\[W_k(F)=\left(\Span_{\Rb}(w_1),\Span_{\Rb}(w_1,w_2),\Span_{\Rb}(w_1,\dots,w_{D-2}),\Span_{\Rb}(w_1,\dots,w_{D-1})\right).\] 
i.e. $W_k(F)$ does not depend on the choice of the basis compatible with $F$. Then define 
$$
\zeta:=W_k \circ \xi : \Lambda(\Gamma) \rightarrow \hat{\Fc}.$$ 
It is clear that $\zeta$ is continuous and $E^k\circ\rho$-equivariant. 

To show that $\zeta$ is transverse, fix $x, y \in \Lambda(\Gamma)$ distinct. Then fix a basis $\Bc=(v_1,\dots,v_d)$ of $\Rb^d$ such that $v_i\in \xi(x)^i\cap \xi(y)^{d-i+1}$ for all $i=1,\dots,d$. 
Then writing $\zeta(x)$ and $\zeta(y)$ in terms of the basis $( v_{i_1}\wedge\dots\wedge v_{i_k})_{1\le i_1<\dots<i_k\le d}$ shows that $\zeta(x)$ and $\zeta(y)$ are transverse. 

Next we show that $\zeta$ is dynamics preserving. Fix a sequence $\{\gamma_n\}$ in $\Gamma$ such that $\gamma_n \rightarrow x$ and $\gamma_n^{-1} \rightarrow y$. Since $\rho$ is $P_\Delta$-transverse, Lemma~\ref{lem: KAK} implies that $\alpha_k(\rho(\gamma_n)) \rightarrow \infty$ for all $k \in \Delta$, 
$$
U_\Delta(\rho(\gamma_n)) \rightarrow \xi(x) \quad \text{and} \quad U_\Delta(\rho(\gamma_n^{-1})) \rightarrow \xi(y).
$$
If we identify $\bigwedge^k\Rb^d$ with $\Rb^D$ using the basis 
\[
(e_{i_1}\wedge\dots\wedge e_{i_k})_{1\le i_1<\dots<i_k\le d}
\] 
then one can check that $\alpha_i(E^k \circ \rho(\gamma_n)) \rightarrow \infty$ for $i \in \{1,2,D-1,D-2\}$, 
$$
U_{1,2,D-2,D-1}(E^k \circ \rho(\gamma_n))=W_k\left(U_\Delta(\rho(\gamma_n))\right) \rightarrow \zeta(x) \quad \text{and} \quad U_{1,2,D-2,D-1}(E^k \circ \rho(\gamma_n^{-1})) \rightarrow \zeta(y).
$$
Hence by Lemma~\ref{lem: KAK}, 
$$
E^k \circ\rho(\gamma_n) (F) \rightarrow \zeta(x)
$$
for all $F \in \hat{\Fc}$ transverse to $\zeta(y)$.

Finally we show that $E^k\circ\rho$ is $(1,1,2)$-hypertransverse. Fix $x_1,x_2,x_3 \in \Lambda(\Gamma)$ pairwise distinct. Then there is a basis $\Bc=(b_1,\dots,b_d)$ of $\Rb^d$ such that $b_i\in \xi(x_1)^i\cap \xi(x_3)^{d-i+1}$ for all $i=1,\dots,d$, and $\xi(x_2)=u(\xi(x_3))$ for some $u\in U^+_{>0}(\Bc)$. Observe that 
\begin{itemize}
\item $\zeta^1(x_1)=\Span_{\Rb}(b_1\wedge\dots\wedge b_k)$, 
\item $\zeta^{\ell-2}(x_3)=\Span_{\Rb}(b_{i_1}\wedge\dots\wedge b_{i_k}:\{i_1,\dots,i_k\}\neq \{1,\dots,k\},\{1,\dots,k-1,k+1\})$
\item $\zeta^1(x_2)=\Span_{\Rb}\left(\sum_{1\le i_1<\dots<i_k\le d}u^{i_1,\dots,i_k}_{d-k+1,\dots,d}\,b_{i_1}\wedge\dots\wedge b_{i_k}\right),$ where $u^{i_1,\dots,i_k}_{j_1,\dots,j_k}$ denotes the minor given by the rows $1\le i_1<\dots<i_k\le d$ and the columns $1\le j_1<\dots<j_k\le d$ of the $d\times d$ matrix representing $u$ in the basis $\Bc$. 
\end{itemize}
In particular, $u^{1,\dots,k-1,k+1}_{d-k+1,\dots,d}\neq 0$, so
\[\zeta^1(x_1)+\zeta^1(x_2)+\zeta^{D-2}(x_3)\]
is a direct sum.
\end{proof}

We may now assemble the  first half of the proof of Corollary \ref{PSgen} which we restate here (the other half appearing in Theorem~\ref{thm:main in body} below) . 

\begin{corollary}\label{PSgen in body}
 If $\Gamma \subset \PSL(2,\Rb)$ is discrete and $\rho:\Gamma\to\PSL(d,\mathbb R)$ is a Hitchin representation, then
$$\delta^{\alpha_k}(\rho)\le 1$$
for all $k\in\{1,\ldots,d-1\}$.
Moreover, if $\Gamma$ is a lattice, then $\delta^{\alpha_k}(\rho)=1$ for all $k\in\{1,\ldots,d-1\}$.
\end{corollary} 

\begin{proof}
By Proposition \ref{prop: Hitchin transverse}, $E^k\circ\rho$ is $(1,1,2)$-hypertransverse. Theorem \ref{thm:BishopJones} implies that 
$$\mathrm{dim}_{H}\left(\Lambda_{1,c}(E^k\circ\rho(\Gamma))\right)=\delta^{\alpha_1}\left(E^k\circ\rho\right)=\delta^{\alpha_k}(\rho).$$
Let $W_k:\Fc\to \hat{\Fc}$ be as in the proof of Proposition~\ref{prop: Hitchin transverse}. Since $W_k$ is smooth and $\Lambda_{1,c}(E^k\circ\rho(\Gamma))=W_k^1(\Lambda_{k,c}(\rho(\Gamma)))$,
Proposition \ref{hitchinlimitset} implies that 
$$\mathrm{dim}_{H}\left(\Lambda_{1,c}(E^k\circ\rho(\Gamma))\right)\le 1.$$
Moreover, if $\Gamma$ is a lattice, then 
$\Lambda_1(E^k\circ\rho(\Gamma))$ is a curve and  $\Lambda_1(E^k(\rho(\Gamma)))\setminus\Lambda_{1,c}(E^k(\rho(\Gamma))$ is countable, so  
$\mathrm{dim}_{H}\left(\Lambda_{1,c}(E^k\circ\rho(\Gamma))\right)=1$.
\end{proof}

We next show that  $\mathrm{dim}_{H}\left(\Lambda_1(E^k\circ\rho(\Gamma))\right)<1$ if $\Gamma$ is geometrically finite but not a lattice. We make use of
the following generalization of a result of Labourie--McShane \cite{labourie-mcshane} which will be established in the appendix.

\begin{proposition}\label{prop: double hitchin}
Suppose that $\Gamma\subset\PSL(2,\Rb)$ is geometrically finite and torsion-free. If $\rho:\Gamma\to\PSL(d,\Rb)$ is a Hitchin representation, then
there exists a torsion-free lattice $\Gamma^D \subset \PSL(2,\Rb)$ so that $\Gamma\subset\Gamma^D$ and a Hitchin representation $\rho^D:\Gamma^D\to\PSL(d,\mathbb R)$ so
that $\rho=\rho^D|_\Gamma$.
\end{proposition}

We then use an argument due to Furusawa \cite{furusawa} to show that the critical exponent drops when one passes to an infinite index geometrically finite subgroup.
Corollary \ref{PSgen2} then follows from Propositions \ref{prop: double hitchin} and \ref{prop: entropy drop} and Corollary \ref{PSgen}.

\begin{proposition}\label{prop: entropy drop}
Suppose that  $\Gamma\subset\PSL(2,\Rb)$ is discrete and $G\subset\Gamma$ is a geometrically finite, infinite index subgroup. 
If $\rho:\Gamma\to\PGL(d,\Rb)$ is a $P_k$-transverse representation, then $\delta^{\alpha_k}(\rho|_{G})<\delta^{\alpha_k}(\rho)$. 
\end{proposition}

\begin{remark} After this paper was submitted for publication, a more general version of this proposition was established, see~\cite[Thm.\ 4.1]{CZZ3}
\end{remark} 

\begin{proof}
Since $G\subset\Gamma$ is geometrically finite and infinite index, $\Lambda(G)$ is a proper closed subset of $\Lambda(\Gamma)$. 
(If $C(G)$ is the convex hull of $\Lambda(G)$ in $\mathbb H^2$,
then $G \backslash C(G)$ has finite area. But if $\Lambda(G)=\Lambda(\Gamma)$, then $\Gamma$ preserves $C(G)$ which contradicts the fact that $\Gamma$ is discrete and $G$
has infinite index.)

We may assume, without loss of generality, that $G$ is torsion-free.
Let $D$ be a finite-sided convex fundamental polygon for the action of $G$ on $\mathbb H^2$ and let $I$ be a component of the intersection
of the closure of $D$ with $\partial\mathbb H^2$ whose interior has a non-empty intersection with $\Lambda(\Gamma)$. Since fixed point pairs of hyperbolic elements are dense in $\Lambda(\Gamma)\times\Lambda(\Gamma)$, we may find
a hyperbolic element $\gamma\in\Gamma$, both of whose fixed points lie in the interior of $I$. By passing to a power of $\gamma$, if necessary, we may assume
that there exists disjoint closed geodesic half-planes $H^+$ and $H^-$ centered at $\gamma^+$ and $\gamma^-$ 
and contained in the interior of $D$ so that $\gamma(H^-)=\mathbb H^2-\mathrm{int}(H^+)$, see Figure \ref{fig: double}. Klein's combination theorem then implies that 
$\hat D=D-\mathrm{int}(H^+\cup H^-)$ is a fundamental domain for $\hat\Gamma=\langle G,\gamma\rangle= G*\langle\gamma\rangle$.

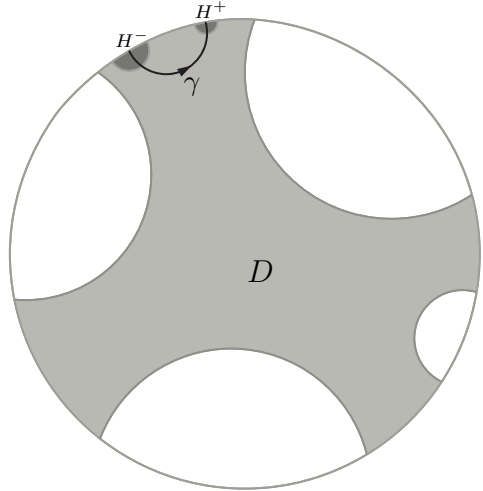
\begin{figure}[h]
    \centering
    \includegraphics[width=0.4\textwidth]{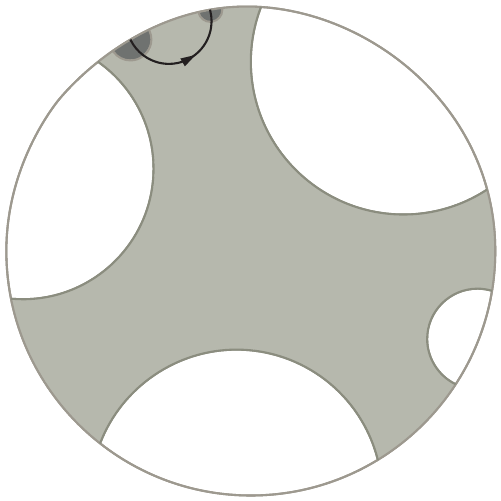}
    \put (-114,153){$\gamma$}
    \large
\put (-90,80){$D$}
\tiny
\put (-110,180){$H^+$}
\put (-140,169){$H^-$}
    \caption{$D$, $H^-$, and $H^+$.}
    \label{fig: double}
\end{figure}

Fix $b_0\in \hat D$. Let $\iota_{b_0} : \Hb^2 -\{b_0\} \rightarrow \partial \Hb^2$ be the map so that $z$ is contained in the geodesic ray starting at $b_0$ and limiting to $\iota_{b_0}(z)$. Choose $\epsilon>0$ small enough so that 
\begin{align}\label{eqn: projection bound two}
\d_{\partial\Hb^2}(\iota_{b_0}(H^+\cup H^-),\Lambda(G))\ge 2\epsilon.
\end{align}
There is a finite subset $F\subset G$ such that 
\begin{align}\label{eqn: projection bound one}
\d_{\partial\Hb^2}(\iota_{b_0}(\eta(p)),\Lambda(G))<\epsilon
\end{align}
for all $p\in D$ and all $\eta\in G-F$. Hence, by the residual finiteness of $G$, we may replace $G$ with a finite index subgroup to assume that 
\eqref{eqn: projection bound one} holds for all $\eta\in G-\id$ and all $p\in D$.

Let $\Gamma(n)\subset\Gamma$ denote the set of elements of the form $\eta_1\gamma_1\dots \eta_n\gamma_n$, where each $\eta_i$ is a non-trivial element of $G$, 
and each $\gamma_i$ is a  non-trivial element of $\langle\gamma\rangle$. (We adopt the convention that $\Gamma(0)=\{\id\}$.) 

\begin{lemma}\label{lem: ping pong}
For all $\zeta\in\bigcup_{n>0}\Gamma(n)$, $\eta\in G-\id$, and $m\ne 0$, we have
\[\d_{\partial\Hb^2}(\iota_{b_0}(\zeta^{-1}(b_0),\iota_{b_0}(\eta(b_0))\ge \epsilon\,\text{ and }\,\d_{\partial\Hb^2}(\iota_{b_0}(\eta^{-1}\zeta^{-1}(b_0),\iota_{b_0}(\gamma^m(b_0))\ge \epsilon.\]
\end{lemma}

\begin{proof}
Since $b_0\in\hat D$ ,$\gamma^m(b_0)\in H^+\cup H^-$ if $m\ne 0$, so \eqref{eqn: projection bound two} gives
\[\d_{\partial\Hb^2}(\iota_{b_0}(\gamma^m(b_0)),\Lambda(G))\ge 2\epsilon.\]
Moreover, by \eqref{eqn: projection bound one}, 
\[\d_{\partial\Hb^2}(\iota_{b_0}(\eta(b_0)),\Lambda(G))<\epsilon\]
for all $\eta\in G-\id$.

We may  apply a ping-pong argument to observe that $\zeta^{-1}(b_0)\in H^+\cup H^-$ for all $\zeta\in\bigcup_{n>0}\Gamma(n)$, so 
\[\d_{\partial\Hb^2}(\iota_{b_0}(\zeta^{-1}(b_0)),\Lambda(G))\ge 2\epsilon.\]
Then since $\zeta^{-1}(b_0)\in D$, \eqref{eqn: projection bound one} implies that
\[\d_{\partial\Hb^2}(\iota_{b_0}(\eta^{-1}\zeta^{-1}(b_0)),\Lambda(G))<\epsilon\]
for all $\eta\in G-\id$. The lemma follows easily from these four inequalities.
\end{proof}

For any subset $S\subset\Gamma$, denote
\[Q^{\alpha_k}_S(s)=\sum_{\gamma\in S}e^{-s\alpha_k(\kappa(\rho(\gamma)))}.\]

\begin{lemma} \label{lem: geom series}
There is some $C>0$ such that for all integers $n\ge 0$, 
\[Q^{\alpha_k}_{\Gamma(n)}(s)\ge \left(e^{-2Cs}Q^{\alpha_k}_{G-\id}(s)Q^{\alpha_k}_{\langle\gamma\rangle-\id}(s)\right)^n\]
\end{lemma}

\begin{proof}
We prove this by induction. For the base case $n=0$ is trivial. For the inductive step, suppose that $n>1$ and note that
\[Q^{\alpha_k}_{\Gamma(n)}(s)=\sum_{(\zeta,\eta,\gamma^m)\in\Gamma(n-1)\times(G-\id)\times(\langle \gamma\rangle-id)}e^{-s\alpha_k(\kappa(\rho(\zeta\eta\gamma^m)))}\]
By Lemma \ref{lem: ping pong} and Lemma \ref{lem:upper_bd_on_sk}, there is some $C>0$ such that
\begin{align*}
\alpha_k (\kappa(\rho(\zeta\eta\gamma^m)))&\le \alpha_k(\kappa(\rho(\zeta\eta)))+\alpha_k(\kappa(\rho(\gamma^m)))+C\\
&\le \alpha_k(\kappa(\rho(\zeta)))+\alpha_k(\kappa(\rho(\eta)))+\alpha_k(\kappa(\rho(\gamma^m)))+2C
\end{align*}
for all $\zeta\in\bigcup_{n>0}\Gamma(n)$, $\eta\in G-\id$, and $m\ne 0$. Thus, 
\begin{align*}
Q^{\alpha_k}_{\Gamma(n)}(s)&\ge \sum_{\zeta\in\Gamma(n-1)}\sum_{\eta\in G-\id}\sum_{m\ne 0}e^{-2Cs}e^{-s\alpha_k(\kappa(\rho(\zeta)))}e^{-s\alpha_k(\kappa(\rho(\eta)))}e^{-s\alpha_k(\kappa(\rho(\gamma^m)))}\\
&= e^{-2Cs}Q^{\alpha_k}_{\Gamma(n-1)}(s)Q^{\alpha_k}_{G-\id}(s)Q^{\alpha_k}_{\langle\gamma\rangle-\id}(s)\\
&\ge \left(e^{-2Cs}Q^{\alpha_k}_{G-\id}(s)Q^{\alpha_k}_{\langle\gamma\rangle-\id}(s)\right)^n, 
\end{align*}
where the last inequality is a consequence of the inductive hypothesis.
\end{proof}

Proposition \ref{crit exp divergence} implies that  $Q^{\alpha_k}_{G-\id}(\delta^{\alpha_k}(\rho|_G))=\infty$. Thus, there is some $s_0>\delta^{\alpha_k}(\rho|_G)$ such that 
\[e^{-2Cs_0}Q_{G-\id}^{\alpha_k}(s_0)Q_{\langle\gamma\rangle-\id}^{\alpha_k}(s_0)\ge 1.\]
By Lemma \ref{lem: geom series}, 
\[\sum_{n=0}^\infty\left(e^{-2Cs}Q^{\alpha_k}_{G-\id}(s)Q^{\alpha_k}_{\langle\gamma\rangle-\id}(s)\right)^n\le \sum_{n=0}^\infty Q^{\alpha_k}_{\Gamma(n)}(s)\le Q_{\rho}^{\alpha_k}(s)\]
for all $s\ge 0$. Hence, $Q_{\rho}^{\alpha_k}(s_0)=\infty$, which implies that
\[\delta^{\alpha_k}(\rho)\ge s_0>\delta^{\alpha_k}(\rho|_G).\qedhere\]
\end{proof}

\subsection{Proof of entropy rigidity}
We are now ready to prove our main result, Theorem~\ref{thm:main}, which we restate here for reference.

\begin{theorem}\label{thm:main in body}
If $\Gamma \subset \PSL(2,\Rb)$ is geometrically finite, $\rho : \Gamma \rightarrow \PSL(d, \Rb)$ is Hitchin and
$\phi = \sum c_j \alpha_j \in \left(\mathfrak{a}^*\right)^+$, then 
\begin{align*}
h^{\phi}(\rho) \le \frac{1}{c_1+\dots + c_{d-1}}.
\end{align*}
Moreover, equality occurs if and only if $\Gamma$ is a lattice and either 
\begin{enumerate}
\item
$\phi=c_k\alpha_k$ for some $k$.
\item
$\rho(\Gamma)$ lies in an irreducible image of $\mathsf{PSL}(2,\mathbb R)$.
\item
$d=2n-1$, the Zariski closure of $\rho(\Gamma)$ is conjugate to $\mathsf{PSO}(n,n-1)$ and  $\phi=c_k\alpha_k+c_{d-k}\alpha_{d-k}$ for some $k$.
\item
$d=2n$, the Zariski closure of $\rho(\Gamma)$ is conjugate to $\mathsf{PSp}(2n,\mathbb R)$ and $\phi=c_k\alpha_k+c_{d-k}\alpha_{d-k}$ for some $k$.
\item
$d=7$, the Zariski closure of $\rho(\Gamma)$ is conjugate to $\mathsf{G}_2$ and $\phi=c_1\alpha_1+c_3\alpha_3+c_4\alpha_4+c_6\alpha_6$, or $\phi=c_2\alpha_2+c_5\alpha_5$.
\end{enumerate}
\end{theorem} 

\begin{proof} 
Corollary \ref{PSgen} and \ref{PSgen2}, and Proposition \ref{entropy and critical exponent} imply that for all $k\in\Delta$,
$$h^{\alpha_k}(\rho)=\delta^{\alpha_k}(\rho) \le 1$$
with equality if and only if $\Gamma$ is a lattice.

Since $h^{a\phi}(\rho)=\frac{1}{a}h^\phi(\rho)$ for all $\phi\in\mathcal B^+(\rho)$ and $a>0$, it follows from Theorem \ref{quint indicator} part (3) that, if $\phi=\sum c_j\alpha_j\ne 0$ and $c_j\ge 0$ for all $j$, then
$$
h^{\phi}(\rho) \le \frac{1}{\frac{c_1}{h^{\alpha_1}(\rho)}+\dots + \frac{c_{d-1}}{h^{\alpha_{d-1}}(\rho)}} \le \frac{1}{c_1+\dots + c_{d-1}}.
$$
In particular, if $\Gamma$ is not a lattice, then
$$h^{\phi}(\rho) < \frac{1}{c_1+\dots + c_{d-1}}.$$

Now, suppose that $\Gamma$ is a lattice. Then
$$h^{\phi}(\rho) \le \frac{1}{c_1+\dots + c_{d-1}}.$$
Let $T_\phi$ denote the simplex in $\mathfrak a^*$ whose vertices are the set $\{\alpha_j:c_j\neq 0\}$. Theorem \ref{quint indicator} part (3) implies that if equality holds then $T_\phi\subset \mathcal Q_\Delta(\rho)$, and so by Theorem \ref{quint indicator} part (4),
$\alpha_i(\nu(\rho(\gamma)))=\alpha_j(\nu(\rho(\gamma)))$  if $c_i\ne 0$ and $c_j\ne 0$.

We recall that Sambarino \cite{sambarino-positive} classified the possible Zariski closures of images of Hitchin representations.

\begin{theorem}\cite[Thm.\ A]{sambarino-positive}\label{thm: Sambarino}
Suppose that $\Gamma\subset\PSL(2,\Rb)$ is a lattice, and $\rho:\Gamma\to\PSL(d,\Rb)$ is a Hitchin representation. Then the Zariski closure of $\rho(\Gamma)$ either  lies in an irreducible image of $\mathsf{PSL}(2,\mathbb R)$ or is
conjugate to either $\PSL(d,\Rb)$, $\mathsf{PSp}(2n,\Rb)$ when $d=2n$, $\mathsf{PSO}(n,n-1)$ when $d=2n-1$, or $\mathsf{G}_2$ when $n=7$.
\end{theorem}

Since the Zariski closure of $\rho(\Gamma)$ is semisimple, a result of Benoist  \cite{benoist-limit-cone} implies that $\mathcal B(\rho)$ is a convex open subset of the limit
cone of its Zariski closure with non-empty interior. Therefore, if $g$ lies in the Zariski closure of $\rho(\Gamma)$, then $\alpha_i(\nu(g))=\alpha_j(\nu(g))$,  if $c_i\ne 0$ and $c_j\ne 0$.

Our result then follows from the following observation about subgroups of $\mathsf{PSL}(d,\mathbb R)$ (see, for example the discussion in \cite[Appendix B]{Zimmer}).

\begin{proposition} Suppose $\mathsf{G} \subset \PSL(d,\Rb)$ is a group and 
$$
S(\mathsf{G})=\{ (i,j) : i < j \text{ and } \alpha_i(\nu(g)) = \alpha_j(\nu(g))  \text{ for all } g \in G\}. 
$$
\begin{enumerate}
\item If $G = \PSL(d,\Rb)$, then $S(\mathsf{G}) = \emptyset$. 
\item If $G$ lies in an irreducible image of $\mathsf{PSL}(2,\mathbb R)$, then $S(\mathsf{G})=\{(i,j) : i < j\}$. 
\item If $d=2n$ and $G=\mathsf{PSp}(2n,\Rb)$, then $S(\mathsf{G})= \{ (k,d-k) : 1 \le k \le n-1\}$. 
\item If $d=2n-1$ and $G=\mathsf{PSO}(n,n-1)$, then $S(\mathsf{G})= \{ (k,d-k) : 1 \le k \le n-1\}$. 
\item If $d=7$ and $G=\mathsf{G}_2$, then $S(\mathsf{G})= \{ (i,j) : i < j \text{ and } i,j \in \{1,3,4,6\} \}\cup\{(2,5)\}$. 
\end{enumerate} 
\end{proposition} 

\end{proof} 


\section{Proof of Corollary~\ref{cor:symmetric rigidity}}


In this section we prove Corollary~\ref{cor:symmetric rigidity} which we restate here. 

\begin{corollary}\label{cor:symmetric rigidity in body} If $\Gamma \subset \PSL(2,\Rb)$ is geometrically finite and $\rho : \Gamma \rightarrow \PSL(d, \Rb)$ is Hitchin, then 
\begin{align*}
\delta_X(\rho) \le 1.
\end{align*}
Moreover, $\delta_X(\rho) =1$ if and only if $\Gamma$ is a lattice and $\rho(\Gamma)$ lies in the image of an irreducible representation $\PSL(2,\Rb)\to\PSL(d,\Rb)$.
\end{corollary}

\begin{proof} Let $x_0 = [ \mathsf{PSO}(d)]\in X$ and
$$
C = (d-1)^2 + (d-3)^2 + \dots + (3-d)^2 + (1-d)^2.
$$
Then 
$$
\d_X( g(x_0), x_0) = \frac{2}{\sqrt{C}} \sqrt{ \sum_{j=1}^d (\log \sigma_j(g))^2 } = \frac{2}{\sqrt{C}} \norm{ \kappa(g)}_2
$$
for all $g \in \PSL(d,\Rb)$.

Let $e_1^*, \dots, e_d^* \in \mathfrak{a}^*$ be the elements with
$$
e_j^*(\mathrm{diag}(A_1,\dots, A_d)) = A_j.
$$
Then define 
\begin{align*}
\phi &= \frac{2}{C}\Big( (d-1)e_1^* + (d-3)e_2^* + \dots + (3-d)e_{d-1}^* + (1-d)e_d^* \Big) \\
& = \frac{2}{C} \left( \sum_{j=1}^{\lfloor\frac{d-1}{2}\rfloor}(d+1-2j)(\alpha_j + \dots + \alpha_{d-j}) + \left(d-2\left\lfloor \frac{d-1}{2} \right\rfloor-1\right) \alpha_{\lfloor d/2 \rfloor} \right).
\end{align*}
Observe that the sum of all the coefficients of 
\[\sum_{j=1}^{\lfloor\frac{d-1}{2}\rfloor}(d+1-2j)(\alpha_j + \dots + \alpha_{d-j}) + \left(d-2\left\lfloor \frac{d-1}{2} \right\rfloor-1\right) \alpha_{\lfloor d/2 \rfloor}\quad\text{ is }\quad\frac{C}{2}.\]
Thus, by Theorem~\ref{thm:main} and Proposition \ref{critandent}, we have 
$$
\delta^{\phi}(\rho) \le 1.
$$ 
Further, 
$$
\phi(\kappa(g)) \le \norm{\phi}_2 \norm{\kappa(g)}_2 =  \frac{2}{\sqrt{C}}\norm{ \kappa(g)}_2 = \d_X( g(x_0), x_0)
$$
for all $g \in \PSL(d,\Rb)$. Hence 
$$
\delta_X(\rho) \le \delta^{\phi}(\rho) \le 1.
$$

If $\delta_X(\rho) = 1$, then $\delta^\phi(\rho) = 1$ and hence Theorem~\ref{thm:main} implies that $\Gamma$ is a lattice and $\rho(\Gamma)$ lies in the image of an irreducible representation $\PSL(2,\Rb)\to\PSL(d,\Rb)$. 

For the other direction, suppose that $\Gamma$ is a lattice and $\rho(\Gamma)$ lies in the image of an irreducible representation $\tau: \PSL(2,\Rb)\to\PSL(d,\Rb)$. Let $\Gamma_1 : = \tau^{-1}(\rho(\Gamma))$. Since $\tau$ induces an isometry $\Hb^2 \hookrightarrow X$, 
$$
\delta_X(\rho) = \delta_{\Hb^2}(\Gamma_1).
$$
Thus it suffices to prove that $\Gamma_1$ is a lattice in $\PSL(2,\Rb)$. 

Since $\Lambda(\Gamma) = \partial \Hb^2$, each element of $\Gamma$ acts non-trivially on $\Lambda(\Gamma)$. Then by the equivariance of the Anosov boundary map, $\rho$ must be faithful. So 
$$
\gamma \in \Gamma \mapsto \tau^{-1}(\rho(\gamma)) \in  \Gamma_1
$$
is an isomorphism. Further, by~\cite[Thm.\ 1.4]{CZZ}, $\alpha \in \Gamma$ is parabolic if and only if $\tau^{-1}(\rho(\alpha))$ is parabolic. So ~\cite[Thm.\ 3.3]{tukia-isomorphism} implies that $\Lambda(\Gamma_1)$ is homeomorphic to $\Lambda(\Gamma) = \partial \Hb^2$. So $\Lambda(\Gamma_1) = \partial \Hb^2$. Hence $\Gamma_1$ is a lattice and 
\[
\delta_X(\rho) = \delta_{\Hb^2}(\Gamma_1) = 1.
\qedhere\]
\end{proof}

\appendix


\section{Doubling Hitchin representations}


In this appendix, we show that if $\Gamma$ is a geometrically finite Fuchsian group uniformizing a finite area hyperbolic 2-orbifold with
non-empty totally geodesic boundary $Q$ and $\rho:\Gamma\to\mathsf{PSL}(d,\mathbb R)$ is
a Hitchin representation, then we may extend $\rho$ to a Hitchin representation of $\Gamma_D$ where $\Gamma_D$ uniformizes the double of $Q$ along its boundary
and $\Gamma\subset\Gamma_D$. Proposition \ref{prop: double hitchin} is an immediate consequence of this result.

When $\Gamma$ is geometrically finite but not a lattice, $\Lambda(\Gamma)$ is a proper subset of  $\partial\Hb^2$. 
Let $\mathcal {C}(\Gamma) \subset \Hb^2$ denote the convex hull of $\Lambda(\Gamma)$, and let $\Sc(\Gamma)$ denote the set of boundary components in $\Hb^2$ of $\mathcal C(\Gamma)$. Note that every $b\in\Sc(\Gamma)$ is the axis of some hyperbolic $\beta_b\in\Gamma$. 
For each $b\in\Sc(\Gamma)$, let $r_b\in\PGL(2,\Rb)$ be the reflection about $b$. 
Let $\Gamma^D\subset\PSL(2,\Rb)$ denote the index two subgroup of the group $\hat\Gamma=\left\langle \Gamma, \{r_b\}_{b\in\Sc(\Gamma)}\right\rangle$
consisting of orientation-preserving isometries. Then $\hat\Gamma \backslash \mathbb H^2$ is the orbifold obtained from $\mathcal{C}(\Gamma)$ by regarding
the boundary components as mirrors and $\Gamma^D \backslash \mathbb H^2$ is its orientable double cover. In particular, $\Gamma^D$ is a lattice.
We refer to $\Gamma^D$ as the \emph{double} of $\Gamma$. 

Labourie and McShane \cite[Cor.\ 9.2.6]{labourie-mcshane} showed that any Hitchin representation of a torsion-free convex cocompact Fuchsian group
extends to a Hitchin representation of its double. We observe that their result generalizes to arbitrary geometrically finite Fuchsian groups. This is probably known to experts, but
we provide a proof for completeness.

\begin{proposition}\label{prop: double}
Let $\Gamma\subset\PSL(2,\Rb)$ be a geometrically finite Fuchsian group that is not a lattice, and let $\rho:\Gamma\to\PSL(d,\Rb)$ be a Hitchin representation.
There exists a Hitchin representation $\rho^D:\Gamma^D\to\mathsf{PSL}(d,\Rb)$ such that $\rho^D|_\Gamma=\rho$.
\end{proposition}

Choose a convex, finite sided fundamental polygon $D$ for the action of $\Gamma$ on $\mathcal{C}(\Gamma)$
so that $D$ intersects each $\Gamma$-orbit in $\Sc(\Gamma)$ exactly once. Let $\mathcal E(\Gamma)$ denote
the finite collection of geodesics in $\Sc(\Gamma)$ which intersect $D$. Notice that  if $\Gamma$ has finite presentation $\langle X\ :\ R\rangle$, then
$$\hat\Gamma=\left\langle X, \{r_b\}_{b\in\mathcal E(\Gamma)}\ :\ R,\ \ \{r_b^2\}_{b\in\mathcal E(\Gamma)}\right\rangle$$
is a finite presentation for $\hat\Gamma$.

Since $\rho$ is a Hitchin representation, there is a continuous, $\rho$-equivariant positive map $\xi:\Lambda(\Gamma)\to\Fc$. For each $b\in\Sc(\Gamma)$, let $R_{b}$ be the projectivization of the linear map that fixes both $\rho(\beta_{b}^+)$ and $\rho(\beta_{b}^-)$, and acts on the line 
$\xi^j(\beta_{b}^+)\cap\xi^{d-j+1}(\beta_{b}^-)$ by scaling by $(-1)^{j-1}$ for all $j\in\{1,\dots,d\}$. We define $\hat\rho:\hat\Gamma\to \mathsf{PGL}(d,\mathbb R)$ by
setting $\hat\rho(r_{b})=R_b$ for all $b\in\Ec(\Gamma)$. Notice that since $R_b^2=I$, $\hat \rho$ is a representation. We then let $\rho^D=\hat\rho|_{\Gamma^D}$.

Let $G\subset\PGL(2,\Rb)$ be the subgroup generated by $\{r_b:b\in\Sc(\Gamma)\}$. By a ping-pong type argument, we see that every element in $G$ can be written uniquely as a reduced word in the alphabet $\{r_b:b\in\mathcal S(\Gamma)\}$, and we deduce the following lemma.

\begin{lemma}\label{lem: S}
If $r_1,r_2\in G$, then one of the following holds:
\begin{itemize}
\item $r_1=r_2$,
\item $r_1=r_2r_b$ for some $b\in\mathcal S(\Gamma)$, in which case $r_1(\mathcal{C}(\Gamma))\cap r_2(\mathcal{C}(\Gamma))=r_1(b)=r_2(b)$,
\item $r_1(\mathcal{C}(\Gamma))\cap r_2(\mathcal{C}(\Gamma))$ is empty.
\end{itemize}
\end{lemma}

If $r\in G$, we may write $r$ uniquely as a reduced word $r_{b_1}\dots r_{b_l}$ for some $b_1,\dots,b_l\in\Sc(\Gamma)$, so
\[\hat\rho(r)=R_{b_1}\dots R_{b_l}\in\PGL(d,\Rb).\]
Then by Lemma \ref{lem: S}, if $r_1,r_2\in G$ are reflections such that the intersection 
\[r_1(\Lambda(\Gamma))\cap r_2(\Lambda(\Gamma))\] 
is non-empty, then either $r_1=r_2$, or $r_1=r_2r_b$ for some $b\in\Sc(\Gamma)$, in which case this intersection is $r_1(b)=r_2(b)$. Hence, we have a well-defined map
\[\xi^D:\Lambda(\Gamma^D)=\bigcup_{r \in G}r(\Lambda(\Gamma))\to\Fc_d\]
given by $\xi^D(r(x))=\hat\rho(r)(\xi(x))$ for all $r\in G$ and $x\in\Lambda(\Gamma)$. It is straightforward to verify that $\xi^D$ is $\rho^D$-equivariant. 
It now suffices to show that $\xi^D$ is positive, since by \cite[Thm.\ 9.2]{CZZ}, a representation of a geometrically finite Fuchsian group with a positive
equivariant limit map is Hitchin.

\begin{lemma}
\label{double positive}
$\xi^D$ is a positive map. 
\end{lemma}

\begin{proof}
It is sufficient to show that $\xi^D$ restricted to $\bigcup_{r\in S'}r(\Lambda(\Gamma))$ is a positive map for any finite subset $S'\subset G$ such that $\bigcup_{r\in S'}r(\mathcal{C}(\Gamma))$ is connected. We prove this by induction on the size of $S'$. The base case where $\# S'=1$ is trivial because positivity is invariant under projective transformations.

For the inductive step, let $S''\subset S'$ be a subset such that $\#S''=\#S'-1$ and $\bigcup_{r\in S''}r(\mathcal{C}(\Gamma))$ is connected. Let $r_1\in S'-S''$. 
Since $\bigcup_{r\in S'}r(\mathcal{C}(\Gamma))$ is connected, there is some $r_2\in S''$ such that $r_1(\mathcal{C}(\Gamma))$ and 
$r_2(\mathcal{C}(\Gamma))$ share a common boundary component $d$, which lies in the $\Gamma$-orbit of  a geodesic in $\mathcal E(\Gamma)$.
Let $y_1,y_2\in\partial\Hb^2$ be the endpoints of $d$, with notation chosen so that $y_1<z_1<y_2<z_2$ for some (any) $z_1\in r_1(\Lambda(\Gamma))$ and $z_2\in r_2(\Lambda(\Gamma))$. 

Pick any $\ell$-tuple of points $x_1<\dots<x_\ell<x_1$ in $\bigcup_{r\in S'}r(\Lambda(\Gamma))$ with $\ell\ge 3$. 
By adding points and cyclically permuting, we  may assume that there is some $t\in\{2,\dots,\ell\}$ such that $y_1=x_1$, $y_2=x_t$, and
$$\{x_i:1\le i\le t\}\subset r_1(\Lambda(\Gamma)), \quad
\{x_i:t\le i\le \ell\}\in \bigcup_{r\in S''}r(\Lambda(\Gamma)) \quad \text{and} \quad x_\ell\in r_2(\Lambda(\Gamma)).$$
By the inductive hypothesis, the restrictions of  $\xi_D$ to $\bigcup_{r\in S''}r(\Lambda(\Gamma))$ and $r_1(\Lambda(\Gamma))$ are positive maps. 
Thus, 
\begin{equation}
\label{enough positivity}
\left(\xi^D(x_1),\dots,\xi^D(x_t)\right)\,\text{ and }\,\left(\xi^D(x_t),\dots,\xi^D(x_\ell),\xi^D(x_1)\right)
\end{equation}
are positive tuples of flags. In the cases when $t=2$ or $t=\ell$, it is immediate that $(\xi(x_1),\dots,\xi(x_\ell))$ is positive.

Suppose now that $t\in\{3,\dots,\ell-1\}$. Let $r_d$ denote the reflection in $d$ and notice that $r_d\in\hat\Gamma$ and $r_d(x_\ell)\in r_1(\Lambda(\Gamma))$. Since $\xi^D$ is positive on $r_1(\Lambda(\Gamma))$, there is an ordered basis $\Bc$ associated to $(\xi^D(x_1),\xi^D(x_t))$ and $u,w\in U_{>0}(\Bc)$ so that $\xi^D(x_2)=u\xi^D(x_t)$ and $R_d(\xi^D(x_\ell))=\xi^D(r_d(x_\ell))=w\xi^D(x_t)=wR_d\xi^D(x_t)$. By the cofactor formula for computing minors, observe that $(R_dw R_d)^{-1}=R_dw^{-1} R_d\in U_{>0}(\Bc)$. Thus, 
$$\left(\xi^D(x_1),\xi^D(x_2),\xi^D(x_t),\xi^D(x_\ell)\right)=\left(\xi^D(x_1),u\xi^D(x_t),\xi^D(x_t),R_dw R_d\xi^D(x_t)\right)$$
is a positive quadruple of flags (see Lemma 7.5 of \cite{CZZ}).

Combining this  with (\ref{enough positivity}), Lemma \ref{lem: basic positivity}(5) implies that 
\[(\xi^D(x_1),\xi^D(x_2),\dots,\xi^D(x_\ell))\]
is positive. This completes the proof of Lemma \ref{double positive} and thus the proof of Proposition \ref{prop: double}.
\end{proof}



\begin{thebibliography}{100}
\begin{scriptsize}

\bibitem{AMS} H. Abels, G. Margulis and G. Soifer, ``Semigroups containing proximal linear maps,'' 
{\em Israel J. Math.} {\bf 9}(1995), 1--30.

\bibitem{beardon-maskit} A. Beardon and B. Maskit, ``Limit points of Kleinian groups and finite sided fundamental polyhedra,''
{\em Acta Math.} {\bf 132}(1974), 1--12.

\bibitem{benoist-limit-cone} Y. Benoist, ``Propri\'et\'es asymptotiques des groupes lin\'eaires,'' {\em G.A.F.A.}, {\bf 7} (1997), 1--47.

\bibitem{benoist-quint-book} Y. Benoist and J.F. Quint, {\em Random Walks on Reductive Groups}, Springer-Verlag, 2016.

\bibitem{bishop-conical} C. Bishop, ``On a theorem of Beardon and Maskit,''
{\em Ann. Acad. Sci. Fenn. Math.} {\bf 21}(1996),  383--388.

\bibitem{bishop-jones} C. Bishop and P. Jones, ``Hausdorff dimension and Kleinian groups,''  {\em Acta Math.} {\bf 179}(1997),  1---39.

\bibitem{bishop-peres} C. Bishop and Y. Peres, {\em Fractals in Probability and Analysis}, Cambridge University Press, 2016. 


\bibitem{BCZZ} P.-L. Blayac, R. Canary, F. Zhu, and A. Zimmer, ``Counting, mixing and equidistribution for GPS systems with applications to relatively Anosov groups,'' preprint, arXiv:2404.09718.

\bibitem{BPS} J. Bochi, R. Potrie and A. Sambarino, ``Anosov representations and dominated splittings,'' {\em J.E.M.S.} {\bf 21}(2019), 3343--3414.

\bibitem{bonahon-ends} F. Bonahon, ``Bouts des vari\'et\'es hyperboliques de dimension 3,''
\emph{Annals of Math.} \textbf{124} (1986), 71--158.


\bibitem{bowditch-relhyp} B. Bowditch, ``Relatively hyperbolic groups,'' {\em Int. J. Alg. Comp.} {\bf 22}(2012),  66 pp. 

\bibitem{BCKM} H. Bray, R. Canary, L.-Y. Kao, and G. Martone, ``Counting, equidistribution and entropy gaps at infinity with applications to cusped Hitchin representations,'' {\em J. Reine Angew. Math.} {\bf 791}(2022), 1--51.
                    

\bibitem{BCKM2} H. Bray, R. Canary, L.-Y. Kao, and G. Martone, ``Pressure metrics for cusped Hitchin components,'' {\em Adv. Math.} {\bf 435}(2023), 109352.


\bibitem{BCLS} M. Bridgeman, R. Canary, F. Labourie and A. Sambarino, ``The pressure metric for 
Anosov representations,'' {\em G.A.F.A.}, {\bf 25}(2015), 1089--1179.


\bibitem{bowen-series} R. Bowen and C. Series, ``Markov maps associated with Fuchsian groups," {\em Publ. Math.,
Inst. Hautes Etud. Sci.} {\bf 50} (1979), 153--170.

\bibitem{brooks-covers} R. Brooks, ``The bottom of the spectrum of a Riemannian covering,'' {\em J. Reine Angew. Math.} 
{\bf 357}(1985), 101--114.

\bibitem{burger} M. Burger, ``Intersection, the Manhattan curve and Patterson-Sullivan theory in rank 2," {\em Internat. Math. Res. Notices} {\bf 7}(1993), 217--225.

\bibitem{canary-laplace} R. Canary, ``On the Laplacian and the geometry of hyperbolic 3-manifolds,'' {\em J. Diff. Geom} {\bf 36}(1992), 349--367.


\bibitem{CLS} R. Canary, M. Lee, A. Sambarino and M. Stover, ``Amalgam Anosov
representations,'' {\em Geom. Top.} {\bf 21}(2017), 215--251.


\bibitem{CZZ} R. Canary, T. Zhang and A. Zimmer, ``Cusped Hitchin representations and Anosov representations of geometrically finite Fuchsian groups,'' {\em Adv. Math.} {\bf 404}(2022), 108439.

\bibitem{CZZ3} R. Canary, T. Zhang and A. Zimmer, ``Patterson-Sullivan measures for transverse subgroups,'' {\em J. Mod. Dyn.} {\bf 20}(2024), 319--377.

\bibitem{CZZ4} R. Canary, T. Zhang and A. Zimmer, ``Patterson-Sullivan measures for relatively Anosov groups,'' {\em Math. Ann.} {\bf 392} (2025), 2309--2363.

\bibitem{CLT} D. Constantine, J.-F. Lafont and D. Thompson, ``Strong symbolic dynamics for geodesic flow
on $\mathrm{CAT}(-1)$ spaces and other metric Anosov flows,''
{\em J. \'Ecole Poly. Math.} {\bf 7}(2020), 201--231.


\bibitem{AIF} M. Crampon and L. Marquis, ``Finitude g\'eom\'etrique en g\'eom\'etrie de Hilbert,''
{\em Ann. Inst. Four.}  {\bf 64}(2014), 2299--2377.

\bibitem{DOP} F. Dal'bo, J.-P. Otal, and M. Peign\'e, ``S\'eries de Poincar\'e des groupes g\'eom\'etriquement
finis,'' {\em Israel J. Math.} {\bf 118}(2000), 109--124.

\bibitem{DGK} J. Danciger, F. Gu\'eritaud and F. Kassel, ``Convex cocompact actions in real projective
geometry,'' preprint, arXiv:1704.08711, to appear in Ann. Sci. \'Ec. Norm. Sup\'er.

\bibitem{dalbo-peigne} F. Dal'bo and M. Peign\'e, ``Some negatively curved manifolds with cusps, mixing and counting,''
{\em J. Reign. Angew. Math.}. {\bf 497}(1998), 141--169.

\bibitem{fock-goncharov} V. Fock and A. Goncharov, ``Moduli spaces of local systems and higher Teichm\"uller theory,''
{\em Publ. Math. I.H.E.S.} {\bf 103}(2006), 1--211.

\bibitem{furusawa} H. Furusawa, ``The exponent of convergence of Poincar\'e series of combination groups,''
{\em Tohoku Math. J.} {\bf 43}(1991),  1--7.

\bibitem{GMT} O. Glorieux, D. Montclair and N. Tholozan, ``Hausdorff dimension of limit sets for projective Anosov representations,'' {\em J. \'Ec. polytech., Math.} {\bf 10}(2023), 1157--1193.

\bibitem{GGKW} F. Gu\'eritaud, O. Guichard, F. Kassel and A. Wienhard, ``Anosov representations and proper actions,''
{\em Geom. Top.} {\bf 21}(2017), 485--584.

\bibitem{guichard-wienhard} O. Guichard and A. Wienhard, ``Anosov representations: Domains of discontinuity and
applications,'' {\em Invent. Math.} {\bf 190}(2012), 357--438.

\bibitem{hitchin} N. Hitchin, ``Lie groups and Teichm\"uller space,''
{\em Topology} {\bf 31}(1992), 449--473. 

\bibitem{islam-zimmer-lms} M. Islam and A. Zimmer, ``A flat torus theorem for convex co-compact actions of projective linear groups,''
{\em Jour. L.M.S.} {\bf 103}(2021), 470--489.

\bibitem{kapovich-leeb} M. Kapovich and B. Leeb, ``Relativizing characterizations of Anosov subgroups, I,''
(with an appendix by Gregory A. Soifer). {\em Groups Geom. Dyn.} {\bf 17}(2023), 1005--1071.

\bibitem{KLP} M. Kapovich, B. Leeb and J. Porti, ``Anosov subgroups: Dynamical and geometric characterizations,''
{\em Eur. Math. J.} {\bf 3}(2017), 808--898.

\bibitem{KMO} D. Kim, Y. Minsky and H. Oh, ``Tent property and directional limit sets for self-joinings of hyperbolic manifolds,'' {\em Geom. Ded.} {\bf 218}(2024), Paper No. 14, 18 pp.

\bibitem{KTZ} S. Kim, S.P. Tan and T. Zhang, ``Weakly positive and directed Anosov representations,'' {\em Adv. Math.} {\bf 408}(2022), 108611.

\bibitem{labourie-invent} F. Labourie, ``Anosov flows, surface groups and curves in projective space,''
{\em Invent. Math.} {\bf 165}(2006), 51--114.

\bibitem{labourie-mcshane} F. Labourie and G. McShane, ``Cross ratios and identities for higher Thurston theory,''
{\em Duke Math. J.} {\bf 148}(2009),  279--345.

\bibitem{ledrappier-sarig} F. Ledrappier and O. Sarig, ``Fluctuations of ergodic sums for horocycle flows on $\mathbb Z^d$-covers of
finite volume surfaces,''  {\em Disc. Cont. Dyn. Sys.} {\bf 22}(2008), 247--325.

\bibitem{lusztig} G. Lusztig, ``Total positivity in reductive groups,'' in {\em Lie theory and geometry}, Birkh\"auser, 1994, 531--568. 

\bibitem{marquis-around} L. Marquis, ``Around groups in Hilbert geometry,'' in {\em Handbook of Hilbert
Geometry}, European Mathematical Society Publishing House, 2014, 207--261.

\bibitem{MU} D. Mauldin and M. Urbanski, {\em Graph directed Markov systems}, Cambridge Tracts in Mathematics, 
vol. 148, Cambridge University Press, Cambridge, Cambridge, 2003.


\bibitem{potrie-sambarino} R. Potrie and A. Sambarino, ``Eigenvalues and entropy of a Hitchin representation,''
{\em Invent. Math.}  {\bf 209},(2017),  885--925.

\bibitem{PSW1} M. Pozzetti, A. Sambarino and A. Wienhard, ``Conformality for a robust class of non-conformal attractors,''
{\em J. Reine Angew. Math} {\bf 74}(2021), 1--51.

\bibitem{quint-ps} J.F. Quint, ``Mesures de Patterson-Sullivan en rang sup\'erieur,'' 
{\em G.A.F.A} {\bf 12}(2002), 776--809.

\bibitem{sambarino-quantitative} A. Sambarino, ``Quantitative properties of convex representations,''
{\em Comm. Math. Helv.} {\bf 89}(2014), 443--488.

\bibitem{sambarino-hyperconvex} A. Sambarino, ``Hyperconvex representations and exponential growth,''
{\em Erg. Thy. Dyn. Sys.} {\bf 34}(2014), 986--1010.


\bibitem{sambarino-positive} A. Sambarino, ``Infinitesmal Zariski closures of positive representations,'' {\em J. Diff. Geom.} {\bf 128}(2024), 861--901.

\bibitem{sambarino-dichotomy} A.Sambarino, ``A short note on an ergodic dichotomy,'' {\em Erg. Thy. Dyn. Sys.} {\bf 44}(2024) 236--289.

\bibitem{sarig-2003} O. Sarig, ``Existence of Gibbs measures for countable Markov shifts,''
{\em Proc. A.M.S.} {\bf 131}(2003), 1751--1758. 

\bibitem{sarig-2009} O. Sarig, {\em Lecture notes on thermodynamic formalism for topological Markov shifts}, 2009.


\bibitem{stadlbauer} M. Stadlbauer, ``The return sequence of the Bowen-Series map
for punctured surfaces,'' {\em Fund. Math.} {\bf 182}(2004), 221--240.

\bibitem{sullivan-ergodic} D. Sullivan, ``The ergodic theory at infinity of an arbitrary discrete group of hyperbolic motions,'' in
{\em Riemann surfaces and related topics: proceedings of the 1978 Stony Brook conference}, I. Kra and B. Maski, ed., Princeton University
Press, 1981, 465--496.

\bibitem{sullivan-gf} D. Sullivan, ``Hausdorff measures old and new, and limit sets of geometrically finite Kleinian groups,'' {\em Acta Math.} {\bf 153}(1984),  259--277. 

\bibitem{SWZ} Z. Sun, A. Wienhard and T. Zhang, ``Flows on the $\mathsf{PGL}(V)$-Hitchin component,'' {\em G.A.F.A.} {\bf 30}(2020), 588--692.


\bibitem{tukia-isomorphism} P. Tukia, ``On isomorphisms of geometrically finite M\"obius groups,'' {\em Inst. Hautes \'Etudes Sci. Publ. Math.}  {\bf 61}(1985), 171--214. 

\bibitem{feng} F. Zhu, ``Relatively dominated representations,'' {\em Ann. Inst. Four.} {\bf 71}(2021), 2169--2235.


\bibitem{Zimmer} A. Zimmer, ``Projective Anosov representations, convex cocompact actions, and rigidity," {\em J. Diff. Geom.} {\bf 119}(2021), 513--586.


\end{scriptsize}
\end{thebibliography}
\end{document}